\documentclass[11pt]{article}

\usepackage[margin=1in]{geometry}
\usepackage[utf8]{inputenc}
\usepackage[T1]{fontenc}
\usepackage{amsmath,amssymb,amsthm,mathtools}
\usepackage{bm}
\usepackage{booktabs}
\usepackage{enumitem}
\usepackage{graphicx}
\usepackage{listings}
\usepackage{placeins}
\usepackage{wrapfig}
\usepackage[dvipsnames]{xcolor}
\usepackage{microtype}
\usepackage{cite}
\usepackage{titletoc}

\input{macros.tex}
\numberwithin{equation}{section}
\usepackage[boxruled]{algorithm2e}

\AddToHook{selectfont}{%
\abovedisplayskip=5pt\relax
\belowdisplayskip=5pt\relax
\abovedisplayshortskip=2pt\relax
\belowdisplayshortskip=2pt\relax
}

\usepackage[colorlinks=true,linkcolor=RoyalPurple,citecolor=Blue,urlcolor=Violet,linktoc=page]{hyperref}
\title{Benign population landscapes do not transfer to empirical risk minimization in nonsmooth optimization
}
\author{%
Penghao Wang\thanks{Edwardson School of Industrial Engineering, Purdue University.}
\and
Mateo D\'iaz\thanks{Department of Statistics \& Data Science Institute, University of Chicago. MD was partially supported by NSF awards CCF 2442615 and DMS 2502377, and a Sloan Research Fellowship.}
\and
Liwei Jiang\footnotemark[1]
}
\date{}

\begin{document}

\maketitle
\begin{abstract}
  Benign landscapes provide a geometric explanation for the success of iterative methods on nonconvex problems arising in data science.
The landscape is benign if every local minimizer is global and every other critical point is a strict saddle. For methods that avoid strict saddles, convergence to a critical point therefore implies convergence to a global minimizer.
% The landscape is benign whenever local minimizers are global and all other critical points are strict saddles. When this happens, saddle-point avoidance guarantees turn convergence to a critical point into convergence to a global minimizer.
For smooth stochastic optimization, a common strategy for establishing such landscapes proceeds in two steps: $(i)$ establish benign geometry for the population objective---which is often easier to analyze---and $(ii)$ transfer this geometry to the empirical objective using concentration of gradients and Hessians. Recent work has extended saddle-point avoidance guarantees to weakly convex nonsmooth optimization through the notion of active strict saddles. We show, however, that benign population geometry need not transfer to empirical objectives in this setting. Specifically, for robust phase retrieval with Gaussian measurements---a well-studied and highly structured problem---the population objective is benign. Nevertheless, the empirical objective has sharp spurious local minima with probability bounded away from zero as the sample size tends to infinity. This asymptotic probability lower bound is independent of the dimension. %Numerical experiments support these findings.

%%% Local Variables:
%%% mode: LaTeX
%%% TeX-master: "main"
%%% End:

\end{abstract}

%\tableofcontents

\section{Introduction}\label{sect:intro}
Simple iterative methods, such as gradient descent, have proven to be highly effective at solving nonconvex optimization problems stemming from data science, machine learning, and signal processing. A working geometric explanation for this success is the benign landscape of the objective: all local minima are global, and every other critical point is a strict saddle---meaning that its Hessian has a strictly negative eigenvalue. Under this geometry, many suitably tuned first-order methods escape strict saddles and consequently recover global minimizers whenever they converge \cite{lee2019first,mucsat2025gradient}. Further, stochastically perturbed variants of gradient descent enjoy strong nonasymptotic guarantees \cite{jin2017escape}. This geometric explanation is supported by a number of high-impact problems that exhibit benign landscapes under suitable assumptions \cite{sun2018geometric,ge2017no,boumal2016nonconvex,wang2020efficient,zhu2020global,ge2015escaping,mcrae2026phase,mcrae2024orthogonal,mcrae2025synchronization}. %, including smooth phase retrieval~\cite{sun2018geometric}, low-rank matrix sensing and completion \cite{ge2017no}, phase synchronization \cite{boumal2016nonconvex}, clustering \cite{wang2020efficient}, shallow linear neural networks \cite{zhu2020global}, and orthogonal tensor decomposition \cite{ge2015escaping}.

Most of these examples arise from statistical models and admit both population and empirical objectives. Formally, for samples $z_1, \dots, z_n$ from some distribution, these objectives are defined as
\begin{equation*}
  \F(x) := \EE_z f(x, z), \qquad \text{and} \qquad \Fn(x) := \frac1n \sum_{i=1}^n f(x, z_i),
  \end{equation*}
where $f$ is problem specific. Although the objective being minimized is often $\Fn$, the population objective $\F$ tends to be simpler to analyze. This motivates a common two-stage approach to establish benign geometry: first, show that it holds for the population objective $\F$ and then prove that this geometry persists for the empirical objective $\Fn$ whenever $n$ is sufficiently large. The second step is carried out by exploiting the concentration of gradients and Hessians of $\Fn$ around those of $\F$.

The seminal work of Davis and Drusvyatskiy \cite{davis2019active} extended saddle-avoidance guarantees to nonsmooth weakly convex optimization---a broad class encompassing many problems arising in data science \cite{davis2019modelbased}. Unlike in the smooth setting, a classical Hessian is not available for these functions. To overcome this obstacle, they introduce the notion of an active strict saddle. In a nutshell, these are critical points lying on a manifold along which the objective is smooth---and has a direction of strictly negative curvature---and away from which the function grows sharply; see Figure~\ref{fig:active-strict-saddle} for an illustration and Section~\ref{sec:prelims} for a formal definition. Just as in the smooth setting, several first-order methods are guaranteed to avoid these active strict saddles \cite{davis2019active,davis2026active,bianchi2024stochastic,davis2022escaping, qiu2026stochastic}. These results provide a reason to seek nonsmooth objectives with benign landscapes and raise a natural question:
\mybox{\centering \textit{Does benign population geometry transfer to the empirical landscape\\in nonsmooth optimization?}}

% Keep the compact spacing local to Figure 1 and its wrapped text.
\begingroup
\setlength{\intextsep}{3pt}
\begin{wrapfigure}{r}{0.40\textwidth}
    \centering
    \includegraphics[width=\linewidth]{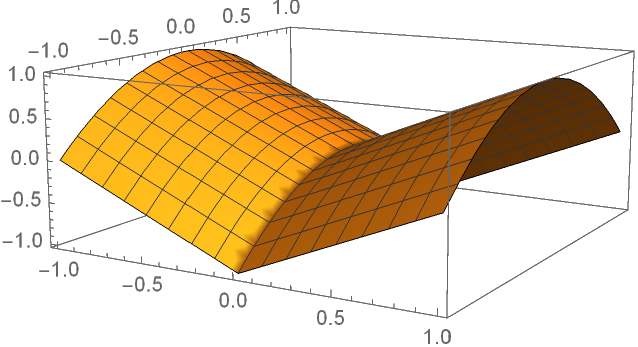}
    \caption{The function $(x,y)\mapsto|x|-y^2$ has an active strict saddle at zero, with the $y$-axis as the active manifold.}
    \label{fig:active-strict-saddle}
\end{wrapfigure}

We answer this question negatively: for general weakly convex problems, this is not the case. To do so, we consider phase retrieval, an inverse problem with a long history in physics and signal processing \cite{fienup1982phase,shechtman2015phase}. Although phase retrieval is commonly formulated over the complex numbers, we study a stylized real version in which the goal is to recover a vector $\xs \in \RR^d$, up to sign, from $n$ nonlinear measurements of the form $\langle a_i, \xs\rangle^2$, where the vectors $a_i$ are i.i.d. standard Gaussians. To solve this problem, we consider the robust formulation
\begin{equation}\label{eq:emp}
  \Fn(x) = \frac{1}{n}\sum_{i =1}^n \left| \dotp{a_i, x}^2 - \dotp{a_i, \xs}^2 \right|.
\end{equation}
%We normalize $\norm{\xs}=1$ by scaling, as explained in Section~\ref{sec:prelims}.
This objective is highly structured: it is weakly convex, a single-index model, and exhibits composite structure, i.e., it is a composition of a convex function with a smooth map \cite{duchi2019solving}. Further, this formulation and related low-rank matrix recovery formulations have been studied extensively and serve as standard testbeds for theory and algorithms in weakly convex optimization \cite{davis2020nonsmooth,davis2018subgradient,davis2019modelbased,diaz2019blind,charisopoulos2021bilinear,charisopoulos2021lowrank,ruan2026uniform}.
\par\endgroup
  %it has beenthat has been studied extensively studied in the literature of weakly-convex optimization \cite{}.

\paragraph{Population objective.} The population counterpart of~\eqref{eq:emp} and its first-order landscape were characterized in \cite{davis2020nonsmooth}. For a unit-norm $\xs$, the critical points of $\F(x) = \EE |\dotp{a, x}^2 -\dotp{a, \xs}^2|$ are \begin{equation}
    \left\{\pm \xs\right\} \,\cup \,\{0\} \,\cup \,\cS \qquad \text{with} \qquad \cS \coloneqq \left\{\tau u \mid \|u\| =1 \text{ and } u \perp \xs \right\},
\end{equation}
where $\tau \approx 0.4416$ is the unique solution of a certain equation. That is, the critical points correspond to the true underlying solutions $\pm \xs$, zero, and a $(d-2)$-dimensional sphere orthogonal to the solution. We show here that $\F$ is twice continuously differentiable on $\RR^d \setminus \{\pm \xs\}$. %Moreover,
%\begin{equation}
%\nabla^2 \F(0) = -2 I, \qquad \text{and} \qquad \nabla^2 \F(\tau u) = \frac{16}{\pi} \cdot \frac{\tau}{(1+\tau^2)^2} \cdot \left(uu^\top - {\xs}{\xs}^\top\right) \quad \text{for all }\tau u \in \cS.
%\end{equation}
Further, zero and the vectors on the sphere $\cS$ are all standard saddles. Close to the ground-truth solutions, the population objective is sharp, i.e., the objective can be bounded below by a constant multiple of the distance to the minimizer set. % there is a constant $\mu > 0$ such that for any $x$ near $\xs$ we have
%\begin{equation}
%\F(x) - \F(\xs) \ge \mu \cdot \|x - \xs\|.
%\end{equation}
In turn, the landscape is benign: the only local minimizers are the ground-truth solutions, and all other critical points are (standard) strict saddles---which are consequently active strict saddles with an active manifold corresponding to the full domain.

\paragraph{Empirical landscape.} The geometry of $\Fn$ differs substantially from that of its population counterpart. We establish three complementary results.
\begin{enumerate}
        \item (\textbf{Active strict saddles}) Suppose $n \geq 2d$. Then, almost surely,
\begin{quote} every critical point of $\Fn$ is either a strict local minimizer or an active strict saddle. \end{quote}
        This dichotomy classifies the critical points, but leaves open whether every local minimizer is global. Our next result shows that this need not be the case.

        \item (\textbf{Expected number of spurious minimizers}) Let $M_{d,n}$ denote the number of sharp spurious local minimizers. Here, sharpness means that the objective grows at least linearly with distance from each such minimizer in a sufficiently small neighborhood. For every fixed dimension $d \geq 4$, we show that
$$
\EE M_{d, n} = \Omega\left(n^{d/2-1}\right) \quad \text{as } n \to \infty.
$$
        Thus, the expected number of spurious minimizers grows at least polynomially with the sample size. However, this alone does not rule out their absence with high probability: the growing expectation could be driven by rare events with very many spurious minimizers.
        \item (\textbf{Probability of having a spurious minimizer}) We therefore establish a stronger guarantee:
\[
\liminf_{n \to \infty} \PP(\Fn \text{ has a spurious sharp local minimum}) > 0.136, \quad \text{for any fixed } d\geq4.
\]
        This lower bound is independent of dimension. Moreover, the spurious minimizers constructed in our proof approach the population saddle set $\cS$ as $n\to\infty$. Thus, sharp empirical local minima can accumulate at population strict saddles, even as the sample size grows.
        \end{enumerate}

        \paragraph{Intuition.} Our results may initially seem counterintuitive in light of the population-to-empirical landscape transfer familiar from smooth optimization. The key difference is that nonsmooth empirical perturbations can introduce linear growth that dominates the negative curvature of the population objective. To illustrate this point, consider the one-dimensional empirical objective
        $$H_n(t)=-t^2+\frac{1}{n}\sum_{i=1}^{n}\alpha_i|t|,$$
        where $\alpha_1,\dots,\alpha_n$ are i.i.d. standard Gaussian random variables. The population objective is $\EE H_n(t)=-t^2$, so the origin is a strict maximizer and hence a strict saddle. Yet, for every finite $n$, the empirical coefficient $n^{-1}\sum_{i=1}^n\alpha_i$ is positive with probability $1/2$. On this event, the linear term dominates the negative quadratic term sufficiently close to the origin, creating a sharp local minimum. Although the perturbation vanishes as $n$ grows, the probability of creating this local minimum remains unchanged. This illustrates the mechanism behind the sharp local minima that arise near the population saddle set $\cS$.
%        \paragraph{Proof technique.} The nonsmooth residuals partition the space by hyperplanes, and intersections of $d$ such hyperplanes form candidate vertices. We first characterize when one of these vertices is a sharp local minimizer: the contribution of the inactive measurements must lie strictly inside the parallelotope generated by the active measurements. We then count a regular family of such vertices in an $n^{-1/2}$-neighborhood of the population saddle sphere $\cS$. A weighted Kac--Rice formula, together with a multivariate local limit theorem for the inactive contribution, yields the first-moment asymptotics for $N_{d,n}$. For the second moment, we apply a two-point Kac--Rice formula and separately control pairs with overlapping active measurements and pairs whose directions are nearly equal or antipodal. The remaining limiting expression is governed by the covariance of a Gaussian random field on $\cS$; quantitative angular correlation bounds make the ratio of the first two moments uniform in the dimension. Finally, the Paley--Zygmund inequality converts this moment ratio into a dimension-independent lower bound on the probability of observing a sharp spurious local minimizer.

        \paragraph{Outline.} The remainder of this section covers related work. Section~\ref{sec:prelims} introduces the necessary background for the rest of the paper. Section~\ref{sec:pop} presents our results on the population and empirical landscapes, and Section~\ref{sec:proof-strategy} outlines the proof strategy. Section~\ref{sec:experiments} studies algorithmic recovery from perturbations of the population saddle sphere. Finally, Section~\ref{sec:conclusion} summarizes our conclusions and discusses limitations and future directions. All technical proofs are deferred to the appendix.

\subsection*{Related work}

\paragraph{Population-to-empirical transfer.}
For smooth losses, \cite{foster2018uniform, mei2018landscape} established general results for the uniform concentration of gradients and Hessians under mild assumptions. For weakly convex optimization, \cite{davis2022graphical} showed the graphical convergence of subdifferentials for weakly convex stochastic problems. Similarly, \cite{ruan2026uniform} studied uniform convergence of the subdifferentials for a more structured family of weakly convex functions and used it to analyze robust phase retrieval \cite{davis2020nonsmooth}. To our knowledge, our result is the first to identify a second-order obstruction in this setting. Beyond weakly convex optimization, \cite{shapiro2007uniform} established uniform laws of large numbers for random set-valued mappings and applied them to the consistency of stationary points in nonsmooth stochastic programming. On the other hand, \cite{tian2025failure,tian2026lipschitzian} provided negative and positive results for uniform convergence of Clarke subdifferentials of Lipschitz objectives.

\paragraph{Nonsmooth landscapes.}
For symmetric rank-one matrix factorization with entrywise $\ell_1$ loss, the recent works \cite{guan2024ell_1, josz2025nonsmooth} established that every second-order stationary point, defined through the second subderivative, is globally optimal. Moreover, \cite{guan2024ell_1} showed that nonglobal stationary points may lack active manifolds and that arbitrarily small linear perturbations can lead to spurious local minima.
More recently, \cite{wang2026absence} developed a deterministic second-order branching criterion that rules out the existence of identifiable manifolds at critical points of finite-max composite problems. They showed that this criterion holds with high probability in overparameterized robust low-rank recovery and almost surely at common interpolators of random minimax regression, suggesting that the absence of identifiable manifolds can be typical in modern nonsmooth models. This line of work studies the failure of manifold identification at a fixed critical point, whereas our obstruction is the creation of sharp empirical local minima near a population saddle.
% For sphere-constrained orthogonal dictionary learning, \cite{bai2019subgradient} use uniform subdifferential concentration to transfer directional recovery bounds to empirical objectives on selected regions.
The asymmetric analogue of this matrix factorization formulation was studied in \cite{josz2022nonsmooth}, which characterized precisely when it has no spurious local minima. For orthogonal dictionary learning, \cite{bai2019subgradient} used uniform concentration of subdifferentials to establish recovery guarantees for a randomly initialized Riemannian subgradient method applied to a sphere-constrained $\ell_1$ objective.
For the nonsmooth amplitude least-squares formulation of phase retrieval, which differs from~\eqref{eq:emp}, McRae~\cite{mcrae2026amplitude} establishes benign empirical landscape guarantees with constant rank overparameterization. The paper leaves open whether the original  formulation without overparameterization has spurious local minima.
The authors of \cite{zheng2024smoothed} established the benign empirical landscape of a smoothed version of \eqref{eq:emp} with a sufficiently small fixed positive smoothing parameter. Intuitively, this landscape seems to smooth out the sharp kinks we identify in the purely nonsmooth objective. The second-order analysis and possible existence of spurious local minima for robust phase retrieval \eqref{eq:emp} were explicitly identified as relevant unresolved problems in \cite{zheng2024smoothed, guan2024ell_1}.% highlighted the need for further classification of its stationary points. Lastly, . %These bounds support recovery without asserting that every empirical local minimum is global. Such positive results underscore the importance of the model and the particular geometric properties being transferred.

\paragraph{Nonsmooth algorithmic guarantees.}
Sharpness, together with Lipschitzness and weak convexity, yields fast local convergence of subgradient and prox-linear methods~\cite{polyak1969minimization, davis2018subgradient,duchi2019solving,charisopoulos2021lowrank}. These convergence results require initialization near the solution set and do not deal with stationary points outside that region. Beyond local convergence, \cite{jin2018local} developed a smoothing algorithm to find approximate second-order stationary points of a smooth population risk from a sufficiently accurate, potentially nonsmooth uniform approximation. To our knowledge, all general-purpose guarantees for convergence to second-order critical points of nonsmooth objectives rely on active strict saddles \cite{davis2019active,bianchi2024stochastic,davis2026active,davis2022escaping,qiu2026stochastic}. Our results suggest that we need to reevaluate this algorithmic explanation if we wish to explain the empirical success of first-order methods on statistical recovery problems.
%Our empirical classification establishes the active-saddle structure directly for the original random-data model. The sharp spurious minima we exhibit show why excluding active saddles from possible limits does not suffice to deduce global recovery.

% \subsection*{AI disclaimer}
%%% Local Variables:
%%% mode: LaTeX
%%% TeX-master: "../main"
%%% End:

\section{Preliminaries}\label{sec:prelims}
In this section, we introduce the notation and background required for the rest of the paper.

\paragraph{Basic notation.}
We use $\dotp{\cdot,\cdot}$ and $\norm{\cdot}$ for the Euclidean inner product and norm, and $\norm{\cdot}_{\mathrm{op}}$ for the operator norm. Set addition and scalar multiplication are understood elementwise: $A+B=\{a+b:a\in A,b\in B\}$ and $Sv=\{sv:s\in S\}$ for $S\subset\R$. Given a differentiable map $g \colon \RR^d \to \RR^m$, we use $\nabla g$ to denote its Jacobian. When $m = 1$, this corresponds to the gradient map, and, further, when $g$ is a twice differentiable function, the symbol $\nabla^2g$ denotes its Hessian map.

\paragraph{Manifolds.} Our treatment of manifolds here is rather brief since we will only focus on Euclidean manifolds, i.e., affine subspaces. We refer the interested reader to the textbooks by Boumal~\cite{boumal2023intromanifolds} and Lee~\cite{lee2012smooth}. A $C^p$ manifold $\cM\subset\R^d$ is locally a level set $G^{-1}(0)$ of a $C^p$ map whose Jacobian has full row rank. Its tangent space at $x$ is $T_{\cM}(x)=\ker \nabla G(x)$. When $\cM$ is an affine subspace, $T_{\cM}(x)$ is equal to the same affine subspace shifted to the origin. We endow $\cM$ with a Riemannian inner product induced by the Euclidean inner product. Given a $C^2$ function $g\colon \cM \to \RR$, we use the symbol $\nabla_{\cM}^2g(x)$ to denote the Riemannian Hessian at $x$. When $\cM$ is affine, $\nabla^2_{\cM}g(x)$ is represented by $V^\top\nabla^2\widetilde g(x)V$, where the columns of $V$ form an orthonormal basis of $T_{\cM}(x)$ and $\widetilde g$ is any $C^2$ extension of $g$ to a neighborhood of $x$.

\paragraph{Nonsmooth analysis.}
A proper closed function $f\colon\R^d\to\R\cup\{+\infty\}$ is $\rho$-weakly convex if $x\mapsto f(x)+\frac{\rho}{2}\norm{x}^2$ is convex, where $\rho\ge0$. For $x\in\dom f$, its subdifferential $\partial f(x)$ consists of the vectors $v\in\R^d$ satisfying
\[
f(y)\ge f(x)+\dotp{v,y-x}-\frac{\rho}{2}\norm{y-x}^2
\qquad\text{for all }y\in\R^d.
\]
We call $x$ a critical point if $0\in\partial f(x)$. This class includes compositions $h\circ G$ with $h$ convex and $L$-Lipschitz and $G$ having a $\beta$-Lipschitz Jacobian: such functions are $L\beta$-weakly convex~\cite{davis2018subgradient}.
For the empirical objective $\Fn$ in~\eqref{eq:emp}, we show in Theorem~\ref{thm:active strict saddle} that $\Fn$ is $2\norm{\frac1n\sum_{i=1}^n a_i a_i^\top}_{\mathrm{op}}$-weakly convex for every realization of the measurements.

Nonsmooth optimization admits several different notions of second order information~\cite{rockafellar1998variational}. We adopt the notion of an active strict saddle introduced by Davis and Drusvyatskiy~\cite{davis2019active}, which underlies saddle-avoidance guarantees for several first-order methods~\cite{bianchi2024stochastic,davis2026active}.
To introduce this notion, we require an auxiliary definition.
\begin{definition}[Active manifold]\label{def:active-manifold}
Consider a closed, weakly convex function $f$, and let $\bar x$ be a critical point. A set $\cM$ containing $\bar x$ is an \emph{active $C^p$ manifold} at $\bar x$ if there is some neighborhood $U$ of $\bar x$ for which the following holds.
\begin{enumerate}
\item (\textbf{Smoothness}) The set $\cM\cap U$ is a $C^p$ manifold and the restriction $f|_{\cM\cap U}$ is $C^p$-smooth.
\item (\textbf{Sharpness}) Subgradients are bounded away from zero off the manifold,
\[
\inf\{\norm{v}:v\in\partial f(x),\ x\in U\setminus\cM\}>0.
\]
\end{enumerate}
\end{definition}

The notion of ``active sets'' has a long history in nonsmooth optimization \cite{alkhayyal1991finite, burke1988identification, burke1990identification,dunn1987convergence,ferris1991finite,flam1992finite, shapiro2003composite, drusvyatskiy2014optimality}. Beyond the definition we follow here, this notion has been modeled in a variety of ways, including identifiable surfaces~\cite{wright1993identifiable}, partial smoothness~\cite{lewis2002active}, $\mathcal{V}\mathcal{U}$-decompositions~\cite{lemarechal2000lagrangian}, and cone reducibility~\cite{bonnans2000perturbation}. Furthermore, these notions have led to a range of algorithmic developments \cite{mifflin2005vu,liang2017activity,lewis2021active,davis2024asymptotic,davis2025local,kong2026lipschitz}. Armed with Definition~\ref{def:active-manifold} we can now introduce active strict saddles.
\begin{definition}[Active strict saddles]
A critical point $\bar x$ of a weakly convex function $f$ is an \emph{active strict saddle} if it admits an active $C^2$ manifold $\cM$ such that $\nabla^2_{\cM}(f|_{\cM})(\bar x)$ has a negative eigenvalue.\footnote{This is equivalent to the parabolic-subderivative definition in~\cite{davis2019active}; see~\cite[Definition~4 and footnote~2]{bianchi2024stochastic}.}
A function $f$ exhibits the \emph{active strict saddle property} if every critical point is either a local minimizer or an active strict saddle.
\end{definition}
Intuitively, at an active strict saddle, the function has a direction of negative curvature along the active manifold and grows sharply away from it. The existence of an active manifold appears crucial for saddle avoidance by simple iterative methods: without one, subgradient methods may converge to critical points that are not local minimizers~\cite{davis2019active}. Finally, a function $f \colon \RR^d \to \RR \cup \{+\infty\}$ with a minimizer set $\cX_\star$ is sharp if there exists a constant $\mu > 0$ such that
\begin{equation}\label{eq:sharp}
  f(x) - \inf f \geq \mu \,\dist(x, \cX_\star) \quad \text{for all }x\in \RR^d.
\end{equation}
A point $\bar x\in\dom f$ is a \emph{sharp local minimizer} if there exist constants $\mu,r>0$ such that
\[
f(x)-f(\bar x)\ge\mu\norm{x-\bar x}
\qquad\text{whenever }\norm{x-\bar x}<r.
\]
Whenever $f$ is $\rho$-weakly convex, $\mu$-sharp and $L$-Lipschitz, then several first order methods converge fast towards a minimum provided they are initialized near enough \cite{davis2018subgradient,drusvyatskiy2019efficiency}.
  % Thus, an active strict saddle has negative curvature along a smooth manifold and sharpness off it.

%%% Local Variables:
%%% mode: LaTeX
%%% TeX-master: "../main"
%%% End:

\section{The second-order nonsmooth landscape of phase retrieval}\label{sec:pop}
In this section, we formally describe our results on the second-order landscape of robust phase retrieval. The first subsection covers the population objective, while the second subsection covers the empirical objective. Throughout, we make the blanket assumption that the signal we wish to recover has norm one, i.e., $\|\xs\| = 1.$ The assumption is made only for convenience, as it simplifies expressions that would otherwise depend on its norm.

\subsection{Population landscape} We start by showing that $\F$ is $C^2$-smooth everywhere but at the solutions and that each critical point is either a global minimum or a saddle. The proof of this proposition appears in Appendix~\ref{pf:population landscape}.%This extends the first-order characterization from \cite{davis2020nonsmooth}.
\begin{proposition}\label{prop:population landscape}
  Fix a unit-norm vector $\xs \in \RR^d$ and consider the population objective $\F(x) = \EE |\dotp{a,x}^2 - \dotp{a,\xs}^2|$ with $a \sim N(0, I_d).$ Then, $\F(\cdot)$ is a closed, weakly convex function, and its critical points correspond to
  \begin{equation}\label{eq:population-landscape}
    \left\{\pm \xs\right\} \,\cup \,\{0\} \,\cup \,\cS \qquad \text{with} \qquad \cS \coloneqq \big\{\tau u \mid \|u\| =1 \text{ and } u \perp \xs \big\},
  \end{equation}
  where $\tau$ is the unique solution of the equation $\arctan \tau + \tau/(1+\tau^2) = \pi/4.$
%The critical set of $\F$ is \[
%\{0,e,-e\}\cup \mathcal{S}\]
The function $\F$ is $C^2$ on $\R^d \setminus \{\pm \xs\}$, and its Hessian evaluates to \begin{equation}\label{eq:population-hessians}
\nabla^2 \F(0)=-2I_d, \quad \text{and} \quad \nabla^2 \F(\tau u)=\frac{16}{\pi} \cdot \frac{\tau}{(1+\tau^2)^2} \cdot (uu^\top-\xs {\xs}^\top) \quad \text{for all } \tau u \in \mathcal{S}.
\end{equation}
Further, $\F$ is $2/\pi$-sharp with global minimizers $\pm\xs$.
\end{proposition}
The characterization of critical points in~\eqref{eq:population-landscape} was first established in~\cite{davis2020nonsmooth}. As direct consequences, the origin is a strict local maximizer, while each point $\tau u\in\cS$ is a saddle with positive curvature in the radial direction $u$, negative curvature along $\xs$, and zero curvature in directions tangent to $\cS$. Since $\F$ is $C^2$ near these nonglobal critical points, each is also an active strict saddle, with active manifold equal to the whole space $\RR^d$. Together with the sharpness of the global minimizers $\pm\xs$, this establishes the benign geometry of the population landscape.
%The saddle sphere has dimension $d-2$. The curvature of $\F$ at $x=\tau u$ is negative along $e$ and positive along $u$. Its tangent directions have zero curvature. Appendix~\ref{pf:population landscape} implies the explicit population calculation for Proposition \ref{prop:population landscape}.

\subsection{Empirical landscape}
In this subsection, we state our main results for the empirical objective $\Fn$, defined in \eqref{eq:emp}.
%We first describe the active manifolds of the empirical objective $\Fn$, defined in \eqref{eq:emp}, and then state its active strict saddle property.
To describe these guarantees, we first define the active index set at a vector $x\in\R^d$ as
\begin{equation}\label{eq:active}
I(x)=\bigl\{i\in\{1,\ldots,n\}:\dotp{a_i,x}^2=\dotp{a_i,\xs}^2\bigr\}.
\end{equation}
As it turns out, when $n\ge 2d$, almost surely every critical point $x\not\in \{\pm \xs\}$ has an active manifold encoded via $I(x)$. The proof of the following theorem appears in Appendix~\ref{pf:active strict saddle}. %; the criticality characterization used in that proof is stated in Section~\ref{sec:proof-strategy}.
\begin{thm}[Active strict saddle property]\label{thm:active strict saddle}
	Fix a unit-norm vector $\xs \in \RR^d$. The empirical objective $\Fn$ in \eqref{eq:emp} is $2\norm{\frac1n\sum_{i=1}^n a_i a_i^\top}_{\mathrm{op}}$-weakly convex for any realization of the measurements $\{a_i\}$. If $n\ge 2d$, then almost surely every critical point $x_\sharp$ of $\Fn$ is either a strict local minimizer or an active strict saddle.
	In the latter case, the active manifold is
    \begin{equation}
      \cM_{x_\sharp} = \left\{ x \in \RR^d \mid \dotp{a_i, x} = \dotp{a_i, x_\sharp} \,\, \text{for all }i \in I(x_\sharp)\right\}.
      \end{equation}
\end{thm}
% The proof appears in Appendix~\ref{pf:active strict saddle}.
Let us make a few remarks. This guarantee, in tandem with the saddle-avoidance results from~\cite{bianchi2024stochastic,davis2026active}, shows that certain stochastically perturbed variants of the subgradient method can converge only to local minimizers of the empirical objective $\Fn$. Since any subset of $d$ Gaussian vectors $\{a_i\}$ is almost surely linearly independent, at every critical point $x_\sharp\notin\{\pm\xs\}$ the index set satisfies $|I(x_\sharp)| \leq d$, and moreover, when equality holds, the manifold reduces to a single point and $x_\sharp$ is a sharp local minimizer. This simple observation will be key in our later developments.

Theorem~\ref{thm:active strict saddle}, together with the benign population landscape described in Proposition~\ref{prop:population landscape}, might give the impression that the empirical objective exhibits benign geometry as well. The next result shows that this is not the case. Recall that $M_{d,n}$ denotes the random variable that counts the number of spurious sharp local minima of $\Fn.$ We defer the proof of this guarantee to Appendix~\ref{app:probability}.

\begin{thm}[\textbf{Existence of spurious sharp local minima}]\label{thm:probability}
Fix a unit-norm vector $\xs \in \RR^d$, consider the empirical objective $\Fn$ in \eqref{eq:emp}, and suppose that $d \geq 4$. Then, there exists a positive constant $c_d$  depending solely on $d$ such that the following hold true.
\begin{enumerate}
  \item\label{thm:probability:1} (\textbf{Bounded moments}) The first moment is asymptotically bounded by
        \begin{equation}
          \liminf_{n \to \infty} \frac{\EE M_{d, n}}{ n^{d/2-1}} \geq c_d .
% \qquad \text{and} \qquad \liminf_{n \to \infty} \EE N_{d,n}^2 \leq \gamma_d \cJ_d \cdot n^{d - 2}.
          \end{equation}
        \item\label{thm:probability:2} (\textbf{Constant probability}) Furthermore, \[
\liminf_{n\to \infty} \PP(F_n \text{ has a sharp spurious local minimum})\ge 0.136.
        \]
        \end{enumerate}
\end{thm}
Remarkably, the limiting probability lower bound is uniform over $d\ge4$. %Thus, for every fixed dimension, sharp spurious local minima persist with probability bounded away from zero even as the number of measurements $n$ tends to infinity!
Our proof constructs these sharp local minima within a shrinking tube around the population saddle sphere $\cS$, so their distance to this saddle set tends to zero as $n$ increases. The next section outlines the proof strategy. Section~\ref{sec:experiments} studies algorithmic recovery from perturbations of $\cS$.
%This suggests a familiar strategy for establishing global convergence to the ground truth from an arbitrary initialization: as in many smooth problems, one seeks to rule out spurious local minima by showing that every local minimizer is global. However, this strategy does not carry over to the nonsmooth formulation considered here: sharp spurious local minimizers persist as the sample size grows. Section~\ref{sec:proof-strategy} explains the counting argument that establishes their increasing expected number and a positive probability of their existence.

%%% Local Variables:
%%% mode: LaTeX
%%% TeX-master: "../main"
%%% End:

% \section{Empirical landscape}\label{sec:emp}
% \subsection{Empirical landscape I: active strict saddle property}
% \input{sections/empirical.tex} 

% \subsection{Empirical landscape II: sharp spurious local minima}\label{sec:probability}
% \input{sections/spurious-minima.tex}

\section{Proof strategy}
\label{sec:proof-strategy}
In this section, we sketch our proof strategy for establishing Theorem~\ref{thm:probability}. A full proof appears in Appendix~\ref{app:probability}. Before starting, let us provide a roadmap. The key idea is to define a random variable $N_{d,n}$ that counts a selected family of sharp spurious local minima inside a tube, denoted $T_n$, around the sphere of saddles $\cS$. This tube shrinks as $n$ goes to infinity. Since $M_{d,n}$ counts all sharp spurious local minima, it follows that
\begin{equation}\label{eq:lower-bound}
M_{d,n}\ge N_{d,n},\qquad \text{and so} \qquad
\PP(M_{d,n}>0)\ge\PP(N_{d,n}>0).
\end{equation}
Thus, it suffices to lower bound the probability that our selected family contains at least one minimum. The proof progresses in three steps: $(i)$ define the family counted by $N_{d,n}$, $(ii)$ compute the asymptotic first moment of $N_{d,n}$, $(iii)$ compute the second moment of $N_{d,n}$ and use the Paley--Zygmund inequality to obtain the desired probability bound. Throughout, we assume that $d\ge4$ is fixed, and let $n$ tend to infinity.

\paragraph{A family of spurious sharp minimizers.}
 The following lemma, proved in Appendix~\ref{app:proof-characterization-critical}, characterizes critical points and sharp local minima. Recall the active index set $I(x)$ in \eqref{eq:active}.
\begin{lemma}[Characterization of critical points and sharp local minima]\label{lem:critical point strategy}
Suppose $n\ge2d$. Then, almost surely, every $x\in\R^d\setminus\{\pm\xs\}$ satisfies $|I(x)|\le d$. Such a point is critical if, and only if, there exist scalars $\lambda_i\in(-1,1)$ with $i\in I(x)$, satisfying
\begin{equation}\label{eq:criticality}
 0=\sum_{i\notin I(x)}\sign(\dotp{a_i, x}^2-\dotp{a_i,\xs}^2)\dotp{a_i, x}a_i
 +\sum_{i\in I(x)}\lambda_i\dotp{a_i, x}a_i.
\end{equation}
Moreover, a critical point $x$ is a sharp local minimum if, and only if, it is critical and $|I(x)|=d$.
\end{lemma}
The lemma implies a simple way to count sharp spurious local minima: enumerate the candidates determined by their active indices and check which are critical. We will use this strategy to bound $N_{d,n},$ but first we will describe the family of sharp local minima that $N_{d,n}$ counts.
For $I=\{i_1<\cdots<i_d\}\subset\{1,\ldots,n\}$ and $s\in\{\pm1\}^d$, the active equations $a_{i_j}^\top x=s_j a_{i_j}^\top\xs$, almost surely, determine a unique candidate
\begin{equation}\label{eq:candidate}
 x_{I,s}=A_I^{-1}\operatorname{diag}(s)A_I\xs.
\end{equation}
Here $A_I$ is the matrix with rows $a_{i_j}^\top$. To define $N_{d,n}$, we fix the pattern $s^\circ=(1,\ldots,1,-1,-1)$, with $d-2$ plus signs, and write $x_I=x_{I,s^\circ}$. Counting sharp minima with only this pattern suffices for a lower bound as in \eqref{eq:lower-bound}. Yet, we will count an even smaller subset of these $x_I.$ %: different index sets give distinct candidates, so the resulting count satisfies $N_{d,n}\le M_{d,n}$ and $\PP(M_{d,n}>0)\ge\PP(N_{d,n}>0)$.

We express the criticality condition~\eqref{eq:criticality} in a form that will help us estimate the probability that a candidate $x_I$ is a sharp local minimum. Define $S_I(x)=\sum_{i \notin I}\sign((a_i^\top x)^2-(a_i^\top\xs)^2)(a_i^\top x)a_i$ and let $Q_I(x)$ have columns $(a_i^\top x)a_i$, $i\in I$, in increasing index order. By Lemma~\ref{lem:critical point strategy}, $x_I$ is a sharp local minimum precisely when
\begin{equation}\label{eq:strategy-criticality}
 - S_I(x_I)\in Q_I(x_I)(-1,1)^d,
\end{equation}
where the set on the right consists of the vectors $Q_I(x_I)\lambda$ with $\lambda\in(-1,1)^d$.

To make the probability of~\eqref{eq:strategy-criticality} easier to estimate, we impose two further restrictions. First, we require the selected measurement vectors $\{a_i\}_{i\in I}$ to have bounded norms and the determinants $|\det A_I|$ and $|\det Q_I(x_I)|$ to be bounded below by positive constants. These bounds prevent the equations defining $x_I$ from becoming nearly dependent and ensure that the set on the right-hand side of~\eqref{eq:strategy-criticality} has volume bounded away from zero. The constants depend only on $d$ and are chosen so that these bounds hold with positive probability. We do not include them here since their precise values are not particularly insightful. Second, we count only candidates in the tube
\begin{equation}\label{eq:tube}
 T_n=\left\{n^{-1/2}z_1\xs+(\tau+n^{-1/2}z_2)u\mid
 u\perp\xs,\ \norm u=1,\text{ and }z\in[-1,1]^2\right\}.
\end{equation}
The set $T_n$ is a tube around the population saddle sphere $\cS$ that shrinks to $\cS$ as $n\to\infty$. The width $n^{-1/2}$ matches the fluctuations exhibited by critical points of the empirical objective \cite{ruan2026uniform}.

\paragraph{The first moment.} For any fixed set $I,$ the probability that $x_I$ is counted in the family described above is exactly the same, thus
\begin{equation}\label{eq:binomial-expansion}
\EE N_{d,n} = \sum_{I} \PP(x_I \text{ is counted}) = \binom nd \PP(x_I \text{ is counted}).
\end{equation}
We therefore fix one set $I$ and focus on computing this probability. Our plan to do so is to integrate over the possible locations of the candidate $x_I$. For a proposed location $x$, write $X_I(x)=A_Ix-\operatorname{diag}(s^\circ)A_I\xs$, so that $X_I(x)=0$ precisely when $x=x_I$. Let $H_I(x)$ be a boolean function that is equal to one precisely when the selected measurements $a_i$ satisfy our bounds on their norms and associated determinants, and the condition~\eqref{eq:strategy-criticality} holds with $x_I$ replaced by $x$, and zero otherwise. The Kac--Rice formula \cite{AzaisWschebor2009}  yields the following identity
\begin{equation}\label{eq:first-moment-kac-rice}
 \PP(x_I\text{ is counted})
 =\int_{T_n}p_{X_I(x)}(0)\,
 \E\bigl[|\det A_I|H_I(x)\bigm|X_I(x)=0\bigr]\,\dd x,
\end{equation}
where $p_{X_I(x)}(0)$ is the probability density of $X_I(x)$ evaluated at zero. 
%The condition $X_I(x)=0$ requires the selected hyperplanes to meet at $x$, the determinant accounts for the change of volume between locations and equation values, and $H_I(x)$ checks whether this intersection is a minimum that we count. Since $X_I(x)$ depends only on the selected measurements, the remaining measurements are still independent under this conditioning.

To describe the asymptotic value of \eqref{eq:first-moment-kac-rice}, we establish the following two estimates, with the first holding uniformly over $x\in T_n$,
\begin{equation}\label{eq:asymptotic-terms}
  \E\bigl[|\det A_I|H_I(x)\bigm|X_I(x)=0\bigr] \asymp_d n^{-d/2}, \qquad \text{and} \qquad \int_{T_n}p_{X_I(x)}(0) dx \asymp_d n^{-1}.
\end{equation}
Here $\asymp_d$ denotes comparison up to positive constants depending only on $d$. Before explaining the rationale behind these two results, let us derive an immediate consequence. Recall that $\binom nd \asymp_d n^d$; combining this fact with \eqref{eq:binomial-expansion}, \eqref{eq:first-moment-kac-rice}, and \eqref{eq:asymptotic-terms} yields
  $
\EE N_{d, n} \asymp_d n^d \cdot n^{-d/2} \cdot n^{-1} = n^{d/2-1}.
$
This establishes the asymptotic lower bound in Theorem~\ref{thm:probability}. A detailed computation accounting for the hidden constants yields two positive constants $\gamma_d$ and $\cI_d$ such that
\begin{equation}\label{eq:first-moment}
 n^{-d/2+1} \EE N_{d, n} \to \gamma_d \cI_d  \qquad \text{as }n \to \infty.
  \end{equation}
Here $\gamma_d>0$ accounts for the selected measurements, while $\mathcal I_d>0$ accounts for the Gaussian density across the tube. Their precise expressions and the full calculation appear in Appendix~\ref{sec:first moment asymptotic}.

Let us quickly summarize the arguments behind \eqref{eq:asymptotic-terms}. For a fixed location $x\in T_n$, we condition further on the selected measurement vectors $\{a_i\}_{i\in I}$. This fixes the set $Q_I(x)(-1,1)^d$, while the remaining measurements retain their independent standard Gaussian distributions. When our norm and determinant bounds hold, the determinant factor $|\det A_I|$ is bounded above and below by positive constants independent of $n$, and the set $Q_I(x)(-1,1)^d$ is bounded with volume bounded away from zero. The condition $H_I(x)=1$ then reduces to asking whether $-S_I(x)$ lies in this fixed set. Each summand in $S_I(x)$ has mean $\frac12\nabla\F(x)$. Since $\nabla\F$ vanishes on $\cS$ and is Lipschitz continuous nearby, this mean has norm $O(n^{-1/2})$ on $T_n$, whose points lie within distance $O(n^{-1/2})$ of $\cS$. Summing the $n-d$ terms, the mean of $S_I(x)$ therefore has norm of order at most $\sqrt n$, uniformly over $x\in T_n$. We can show that its standard deviation in every unit direction is also of order $\sqrt n$, uniformly over $x\in T_n$. A local limit theorem implies a Gaussian approximation to the probability of entering this region, yielding a probability of order $n^{-d/2}$ uniformly over the selected measurements satisfying our bounds and the locations in $T_n$. The first estimate in~\eqref{eq:asymptotic-terms} then follows.

For the second estimate in~\eqref{eq:asymptotic-terms}, we consider the volume of $T_n$. The tube restricts the component along $\xs$ to within $n^{-1/2}$ of zero and the length of the perpendicular component to within $n^{-1/2}$ of $\tau$. These two restrictions give its volume order $n^{-1}$. We can prove that the density $p_{X_I(x)}(0)$ is bounded above and below by positive constants on $T_n$, so its integral has order $n^{-1}$.

 \paragraph{A probability bound uniform in dimension.}
 We use the Paley--Zygmund inequality to turn moment estimates into a lower bound on the probability that $N_{d,n}>0$, namely
\begin{equation}\label{eq:paley-zygmund}
 \liminf_{n\to\infty}\PP(N_{d,n}>0)
  \ge \lim_{n\to\infty}\frac{(\E N_{d,n})^2}{\E N_{d,n}^2}.
\end{equation}
Thus, to derive our lower bound we also asymptotically characterize $\E N_{d,n}^2.$ Computing the second moment follows a similar Kac--Rice approach as for the first moment, but the dependence between candidate pairs $(x_I, x_J)$ makes the argument significantly more technical. We defer the proof to Appendix~\ref{sec:second moment asymptotic} and simply state here the resulting estimate
\begin{equation}\label{eq:second-moment}
 n^{-d+2}\E N_{d,n}^2\to \gamma_d^2\mathcal J_d \qquad\text{as }n\to \infty,
\end{equation}
where $\cJ_d > 0$ is finite and depends on $d$. Combining this estimate with \eqref{eq:first-moment} and \eqref{eq:paley-zygmund} implies
\[
 \liminf_{n\to\infty}\PP(N_{d,n}>0)
  \ge \frac{\mathcal I_d^2}{\mathcal J_d},
\]
the $\gamma_d$ factor cancels exactly. By arguing directly with the definitions of $\cI_d$ and $\cJ_d$ we show that this ratio can be further lower bounded by an explicit quantity $c_\star$ that is independent of the dimension. Combining analytic estimates with certified numerical computations, as described in Appendix~\ref{app:numerics}, we show that the value of $c_{\star}$ is lower bounded by $0.1361$.

%%% Local Variables:
%%% mode: LaTeX
%%% TeX-master: "../main"
%%% End:

\section{Numerical experiments}\label{sec:experiments}
In this section, we numerically test our theory. The objective of our experiments is two-fold: $(i)$ to show that for moderate $d$ and $n$, we observe empirically that the probability of having a spurious sharp minimum is large and $(ii)$ to show that despite our negative results, simple subgradient methods do not often get trapped at these minima. We present two numerical experiments supporting these claims. For both experiments we set $\xs = e_{1}$, the first canonical vector. The code to reproduce these experiments can be found at \reproducibilitylink.

% The lower bound in Theorem~\ref{thm:probability} does not quantify the probability that a sharp spurious minimum exists at finite $(d,n)$. Moreover, although certain stochastically perturbed subgradient methods~\cite{bianchi2024stochastic,davis2026active} can avoid active strict saddles almost surely, the existence of spurious local minima leaves open how often such a method recovers a global solution. Therefore, we have two questions to ask:
% \begin{enumerate}
%     \item How often does $\Fn$ have a sharp spurious local minimum at finite $(d,n)$?
%     \item How often does a perturbed subgradient method recover a global solution?
% \end{enumerate}

% We investigate these two questions in the following experiments. By rotational invariance of Gaussian measurements, we always assume $\xs=e_1$ without loss of generality.

\paragraph{Probability of having a spurious sharp minimum.}
We estimate the probability that $\Fn$ has a sharp spurious minimum for $d=4$ and $n\in\{12,16,24,32,48,64\}$. For each $n$, we draw 100 independent problem instances with $n$ i.i.d. Gaussian measurement vectors. For each instance, we enumerate all pairs $(I,s)$ with $|I|=d$ and nonconstant $s\in\{\pm1\}^d$, and test the candidates $x_{I,s}$ defined in~\eqref{eq:candidate}. We count an instance as positive if at least one candidate satisfies the criticality condition~\eqref{eq:criticality} up to the numerical tolerances specified in Appendix~\ref{app:experiments}. Each estimate is the fraction of positive instances, and the shaded band shows $95\%$ confidence intervals.

\begin{figure}[!t]
    \centering
    \begin{minipage}[b]{0.48\linewidth}
        \centering
        \includegraphics[width=\linewidth]{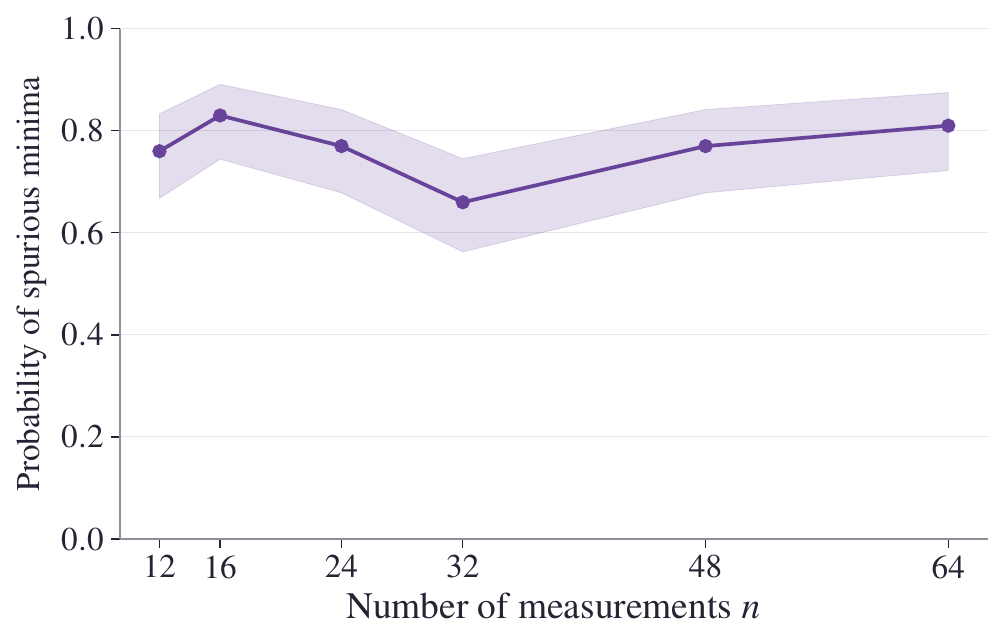}
        \par\smallskip
        (a)
    \end{minipage}
    \hfill
    \begin{minipage}[b]{0.48\linewidth}
        \centering
        \includegraphics[width=0.97\linewidth]{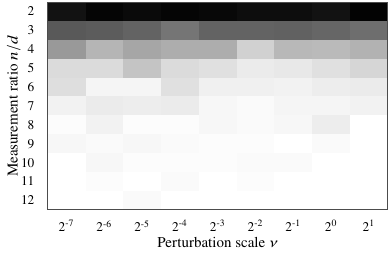}
        \par\smallskip
        (b)
    \end{minipage}
    \caption{(a) Estimated probability of a sharp spurious minimum. (b) Recovery frequency of the perturbed subgradient method. White denotes a recovery fraction of one and black denotes zero.}
    \label{fig:numerical-experiments}
\end{figure}

Figure~\ref{fig:numerical-experiments}(a) shows that the estimates change little with $n$ and remain above $0.5$ across all tested settings. These results are qualitatively consistent with Theorem~\ref{thm:probability}, but suggest that its asymptotic lower bound of $0.136$ is conservative. This may reflect the fact that our proof counts only a strict subset of the sharp spurious minima.

\paragraph{Recovery near the saddle sphere.}
For this experiment, we estimate how often a perturbed subgradient method recovers $\xs$ up to sign. We fix $d=100$ and vary the measurement ratio $n/d\in\{2,3,\ldots,12\}$ and perturbation radius $\nu\in\{2^{-7},\ldots,2^1\}$. For each measurement ratio, we conduct $100$ independent trials. In each trial, we draw $n$ i.i.d. Gaussian measurement vectors. For each radius $\nu$, we draw fresh independent $u\sim\operatorname{Unif}(\cS)$ and $v\sim\operatorname{Unif}(\mathbb{S}^{d-1})$, then initialize $x_0=u+\nu v$. We run the stochastically perturbed subgradient method analyzed in~\cite{davis2026active} for at most $5\times10^5$ updates per initialization, with updates
\[
x_{k+1}=x_k-\frac{0.05}{(k+1)^{0.6}}\bigl(g(x_k)+0.2\xi_k\bigr),
\]
where $g(x_k)\in\partial\Fn(x_k)$ and the vectors $\xi_k\sim\operatorname{Unif}(\mathbb{S}^{d-1})$ are independent. We count a run as recovered if an updated iterate satisfies both $\operatorname{dist}(x_k,\{\pm e_1\})\le10^{-3}$ and $\Fn(x_k)\le10^{-4}$.

% Add \usepackage{graphicx} to your preamble if not already included.

Figure~\ref{fig:numerical-experiments}(b) suggests that the perturbation radius $\nu=\|x_0-u\|$ has little effect on recovery frequency over the tested range. The measurement ratio $n/d$ has a stronger effect, with a sharp increase in recovery near $n/d=4$. Thus, sharp spurious minima near the saddle $\cS$ do not preclude high recovery rates for this method for a sufficiently large number of samples.

%%% Local Variables:
%%% mode: LaTeX
%%% TeX-master: "../iclr2027_conference"
%%% End:

\FloatBarrier

\section{Conclusion and limitations}\label{sec:conclusion}
Our main result shows that even for highly structured weakly convex stochastic objectives, benign population geometry need not transfer to empirical risk minimization. As a consequence, active strict saddle theory alone cannot explain the practical success of iterative methods applied to robust phase retrieval. We conjecture that a similar phenomenon occurs in other nonsmooth formulations for low-rank matrix recovery problems. This opens two directions for future research: $(i)$ to understand which other geometric properties enable these iterative methods to find global solutions for these nonsmooth formulations, and $(ii)$ to characterize general geometric conditions of weakly convex stochastic problems under which population-to-empirical transfer holds. Our analysis has several limitations: it covers only noiseless real Gaussian measurements and takes $n\to\infty$ at each fixed dimension. Although the limiting probability bound is independent of $d$, our proof does not quantify how large $n$ must be relative to $d$ or characterize regimes in which both grow. Extending the analysis to noisy or non-Gaussian measurements would clarify the scope of this obstruction.

%%% Local Variables:
%%% mode: LaTeX
%%% TeX-master: "../main"
%%% End:

% \section{Main results}
% \input{sections/Main_results.tex}

\section{Use of AI in proof discovery}\label{sec:ai}

The last author established Proposition~\ref{prop:population landscape} and
Theorem~\ref{thm:active strict saddle} in 2021. After failing to prove a
population-to-empirical landscape transfer, we questioned whether such a
result could hold and conjectured instead that the empirical objective admits
spurious sharp local minima. We did not make this work public at the time
because the conjecture remained unresolved.

In August 2026, we began interacting with ChatGPT 5.6 Sol and asked it to test
the conjecture numerically. The resulting experiments convinced us that the
conjecture was true, and we began working with the model on a proof. At that
stage, we did not expect to obtain a positive probability lower bound
independent of $d$ as $n\to\infty$ for every fixed $d\geq 4$.

In the subsequent proof discussions, ChatGPT proposed counting spurious sharp
local minima in a tube of width $O(n^{-1/2})$ around the population saddle
sphere. The preliminary argument only yielded a probability lower bound that
tended to zero as $n\to\infty$. We then supplied this argument to
\href{https://github.com/powergiant/qmd-prover}{\texttt{qmd-prover}}, a
proof-development harness developed by Xiao Ma.\footnote{The repository and
its creator are not affiliated with the authors.} We used separate instances
of ChatGPT 5.6 Sol with the ultra reasoning setting as prover and verifier.
These interactions produced draft first- and second-moment arguments based on
a weighted Kac--Rice formula and local limit theorems, leading to the
dimension-independent probability bound in Theorem~\ref{thm:probability}.

In later interactions with GPT-5.6 and GPT-6, we checked the hypotheses of the
probabilistic tools, filled in details for special cases, and reorganized the
dependencies among the lemmas. The authors parsed, completed, and simplified
the resulting arguments and rewrote all proofs in the final manuscript with
human readability in mind. Beyond proof discovery, we used ChatGPT 6.0 Astra
to proofread the manuscript, polish the writing, and implement the numerical
experiments. We reviewed all AI-assisted work and take full responsibility for
the content of this paper.

\bibliographystyle{plain}
\bibliography{reference}

%\clearpage
\appendix
This appendix is rather long, so we start by providing the reader with a roadmap. Appendices~\ref{pf:population landscape}, \ref{pf:active strict saddle}, and~\ref{app:probability} present the proofs of our three main results, Proposition~\ref{prop:population landscape}, Theorem~\ref{thm:active strict saddle}, and Theorem~\ref{thm:probability}, respectively. The proof of the last of these is quite involved and requires several auxiliary results. Appendix~\ref{app:probability-supporting-tools} develops the probabilistic tools used to count spurious sharp local minima, while Appendix~\ref{app:second moment supporting proofs} establishes certain estimates appearing in our characterization of the second moment $\EE N_{d, n}^{2}$.
%The numerical certification of the probability lower bound is described in Appendix~\ref{app:numerics}.
Appendix~\ref{app:auxiliary-lemmas} collects the remaining auxiliary lemmas, and Appendix~\ref{app:experiments} provides the details of the numerical criticality checks in our first experiment.
\ifaiappendix
Finally, Appendix~\ref{app:ai} discusses our use of AI-assisted proof discovery.
\fi

%%% Local Variables:
%%% mode: LaTeX
%%% TeX-master: "../main"
%%% End:

\startcontents[appendix]
\section*{Appendix contents}
\printcontents[appendix]{}{1}{\setcounter{tocdepth}{2}}
%\clearpage
\section{Proof of Proposition~\ref{prop:population landscape}}\label{pf:population landscape}
\begin{proof}
The  set of critical points of $\F$ is characterized in \cite[Theorem~5.3]{davis2020nonsmooth}. The function $\F$ is continuous and hence closed. Using the identity $|u-v| = \max\{u - v, v - u\}$, we derive
\[
\F(x)+\norm{x}^2
=\E\max\{2(a^\top x)^2-(a^\top\xs)^2,(a^\top\xs)^2\},
\]
which is convex since it is the expectation of a maximum of convex quadratics, so $\F$ is closed and weakly convex. Next, we prove the statements about smoothness, the Hessian, and sharpness.

We write \[
\F(x)= \E\left[\abs{(a^\top x)^2- (a^\top\xs)^2}\right]=\E\left[\abs{a^\top (x-\xs)}\abs{a^\top (x+\xs)}\right]= \E[\abs{YZ}],
\]
where $Y=a^\top (x-\xs)$ and $Z=a^\top (x+\xs)$. Both $Y$ and $Z$ are real Gaussian random variables with \[
\E[Y]=\E[Z]=0, \quad \E[Y^2]= \norm{x-\xs}^2, \quad \E[Z^2]=\norm{x+\xs}^2, \quad \text{and} \quad \E[YZ]= \norm{x}^2-1.
\]
Take $x\notin \operatorname{span}\{\xs\}$ and write $x=t\xs+ru$, where $u$ is a unit vector orthogonal to $\xs$ and $r> 0$. Applying \cite[Lemma~1]{IsmailPitman2000} to the normalized $Y$ and $Z$, we obtain \begin{equation}\label{eq:Gaussian absolute identity} \F(x)=\frac{2}{\pi}\left(2r+ (t^2+r^2-1) \arctan \frac{t^2+r^2-1}{2r}\right). \end{equation}
We denote the function above by $\psi(t,r)$. Since $\psi$ is $C^2$ when $r>0$, the chain rule implies that $\F$ is $C^2$ on $\R^d \setminus \operatorname{span}\{\xs\}$. Moreover, for any $x=\tau u$ on the saddle sphere $\mathcal{S}$, a direct computation shows \[
\frac{\partial\psi}{\partial r}(0,\tau)=\frac{\partial^2\psi}{\partial t\,\partial r}(0,\tau)=0
\qquad\text{and}\qquad
\frac{\partial^2\psi}{\partial r^2}(0,\tau)=-\frac{\partial^2\psi}{\partial t^2}(0,\tau)=\frac{16\tau}{\pi(1+\tau^2)^2}.
\]
The chain rule yields
\begin{align*}
\nabla^2 \F(\tau u)
&= \frac{\partial^2\psi}{\partial t^2}(0,\tau)\xs(\xs)^\top
   +\frac{\partial^2\psi}{\partial t\,\partial r}(0,\tau)(\xs u^\top+u(\xs)^\top)\\
&\quad +\frac{\partial^2\psi}{\partial r^2}(0,\tau)uu^\top
   +\frac{1}{\tau}\frac{\partial\psi}{\partial r}(0,\tau)(I_d-\xs(\xs)^\top-uu^\top)\\
&=\frac{16\tau}{\pi(1+\tau^2)^2}(uu^\top-\xs(\xs)^\top).
\end{align*}
Now, let us study smoothness near a point $t_0\xs$ of the axis $\operatorname{span}\{\xs\}$ with $t_0\neq\pm1$. For $x$ near $t_0\xs$, set $t=(\xs)^\top x$, $r=\norm{x-t\xs}$, and $q=\norm{x}^2-1$, which stays away from zero. Using the identity $\arctan s=\sign(s)\,\pi/2-\arctan(1/s)$ for $s\neq0$, and a Taylor expansion of $y-\arctan y$ at $0$, we derive from~\eqref{eq:Gaussian absolute identity} that
\[
\F(x)=|q|+\frac{2}{\pi}\left[2r-|q|\arctan\frac{2r}{|q|}\right]=|q|+r^3h(t,r^2),
\]
where $h$ is smooth near $(t_0,0)$; in particular, the remainder term $r^3h(t,r^2)$ is $O(r^3)$ as $r\to0$. Both expressions remain valid on the axis, where $r=0$ and $\F(x)=|q|\,\E(a^\top\xs)^2=|q|$. Since $t$, $r^2$, and $|q|$ are smooth in $x$, and $r^3=\norm{x-t\xs}^3$ is $C^2$, the function $\F$ is $C^2$ near $t_0\xs$, and hence on $\R^d\setminus\{\pm\xs\}$. At the origin, the same representation implies
\[
\F(x)=1-\norm{x}^2+O(\norm{x}^3)\qquad\text{as }x\to0,
\]
so $\nabla^2\F(0)=-2I_d$.

Finally, $\F\ge 0$ and $\F(\pm\xs)=0$, so $\pm\xs$ are global minimizers.
To establish sharpness, fix $x\neq\pm\xs$ and let $p=\norm{x-\xs}\norm{x+\xs}$. The normalized variables $Y/\norm{x-\xs}$ and $Z/\norm{x+\xs}$ are standard normal with correlation $\rho=(\norm{x}^2-1)/p\in[-1,1]$. The Gaussian product formula~\cite[Lemma~1, equation~(45)]{IsmailPitman2000} therefore implies
\[
\F(x)=\frac{2p}{\pi}\left(\sqrt{1-\rho^2}+\rho\arcsin\rho\right)\ge\frac{2p}{\pi}.
\]
The inequality holds because the expression in parentheses has minimum $1$ at $\rho=0$: its derivative on $(-1,1)$ is $\arcsin\rho$, which is negative for $\rho<0$ and positive for $\rho>0$.
It remains to compare $p$ to the distance to the solution set. The triangle inequality implies $2 \max\{\norm{x-\xs},\norm{x+\xs}\} \geq \norm{x-\xs}+\norm{x+\xs}\ge \norm{2\xs} = 2$. Thus, \[
\begin{aligned}
p
&=\min\{\norm{x-\xs},\norm{x+\xs}\}
  \max\{\norm{x-\xs},\norm{x+\xs}\}\\
&\ge\min\{\norm{x-\xs},\norm{x+\xs}\}\\
&=\dist(x,\{\pm\xs\}).
\end{aligned}
\]
Combining the two bounds with $\min\F=0$, we obtain
\(
\F(x)-\min\F\ge\frac{2}{\pi}\dist(x,\{\pm\xs\})
\).
The same inequality holds, with equality, at $x=\pm\xs$. Hence, $\F$ is $2/\pi$-sharp with global minimizers $\pm x^\star$.
\end{proof}

\section{Proof of Theorem~\ref{thm:active strict saddle}}\label{pf:active strict saddle}
The proof has three parts. First, we characterize the critical points and sharp local minima of $\Fn$ in Lemma~\ref{lem:critical point}. Second, in Proposition~\ref{prop:active manifold}, we use this characterization to establish the existence of an active manifold at each critical point other than $\pm\xs$. Third, we analyze the Hessian along these manifolds to prove the active strict saddle property.

\subsection{Part 1: Characterization of critical points and sharp local minima}\label{app:proof-characterization-critical}

\begin{lemma}[Characterization of critical points and sharp local minima]\label{lem:critical point}
	The empirical function $\Fn$ is $2\norm{\frac1n\sum_{i=1}^n a_i a_i^\top}_{\mathrm{op}}$-weakly convex.
	For every $x\in\R^d$, the subdifferential of $\Fn$ at $x$ is
	\[
	\partial \Fn(x)
	=
	\frac{2}{n}
	\left(
	\sum_{i\notin I(x)}
	\sign\bigl((a_i^\top x)^2-(a_i^\top\xs)^2\bigr)
	(a_i^\top x)a_i
	+
	\sum_{i\in I(x)}
	[-1,1](a_i^\top x)a_i
	\right).
	\]
	Suppose $n\ge 2d$. Then, almost surely, every
	$x\in\R^d\setminus\{\pm\xs\}$ satisfies $\abs{I(x)}\le d$.
	Moreover, such a point $x$ is critical if, and only if, there exist coefficients
	$\lambda_i\in(-1,1)$, $i\in I(x)$, satisfying
	\[
	0
	=
	\sum_{i\notin I(x)}
	\sign\bigl((a_i^\top x)^2-(a_i^\top\xs)^2\bigr)
	(a_i^\top x)a_i
	+
	\sum_{i\in I(x)}
	\lambda_i(a_i^\top x)a_i.
	\]
	Furthermore, such a point $x$ is a sharp local minimum if, and only if, it is critical and $|I(x)|=d$.
\end{lemma}

\begin{proof}
By rotational invariance of the standard Gaussian measurements, we may assume without loss of generality that $\xs=e$, where $e$ is the first standard basis vector.
For every realization of the measurements, $\Fn$ is weakly convex. Indeed, the identity $|s-t|+s=\max\{2s-t,t\}$ implies
\[
\Fn(x)+\frac1n\sum_{i=1}^n(a_i^\top x)^2
=\frac1n\sum_{i=1}^n\max\bigl\{2(a_i^\top x)^2-(a_i^\top e)^2,\,(a_i^\top e)^2\bigr\},
\]
and the right-hand side is convex. Since $\frac1n\sum_{i=1}^n a_i a_i^\top\preceq\norm{\frac1n\sum_{i=1}^n a_i a_i^\top}_{\mathrm{op}}I_d$, it follows that $\Fn$ is $2\norm{\frac1n\sum_{i=1}^n a_i a_i^\top}_{\mathrm{op}}$-weakly convex.
The formula for the subdifferential of $\Fn$ follows directly from the standard chain rule~\cite[Theorem~10.6]{rockafellar1998variational}. For the remainder of the proof, we work on the probability-one event
\begin{equation}\label{eqn:event_eps1}
\ificlr
\begin{aligned}
\Epsilon_1={}&\{a_i^\top e\neq 0\text{ for all }i\in [n]\}\\
&\cap \{(a_i)_{i\in I} \text{ are linearly independent for all }I\subset [n]\text{ with }\abs{I}\le d\}.
\end{aligned}
\else
\begin{aligned}
\Epsilon_1={}&\{a_i^\top e\neq 0\text{ for all }i\in [n]\} \cap \{(a_i)_{i\in I} \text{ are linearly independent for all }I\subset [n]\text{ with }\abs{I}\le d\}.
\end{aligned}
\fi
\end{equation}
\paragraph{Bound on active set size.} For any fixed $I\subset[n]$ with $\abs{I}=d$ and $s\in\{-1,1\}^I$, let $x_{I,s}$ denote the solution of the equations $a_i^\top x=s_ia_i^\top e$, $i\in I$. On $\Epsilon_1$, this solution exists and is unique, and it depends only on $(a_i)_{i\in I}$. Suppose that some $x\neq\pm e$ has more than $d$ active indices. Choose $I\subset I(x)$ with $\abs{I}=d$ and $j\in I(x)\setminus I$. Since every $i\in I$ is active at $x$, there are signs $s_i$ with $a_i^\top x=s_ia_i^\top e$, so $x=x_{I,s}$. Hence the event
\(
\Epsilon_{I,s,j}=\bigl\{x_{I,s}\neq\pm e\text{ and } j\in I(x_{I,s})\bigr\}
\)
occurs. Since $j\notin I$, the vector $a_j$ is stochastically independent of $(a_i)_{i\in I}$ and hence of $x_{I,s}$. Conditional on $(a_i)_{i\in I}$ and on $x_{I,s}\neq\pm e$, the event $j\in I(x_{I,s})$ requires $a_j^\top(x_{I,s}-e)=0$ or $a_j^\top(x_{I,s}+e)=0$. Each is a nonzero linear constraint on a standard Gaussian vector, so $\PP(\Epsilon_{I,s,j})=0$. A union bound over the finitely many triples $(I,s,j)$ shows that, almost surely, $\abs{I(x)}\le d$ for every $x\in\R^d\setminus\{\pm e\}$.

% To prove the bound on the number of active indices, fix an index set $I\subset[n]$ with $\abs{I}=d$ and a choice of signs $(s_i)_{i\in I}$, and consider the equations
% \[
% a_i^\top x=s_i a_i^\top e \qquad \text{with }s_i\in \{-1,1\},\quad\text{and }i\in I.
% \]
% By linear independence, these equations determine a unique solution $x$, which depends only on $(a_i)_{i\in I}$ and the chosen signs. If all signs are $+1$, then $x=e$; if all signs are $-1$, then $x=-e$. For any mixed sign pattern, we have $x\neq\pm e$, since $a_i^\top e\neq0$ for every $i\in I$.

% Now fix a mixed sign pattern and an additional index $j\notin I$. Conditional on $(a_i)_{i\in I}$, the solution $x$ is fixed, while $a_j$ remains an independent Gaussian vector. For $j$ to be active at $x$, we must have either $a_j^\top(x-e)=0$ or $a_j^\top(x+e)=0$. Since both $x-e$ and $x+e$ are nonzero, each of these equalities has conditional probability zero. Thus, almost surely, the solution associated with this choice of $I$ and signs has no additional active index.

% Finally, any point $x\neq\pm e$ with more than $d$ active indices would satisfy the equations above for some choice of $I$ and a mixed sign pattern, together with an additional active equation indexed by some $j\notin I$. There are only finitely many choices of $I$, signs, and $j$. Taking a finite union of the corresponding probability-zero events therefore shows that, almost surely, $|I(x)| \le d$ for every $x \in \RR^d \backslash \{\pm e\}$.

\paragraph{Characterization of criticality.} Recall that $x$ is a critical point if, and only if, there exist coefficients $\lambda_i\in[-1,1]$, $i\in I(x)$, such that
\begin{equation}\label{eq:stationary condition}
	0= \sum_{i\not \in I(x)} \sign((a_i^\top x)^2-(a_i^\top e)^2)(a_i^\top x)a_i +\sum_{i\in I(x)} \lambda_i(a_i^\top x) a_i.
\end{equation}
We will show that, almost surely, every critical point $x\neq\pm e$ satisfies this equation with all active coefficients in $(-1,1)$.
We work on the intersection of $\Epsilon_1$ and the probability-one event that every $x\neq\pm e$ has at most $d$ active indices. For any such point, the vectors $a_i$, $i\in I(x)$, are linearly independent, and $a_i^\top x=\pm a_i^\top e\neq0$ for every $i\in I(x)$. Thus, if $x$ is critical, the coefficients $\lambda_i$ in \eqref{eq:stationary condition} are unique. If $I(x)=\varnothing$, there is nothing to prove. It therefore remains to rule out the existence of a critical point $x\neq\pm e$ for which $\abs{\lambda_k}=1$ for some $k\in I(x)$.
Fix a nonempty index set $I\subset[n]$ with $\abs{I}\le d$, active measurement signs $s\in\{-1,1\}^{I}$, and inactive residual signs $\sigma\in\{-1,1\}^{I^c}$. Define the event
\[
\Epsilon_{I,s,\sigma}:=\left\{
\begin{gathered}
	\text{there is a critical point }x\neq\pm e\text{ with }
	I(x)=I,\quad a_i^\top x=s_i a_i^\top e \text{ for all }i\in I,\\
	\sign\bigl((a_j^\top x)^2-(a_j^\top e)^2\bigr)=\sigma_j \text{ for all }j\notin I,
	\quad\text{and}\quad \abs{\lambda_k}=1\text{ for some }k\in I
\end{gathered}
\right\}.
\]
Since there are only finitely many choices of $I$, $s$, and $\sigma$, it suffices to show that
\(
\PP(\Epsilon_{I,s,\sigma})=0
\).

For each $i\in I$, write $a_i=R_i u_i$, where $R_i=\norm{a_i}$ and $u_i=a_i/\norm{a_i}$. Since the vectors $a_i$ are independent standard Gaussian vectors, the active radii $(R_i)_{i\in I}$ are independent of the active directions $(u_i)_{i\in I}$ and the inactive vectors $(a_j)_{j\notin I}$. Moreover, each $R_i^2$ has a continuous distribution.
The proof proceeds in two steps. First, we show that any critical point with active set $I$, active signs $s$, and inactive signs $\sigma$ must equal a single candidate $\hat{x}$ that is determined without using the active radii $R_{i}$. We then condition on the active directions and inactive vectors and show that $|\lambda_k| =1$ for some $k \in I$ would force one of the continuously distributed radii to take a specific value.

\medskip\noindent\textit{Step 1: critical points are independent of the active radii.} %Determine the candidate critical point without using the active radii.}
The active equations do not depend on the radii and take the form
\[
u_i^\top x=s_i u_i^\top e\qquad\text{for all }i\in I.
\]
 If $\abs{I}=d$, linear independence of the active directions implies that the active equations already determine a unique point $\hat{x}$, which is therefore independent of the radii. For the rest of this step, suppose that $\abs{I}<d$. Choose an orthonormal basis $v_1,\ldots,v_{d-\abs{I}}$ of the orthogonal complement of $\operatorname{span}\{u_{i}\}_{i\in I}$
and let $V=\begin{bmatrix}
	v_1,\ldots, v_{d-\abs{I}}
\end{bmatrix}$.
Also choose a solution $\bar{x}$ of the active equations. Both $V$ and $\bar x$ can be chosen as measurable functions of the active directions and the fixed signs $s$. Every solution of the active equations then has the form
\[
x=\bar{x}+V z\qquad\text{with }z\in\R^{d-\abs{I}}.
\]
If $x$ is a critical point with active set $I$ and signs $(s,\sigma)$, multiplying \eqref{eq:stationary condition} by $V^\top$ eliminates the active terms and leads to
\begin{align}\label{eqn:projected_stationarity}
0=V^\top \left(\sum_{i\not\in I} \sigma_i a_i a_i^\top x + \sum_{i\in I} \lambda_i R_i^2 (u_i^\top x)u_i \right) =V^\top \left(\sum_{i\not\in I} \sigma_i a_i a_i^\top\right) x.
\end{align}
Substituting $x=\bar{x}+Vz$ and setting $H=\sum_{i\not\in I}\sigma_i a_i a_i^\top$ yields \[
V^\top H V z=-V^\top H \bar{x}.
\]
Lemma~\ref{lem:nonsingularity of projected hessian} implies that
\begin{equation*}
	\PP(\det(V^\top H V)=0\mid (u_i)_{i\in I})=0
	\quad\text{almost surely}.
\end{equation*}
Thus, almost surely, any critical point with the specified active set and signs must equal
\[
\hat{x}=\bar{x}-V(V^\top H V)^{-1}V^\top H\bar{x}.
\]
In other words, the set of critical points with active set $I$ and signs $(s,\sigma)$ is contained in $\{\hat{x}\}$, where $\hat{x}$ is a measurable function of $(u_i)_{i\in I}$ and $(a_j)_{j\notin I}$ alone. The candidate $\hat{x}$ need not itself be critical or realize all the specified signs. Even when it does realize them, whether $\hat{x}$ is critical still depends on the active radii through the constraint $\abs{\lambda_i}\le 1$, and Step~2 exploits this dependence.

\medskip\noindent\textit{Step 2: Use the active radii to exclude boundary coefficients.}
When $\abs{I}<d$, the projected stationarity equation~\eqref{eqn:projected_stationarity} shows that $V^\top H\hat{x}=0$, so $H\hat{x}\in\operatorname{span}\{u_i:i\in I\}$. When $\abs{I}=d$, the same inclusion holds because the active directions span $\R^d$. Moreover, the active equations imply that $u_i^\top \hat{x}=s_i u_i^\top e\neq0$ for every $i\in I$. By linear independence, there are therefore unique coefficients $c_i$, $i\in I$, satisfying
\[
\sum_{i\not\in I} \sigma_i a_i a_i^\top \hat{x} + \sum_{i\in I} c_i (u_i^\top \hat{x})u_i=0.
\]
These coefficients are determined entirely by the active directions and inactive vectors. In particular, both $\hat{x}$ and the coefficients $c_i$ are fixed once we condition on $(u_i)_{i\in I}$ and $(a_j)_{j\notin I}$.

Now suppose that a critical point $x$ has active set $I$ and signs $s,\sigma$. By Step~1, we have $x=\hat{x}$. Substituting $a_i=R_i u_i$ into the stationarity condition and comparing the result with the preceding equation, we conclude by uniqueness of the coefficients that $\lambda_i R_i^2=c_i$ for every $i\in I$. Equivalently, the only possible active coefficients are $\lambda_i=c_i/R_i^2$. Since $R_i>0$, we have
\[
\abs{\lambda_i}=1 \iff R_i^2=\abs{c_i}.
\]
Conditional on the active directions and inactive vectors, the value $\abs{c_i}$ is fixed, whereas $R_i^2$ remains continuously distributed. For each $i\in I$, the equality on the right therefore has conditional probability zero. Consequently, almost surely,
\[
\PP(\Epsilon_{I,s,\sigma}\mid (u_i)_{i\in I},(a_i)_{i\notin I})
\le
\sum_{i\in I}
\PP(R_i^2=\abs{c_i}\mid (u_i)_{i\in I},(a_i)_{i\notin I})
=0.
\]
Taking expectations yields $\PP(\Epsilon_{I,s,\sigma})=0$. Finally, taking a finite union over all choices of $I$, $s$, and $\sigma$ shows that, almost surely, every critical point $x\neq\pm e$ satisfies \eqref{eq:stationary condition} with coefficients $\lambda_i\in(-1,1)$ for all $i\in I(x)$.

\paragraph{Characterization of sharp local minima.} We start by proving the contrapositive. Suppose $x\neq\pm e$ is a critical point with $|I(x)|<d$. Choose a nonzero vector $v$ orthogonal to every $a_i$ with $i\in I(x)$. Along $x+tv$, the active residuals remain zero, while the inactive residuals retain their signs for sufficiently small $t$. Thus $t\mapsto\Fn(x+tv)$ is locally smooth, and \eqref{eq:stationary condition} implies that its derivative at zero vanishes. Consequently, $\Fn(x+tv)-\Fn(x)=O(t^2)$, which rules out sharpness.

Conversely, suppose $x$ is critical and $|I(x)|=d$. Set $I=I(x)$, let $\lambda_i\in(-1,1)$ be its active coefficients, and write $b_i(h)=(a_i^\top x)(a_i^\top h)$. Expanding the residuals at $x$ and using \eqref{eq:stationary condition} to cancel the inactive linear terms implies
\[
 \Fn(x+h)-\Fn(x)
 =\frac2n\sum_{i\in I}\bigl(|b_i(h)|-\lambda_i b_i(h)\bigr)
 +O(\norm h^2).
\]
Let $\eta=1-\max_{i\in I}|\lambda_i|>0$. Each summand is at least $\eta|b_i(h)|$. On the event $\Epsilon_1$, the vectors $(a_i^\top x)a_i$, $i\in I$, form a basis, so there is a constant $c_0>0$ such that $\sum_{i\in I}|b_i(h)|\ge c_0\norm h$ for every $h$. Hence
\[
 \Fn(x+h)-\Fn(x)\ge\frac{2\eta c_0}{n}\norm h-C\norm h^2
\]
for some $C\ge0$ and all sufficiently small $h$. Shrinking the neighborhood if necessary yields the required positive linear lower bound and proves sharpness, which completes the proof of Lemma~\ref{lem:critical point}.
\end{proof}

\subsection{Part 2: Existence of active manifolds}

\begin{proposition}\label{prop:active manifold}
	Suppose $n\ge 2d$. Then, almost surely, every critical point
	$x_0 \in \RR^{d} \backslash \{\pm \xs\}$ of $\Fn$ admits the $C^\infty$ active manifold
	\[
	\mathcal{M}_{x_0}
	=
	\bigl\{y\in\R^d:
	a_i^\top(y-x_0)=0
	\text{ for all } i\in I(x_0)\bigr\}.
	\]
	In particular, if $I(x_0)=\varnothing$, then
	$\mathcal{M}_{x_0}=\R^d$.
\end{proposition}

\begin{proof}
By rotational invariance, we may assume without loss of generality that $\xs=e$, where $e$ is the first standard basis vector.
By Lemma~\ref{lem:critical point}, we may work on a probability-one event on which the conditions in~\eqref{eqn:event_eps1} hold and every critical point satisfies the stationarity condition~\eqref{eq:stationary condition} with coefficients $\lambda_i\in(-1,1)$. Fix a critical point $x_0\in \R^d\setminus\{\pm e\}$. We will show that it admits an active manifold. For brevity, write $I=I(x_0)$ and $\mathcal{M}=\mathcal{M}_{x_0}$.

If $I=\varnothing$, then $\Fn$ is $C^\infty$ near $x_0$, so $\R^d$ is an active manifold at $x_0$. For the remainder of the proof, suppose $\abs{I}\ge 1$. The restriction of $\Fn$ to $\mathcal{M}$ is smooth near $x_0$.
It remains to establish sharpness: we must find a neighborhood $U$ of $x_0$ and a constant $c>0$ such that \[ \dist(0,\partial \Fn(x))\ge c \quad \text{for all }x\in U\setminus \mathcal{M}. \]
We first introduce two positive constants. Set
$$\eta = \min_{i \in I}\{1- \abs{\lambda_i}\} >0.$$
Since the vectors $a_i$, $i\in I$, are linearly independent and $\abs{a_i^\top x_0}=\abs{a_i^\top e}>0$ for every $i\in I$, the vectors $(a_i^\top x_0)a_i$, $i\in I$, are also linearly independent. Define
\begin{equation*}
c_0 = \min_{i\in I}\dist\bigl((a_i^\top x_0)a_i,\ \spann\{(a_j^\top x_0)a_j\colon j \in I\setminus\{i\}\}\bigr)>0.
\end{equation*}
We claim that there exists a radius $\delta>0$ such that $\dist(0, \partial \Fn(x))$ is bounded below by a positive constant uniformly over $x \in B_\delta(x_0)\setminus \mathcal{M}$. By continuity, we can choose $\delta$ sufficiently small that, for every $i \in I^c$ and $x \in B_\delta(x_0)$, the sign of $(a_i^\top x)^2-(a_i^\top e)^2$ agrees with its sign at $x_0$. We may also shrink this neighborhood so that, for every $i\in I$, $a_i^\top x$ and $a_i^\top x_0$ have the same sign. Now fix $x \in B_\delta(x_0)\setminus \mathcal{M}$. Since $x\notin \mathcal{M}$, there is an index $i_x \in I$ such that $a_{i_x}^\top x\neq a_{i_x}^\top x_0$. We consider the case
\begin{equation*}
(a_{i_x}^\top x)^2 > (a_{i_x}^\top x_0)^2=(a_{i_x}^\top e)^2.
\end{equation*}
The case with the reverse inequality is analogous. In the displayed case, every element of $\partial \Fn(x)$ can be written as
\begin{equation*}
	\frac{2}{n}\left(\sum_{i \in I\setminus \{i_x\}} \lambda_i' (a_i^\top x)a_i + (a_{i_x}^\top x)a_{i_x}+\sum_{i \in I^c} \sign((a_i^\top x)^2 - (a_i^\top e)^2)(a_i^\top x)a_i\right),
\end{equation*}
where $\lambda_i' \in [-1,1]$. Omitting the factor $\frac{2}{n}$ for the moment, we decompose this expression as follows:
\begin{align*}
		&\sum_{i \in I\setminus \{i_x\}} \lambda_i' (a_i^\top x)a_i + (a_{i_x}^\top x)a_{i_x}+\sum_{i \in I^c} \sign((a_i^\top x)^2 - (a_i^\top e)^2)(a_i^\top x)a_i\\
		&\qquad =\underbrace{\sum_{i \in I\setminus \{i_x\}} \lambda_i' (a_i^\top x_0)a_i + (a_{i_x}^\top x_0)a_{i_x}+\sum_{i \in I^c} \sign((a_i^\top x_0)^2 - (a_i^\top e)^2)(a_i^\top x_0)a_i}_{I_1}\\
		&\qquad \qquad +  \underbrace{\sum_{i \in I\setminus \{i_x\}} \lambda_i' (a_i^\top (x- x_0))a_i +\dotp{a_{i_x},x-x_0}a_{i_x}}_{I_2} + \underbrace{  \sum_{i\in I^c} \sign((a_i^\top x)^2 - (a_i^\top e)^2)(a_i^\top (x- x_0))a_i}_{I_3}.
\end{align*}
By continuity, the norms of $I_2$ and $I_3$ can be made arbitrarily small by choosing $\delta$ sufficiently small. For $I_1$, the stationarity condition at $x_0$ implies that
\begin{align*}
	I_1 & =\sum_{i \in I\setminus \{i_x\}} \lambda_i' (a_i^\top x_0)a_i + (a_{i_x}^\top x_0)a_{i_x}+\sum_{i \in I^c} \sign((a_i^\top x_0)^2 - (a_i^\top e)^2)(a_i^\top x_0)a_i\\
	&= \sum_{i \in I\setminus \{i_x\}} (\lambda_i'-\lambda_i) (a_i^\top x_0)a_i + (1-\lambda_{i_x})(a_{i_x}^\top x_0)a_{i_x}.
\end{align*}
Since the first sum lies in $\spann\{(a_j^\top x_0)a_j \colon j \in I\setminus \{i_x\}\}$, the definitions of $\eta$ and $c_0$ give
\[
\norm{I_1}\ge (1-\lambda_{i_x})\dist\bigl((a_{i_x}^\top x_0)a_{i_x},\ \spann\{(a_j^\top x_0)a_j \colon j \in I\setminus \{i_x\}\}\bigr) \ge \eta c_0 >0.
\]
Consequently, by choosing $\delta$ sufficiently small that $\norm{I_2}+\norm{I_3} \le \frac{\eta c_0}{2}$ and restoring the factor $2/n$, we establish the desired sharpness with $c=\eta c_0/n$.
\end{proof}

\subsection{Part 3: Active strict saddle property}

\begin{proof}[Proof of Theorem~\ref{thm:active strict saddle}]
By rotational invariance, we may assume without loss of generality that $\xs=e$, where $e$ is the first standard basis vector. We work on a probability-one event on which the conclusions of both Lemma~\ref{lem:critical point} and Proposition~\ref{prop:active manifold} hold.
Fix a critical point $x_0\notin \{\pm e\}$. By Proposition~\ref{prop:active manifold}, an active manifold at $x_0$ is
\[
\mathcal{M}_{x_0}=\{y: a_i^\top (y-x_0)=0,i\in I(x_0)\}.
\]
We consider two cases. If $\abs{I(x_0)}=d$, then the last part of Lemma~\ref{lem:critical point} shows that $x_0$ is a sharp, and hence strict, local minimizer of $\Fn$.
 Now, suppose that $\abs{I(x_0)}< d$. For brevity, write $I=I(x_0)$ and $\mathcal{M}=\mathcal{M}_{x_0}$. Points on the manifold can be written as \[
x=x_0+ Vz\qquad\text{with }z\in \R^{d-\abs{I}}, \]
where the columns of $V=\begin{bmatrix} v_1,\ldots, v_{d-\abs{I}}\end{bmatrix}$ form an orthonormal basis of the tangent space of $\mathcal{M}$. Near $x_0$, the restriction of $\Fn$ to $\mathcal{M}$ is quadratic. Define \(
H=\frac{2}{n}\sum_{i\notin I} \sign((a_i^\top x_0)^2-(a_i^\top e)^2) a_i a_i^\top
\).
In these coordinates, the Riemannian Hessian of the restriction at $x_0$ is \(  V^\top H V\).
By Lemma~\ref{lem:nonsingularity of projected hessian}, this Riemannian Hessian is almost surely nonsingular. If it is positive definite, then $x_0$ is a strict local minimizer of $\Fn$ on $\mathcal{M}_{x_0}$, and hence also on $\R^d$ by Lemma~\ref{lem:restricted local minimality}. Otherwise, nonsingularity implies that it has a negative eigenvalue, so $x_0$ is an active strict saddle.

Finally, $\pm e$ are strict local minimizers. Indeed, near $e$, every $a_i^\top x$ has the same sign as $a_i^\top e$. Thus, $\Fn(x)=0$ implies $a_i^\top x=a_i^\top e$ for every $i$, and hence $x=e$ because the measurement vectors span $\R^d$. The same argument applies at $-e$. This proves that $\Fn$ has the active strict saddle property.
\end{proof}

%%% Local Variables:
%%% mode: LaTeX
%%% TeX-master: "../main"
%%% End:

\section{Proof of Theorem~\ref{thm:probability}}\label{app:probability}
We first prove a first-moment result, Theorem~\ref{thm:first moment asymptotic appendix} (Appendix~\ref{sec:first moment asymptotic}), which directly implies part~(\ref{thm:probability:1}) of Theorem~\ref{thm:probability}. Then, together with the second-moment result in Theorem~\ref{thm:second moment asymptotic appendix} (Appendix~\ref{sec:second moment asymptotic}), we prove part~(\ref{thm:probability:2}) of Theorem~\ref{thm:probability} (Appendix~\ref{sec:probability-proof}). We put all supporting tools in Appendices~\ref{app:probability-supporting-tools} and~\ref{app:second moment supporting proofs}, together with explanations of their roles in the proof.
We first fix some common notation in Appendix~\ref{sec:common notation}, which is used consistently in Appendices~\ref{app:probability}, \ref{app:probability-supporting-tools}, and~\ref{app:second moment supporting proofs}.

\subsection{Common notation}\label{sec:common notation}

Throughout Appendices~\ref{app:probability}, \ref{app:probability-supporting-tools}, and~\ref{app:second moment supporting proofs}, fix an integer $d\ge4$ and a deterministic unit vector $\xs\in\R^d$.
We start by defining the notation we will use in the proof. Most of the intuition for what we are about to introduce already appeared in Section~\ref{sec:proof-strategy}. We adopt the following convention throughout. All our asymptotic symbols of the form $O(\cdot)$ and $o(\cdot)$ refer to limits $n\to\infty$. Sometimes we will add subscripts on $O(\cdot)$ or $o(\cdot)$ to highlight the parameters that are held fixed and on which the implicit constants may depend.

We use the following constant \begin{equation}\label{eq:def-kappa}\kappa=8\tau/[\pi(1+\tau^2)^2]>0,\end{equation}
where again $\tau$ is the unique solution of the equation $\arctan \tau + \tau/(1+\tau^2) = \pi/4.$ Throughout, we fix the following active sign pattern, \begin{equation}\label{eq:sign pattern}
s=(s_1,\ldots,s_d)=(\underbrace{1,\ldots,1}_{d-2},-1,-1).
\end{equation}
For any realization of $(a_1,\ldots,a_n)$, let $A$ denote the matrix whose $j$th-row corresponds to the vector $a_j$. For any index set $J=\{j_1<\cdots<j_r\} \subset [n]$, define the submatrix indexed by $J$ as \[
A_J=\begin{pmatrix}
a_{j_1} & \cdots & a_{j_r}
\end{pmatrix}^\top.
\]
For any $x,y\in \R^d$ and a random vector $a\sim N(0,I_d)$, define the objects
\begin{equation}\label{eq:def Z}
Z_x(a):=\sign((a^\top x)^2-(a^\top \xs)^2)(a^\top x)a,
\end{equation}
\begin{equation}\label{eq:def mean cov}
\mu(x):=\mathbb E[Z_x(a)], \quad \Sigma_x:=\operatorname{Cov}(Z_x(a)), \quad\text{and}\quad \Sigma_{x,y}:=\operatorname{Cov}\begin{pmatrix}Z_x(a)\\ Z_y(a)\end{pmatrix}.
\end{equation}
For any pair of matrices $B=\begin{pmatrix} b_1& \cdots & b_d \end{pmatrix}^\top\in\R^{d\times d}$ and $C=\begin{pmatrix}  c_1 &\cdots & c_r \end{pmatrix}^\top\in\R^{r\times d}$, define
\begin{equation}\label{eq:def Q S}
Q(x;B):=\begin{pmatrix} (b_1^\top x)b_1&\cdots &(b_d^\top x)b_d \end{pmatrix},\qquad\text{and}\qquad
S(x;C):=\sum_{i=1}^{r}Z_x(c_i).
\end{equation}
For any index set $I=\{i_1<\cdots<i_d\}\subset [n]$, define
\begin{equation}\label{eq:def X_I}
X_I(x):=\begin{pmatrix}a_{i_1}^{\top}(x-s_1\xs)\\ \vdots\\ a_{i_d}^{\top}(x-s_d\xs)\end{pmatrix}=A_Ix-\diag(s)A_I\xs.
\end{equation}
Almost surely, $A_I$ is invertible, and the solution of $X_I(x)=0$ is denoted by
\begin{equation}\label{eq:def x_I}
x_I=A_I^{-1}\diag(s)A_I \xs.
\end{equation}
By Lemma~\ref{lem:critical point strategy}, almost surely, $x$ is a sharp spurious local minimum with active index set $I$ and active sign pattern $s$ if, and only if, the following holds:
\[ 0\in S(x;A_{I^c})+Q(x;A_I)(-1,1)^d,\qquad\text{and}\qquad X_I(x)=0. \]

Recall the saddle radius $\tau$ and the population saddle sphere $\mathcal S$ from Proposition~\ref{prop:population landscape}, and the tube $T_n$ from~\eqref{eq:tube}. We write $\mathbb S^{d-2}=\{u\colon \norm{u}=1,\ u\perp\xs\}$, so that $\mathcal S=\tau\,\mathbb S^{d-2}$. When $n^{-1/2}<\tau$, every $x\in T_n$ can be written uniquely as $x=\alpha_x\xs+(\tau+\beta_x)u_x$ with $u_x\in\mathbb S^{d-2}$ and $\abs{\alpha_x},\abs{\beta_x}\le n^{-1/2}$; we call $(\alpha_x,\beta_x,u_x)$ the tube coordinates of $x$.
We will restrict to well-conditioned active matrices, using the following lemma. We defer its proof to Appendix~\ref{app:proof regularity constants}.
\begin{lemma}[Choice of regularity constants]\label{lem:regularity constants}
Fix a dimension $d\ge4$. Then, there exist finite deterministic constants $M_d,\lambda_d,\vartheta_d>0$, depending only on $d$, with the following properties. For $x\in\R^d$, define the set
\[
\mathcal{K}_x=\left\{B\in\R^{d\times d} : \max_{1\le j\le d}\norm{b_j}\le M_d,\ \abs{\det B}\ge\lambda_d,\ \text{and}\ \abs{\det Q(x;B)}\ge\lambda_d\right\},
\]
and, for $n\ge d$ and $I\subset[n]$ with $\abs{I}=d$, the conditional expectation
\[
a(x)=\E\bigl[\abs{\det A_I}\abs{\det Q(x;A_I)}1_{\mathcal{K}_x}(A_I)\bigm|X_I(x)=0\bigr].
\]
The function $a(\cdot)$ is independent of the choice of $I$, is continuous near $\mathcal S$, and satisfies
\[
a(\tau u)=\vartheta_d\qquad\text{for every }u\in\mathbb S^{d-2}.
\]
Moreover, for any fixed $d$,\[
\sup_{x\in T_n}\abs{a(x)-\vartheta_d}\to0\qquad\text{as }n\to\infty.
\]
\end{lemma}
We fix the constants $M_d$, $\lambda_d$, and $\vartheta_d$ from Lemma~\ref{lem:regularity constants} for the rest of the proof.
We count only those sharp spurious local minimizers that also lie in the tube $T_n$ around the population saddle sphere and whose active matrix satisfies $A_I\in\mathcal K_{x_I}$, with $\mathcal K_x$ as in Lemma~\ref{lem:regularity constants}. Formally, for a set $I\subset [n]$ with $\abs{I}=d$, define the random indicator function
\begin{equation}\label{eq:def N}
\begin{aligned}
\chi_I^n(a_1,\ldots,a_n) &= 1_{\{x_I\in T_n\}} 1_{\mathcal{K}_{x_I}}(A_I) 1_{\{-S(x_I;A_{I^c})\in Q(x_I;A_I)(-1,1)^d\}},\quad\text{and}\\
N_{d,n}(a_1,\ldots,a_n)&=\sum_{\abs{I}=d} \chi_I^n.
\end{aligned}
\end{equation}
Then $\chi_I^n=1$ precisely when the candidate solution $x_I$ is a sharp spurious local minimum lying in $T_n$ and satisfying $A_I\in\mathcal{K}_{x_I}$. Thus $N_{d,n}$ counts the total number of such sharp spurious local minima for a given realization.

%%%%%%%%%%%%%%%%%%%%%%%%% First moment
\subsection{First moment characterization}\label{sec:first moment asymptotic}\label{pf:first moment asymptotic}
For the next result, recall that $M_d$, $\lambda_d$, and $\vartheta_d$ are the constants provided by Lemma~\ref{lem:regularity constants}, and $\Sigma_x$ and $N_{d,n}$ were defined in~\eqref{eq:def mean cov} and~\eqref{eq:def N}, respectively.
\begin{thm}\label{thm:first moment asymptotic appendix}
For any $d \geq 4,$ we have that
\[ n^{-d/2+1} \E[N_{d,n}] \rightarrow \gamma_d \mathcal{I}_{d} \quad \text{as }n\to\infty \]
with
\[
\gamma_d=\frac{2^d\tau^{d-2}\vartheta_d}{d! [2\pi(1+\tau^2)]^{d/2}}, \qquad\text{and}\qquad \mathcal{I}_d = \int_{\mathbb{S}^{d-2}}\int_{[-1,1]^2} \varphi_{\scriptscriptstyle \Sigma_{\tau u}}(-\kappa z_1 \xs+\kappa z_2 u) \,\dd z \dd \sigma(u),
\]
where $\dd\sigma$ is the surface measure on $\mathbb S^{d-2}$, and $\varphi_\Sigma$ is the density of $N(0,\Sigma)$.
\end{thm}
\begin{proof}
 By rotational invariance, we may assume without loss of generality that $\xs=e$, where $e$ is the first standard basis vector.
Recall that $N_{d,n}=\sum_{\abs{I}=d}\chi_I^n$. Since the Gaussian rows are i.i.d.\ and every ordered active set uses the same sign pattern, all the indicators $\chi_I^n$ have the same expectation. Hence, for any fixed $I\subset[n]$ with $\abs{I}=d$,
\[
\E[N_{d,n}]=\sum_{\abs{J}=d}\E[\chi_J^n]=\binom{n}{d}\E[\chi_I^n].
\]
Since $\binom{n}{d}=n^d/d!\,(1+O_d(n^{-1}))$, it suffices to prove
\begin{equation}\label{eq:first moment for one active set}
\E[\chi_I^n]=d!\,\gamma_d\mathcal I_d\,n^{-d/2-1}+o_d(n^{-d/2-1}).
\end{equation}
In what follows we decompose the expectation of $\chi_I^n$ into a leading term and a negligible remainder, and then evaluate the leading term in tube coordinates. To do so, we leverage the Kac-Rice formula (Lemma~\ref{lem:Kac-Rice for first and second moment}(\ref{lem:Kac-Rice for first and second moment:I})), which gives
\begin{align*}
\E[\chi_I^n]&=\int_{T_n} p_{X_I(x)}(0)\E\bigl[\abs{\det A_I}1_{\mathcal{K}_x}(A_I) 1_{-S(x;A_{I^c})\in Q(x;A_I)(-1,1)^d}\bigm|X_I(x)=0 \bigr] \dd x\\
&=\int_{T_n} p_{X_I(x)}(0)\E\bigl[\abs{\det A_I} I_1\bigm|X_I(x)=0 \bigr] \dd x,
\end{align*}
where $p_{X_I(x)}$ is the density of the $\R^d$ Gaussian vector $X_I(x)$, and
\[
I_1=\E\left[1_{\mathcal{K}_x}(A_I)1_{-S(x;A_{I^c})\in Q(x;A_I)(-1,1)^d}\middle|A_I,X_I(x)=0\right].
\]
We use a local central limit theorem, Lemma~\ref{lem:limit theorem for first and second moment}(\ref{lem:limit theorem for first and second moment:I}), to approximate \begin{equation}\label{eq:limit theorem for first moment}
 I_1 =\underbrace{\frac{2^d}{(n-d)^{d/2}} \abs{\det Q(x;A_I)} 1_{\mathcal{K}_x}(A_I) \varphi_{\scriptscriptstyle\Sigma_x}(-\sqrt{n-d}\mu(x))}_{I_2}+\underbrace{\frac{1}{(n-d)^{d/2}} 1_{\mathcal{K}_x}(A_I)r_n(x;A_I)}_{I_3},
\end{equation}
where $\varphi_{\scriptscriptstyle\Sigma}$ is the density of $N(0,\Sigma)$. To justify this approximation, note that $X_I(x)$ depends only on $A_I$. Since $A_{I^c}$ is independent of $A_I$, its rows remain independent standard Gaussian vectors under the conditioning in $I_1$.
The remainder $I_3$ is negligible, in particular, Lemma~\ref{lem:limit theorem for first and second moment}(\ref{lem:limit theorem for first and second moment:I}) ensures that the function $r_n$ satisfies the uniform bound
\begin{equation}\label{eq:remainder first moment}
\sup_{x\in T_n}\sup_{A\in \R^{d\times d}} 1_{\mathcal{K}_x}(A) \abs{r_n(x;A)} \le \eta_n
\end{equation}
for some deterministic sequence $\eta_n\to0$, so that $\abs{I_3}\le(n-d)^{-d/2}\eta_n$ for all $x\in T_n$.

Substituting $I_1=I_2+I_3$ into the integral above, we split the contribution of this active set into the leading term $I_4$ and the remainder $I_5$:
\ificlr
\begin{align*}
\E[\chi_I^n]
&= \underbrace{\int_{T_n}p_{X_I(x)}(0) \E\bigl[\abs{\det A_I} I_2\bigm|X_I(x)=0\bigr] \dd x }_{I_4} \\
&\quad+ \underbrace{\int_{T_n}p_{X_I(x)}(0) \E\bigl[\abs{\det A_I} I_3\bigm|X_I(x)=0\bigr] \dd x }_{I_5}.
\end{align*}
\else
\begin{align*}
\E[\chi_I^n]
&= \underbrace{\int_{T_n}p_{X_I(x)}(0) \E\bigl[\abs{\det A_I} I_2\bigm|X_I(x)=0\bigr] \dd x }_{I_4}+ \underbrace{\int_{T_n}p_{X_I(x)}(0) \E\bigl[\abs{\det A_I} I_3\bigm|X_I(x)=0\bigr] \dd x }_{I_5}.
\end{align*}
\fi
To control the remainder term $I_5$, we use \eqref{eq:remainder first moment} together with the following bounds, valid for $x\in T_n$ and $A_I\in\mathcal{K}_x$:
\begin{gather*}
\abs{\det A_I}\le \prod_{i\in I} \norm{a_i} \le M_d^d,\qquad \abs{\det Q(x;A_I)}\le \prod_{i\in I} \norm{(a_i^\top x)a_i} =O_d(1),\\
\operatorname{Vol}(T_n)=O_d(n^{-1}),\qquad\text{and}\qquad p_{X_I(x)}(0)=(2\pi)^{-d/2}\prod_{i=1}^d \frac{1}{\norm{x-s_i e}} =O_d(1).
\end{gather*}
The first three bounds follow from our construction. To derive the last relation, i.e., the density formula, write $I=\{i_1<\cdots<i_d\}$. For fixed $x$, the coordinates $a_{i_j}^\top(x-s_j e)$ of $X_I(x)$ are independent $N(0,\norm{x-s_j e}^2)$ variables, so their densities at zero, $(2\pi)^{-1/2}\norm{x-s_j e}^{-1}$, multiply to produce the displayed expression. Moreover, on $T_n$, we have $\norm{x-s_j e}\ge\tau-n^{-1/2}\ge\tau/2$ for sufficiently large $n$, which implies the uniform $O_d(1)$ bound.
Since $(n-d)/n\to1$ and $\eta_n\to0$, these bounds imply
\begin{equation}\label{eq:remainder integral first moment}
I_5 = o_d(n^{-d/2-1}).
\end{equation}

It remains to evaluate the leading term
\[I_4= 2^d (n-d)^{-d/2} \int_{T_n}p_{X_I(x)}(0) a(x) \varphi_{\scriptscriptstyle\Sigma_x}(-\sqrt{n-d}\mu(x)) \dd x.
\] % where \[ a(x)=\E\left[\abs{\det A_I}\abs{\det Q(x;A_I)}1_{\mathcal{K}_x}(A_I)\;\middle|\; X_I(x)=0 \right]. \]
By Lemma~\ref{lem:regularity constants}, $a(\tau u)=\vartheta_d>0$ for every $u\in\mathbb S^{d-2}$. In what follows, we do a change of variables to integrate over the tube. For sufficiently large $n$ such that $n^{-1/2}<\tau$, parametrize $T_n$ by
\[
x=x_{n,z,u}=\frac{z_1}{\sqrt n}e+\left(\tau+\frac{z_2}{\sqrt n}\right)u\qquad\text{with } z\in[-1,1]^2,\quad\text{and}\quad u\in\mathbb S^{d-2}.
\]
Hence
\(
\dd x=\frac{1}{n}\left(\tau+\frac{z_2}{\sqrt n}\right)^{d-2}\dd z\,\dd\sigma(u),
\)
where $\dd\sigma$ is surface measure on $\mathbb S^{d-2}$. Therefore,
\ificlr
\begin{align*}
I_4 &= \frac{2^d}{n(n-d)^{d/2}}
\int_{\mathbb S^{d-2}}\int_{[-1,1]^2}
p_{X_I(x_{n,z,u})}(0)\,a(x_{n,z,u})\\
&\qquad\times
\varphi_{\scriptscriptstyle\Sigma_{x_{n,z,u}}}\!\bigl(-\sqrt{n-d}\,\mu(x_{n,z,u})\bigr)
\bigl(\tau+z_2/\sqrt n\bigr)^{d-2}\,\dd z\,\dd\sigma(u).
\end{align*}
\else
\[
I_4 = \int_{\mathbb S^{d-2}}\int_{[-1,1]^2}
\frac{2^d\bigl(\tau+z_2/\sqrt n\bigr)^{d-2}}{n(n-d)^{d/2}}\,
p_{X_I(x_{n,z,u})}(0)\,a(x_{n,z,u})\,
\varphi_{\scriptscriptstyle\Sigma_{x_{n,z,u}}}\!\bigl(-\sqrt{n-d}\,\mu(x_{n,z,u})\bigr)
\,\dd z\,\dd\sigma(u).
\]
\fi
We now let $n\to\infty$. Since $e\perp u$ and $\norm{e}=\norm{u}=1$, a direct calculation shows that
\[
\norm{x_{n,z,u}-s_i e}^2
=1+\tau^2+\frac{2(\tau z_2-s_i z_1)}{\sqrt n}+\frac{z_1^2+z_2^2}{n}
\to1+\tau^2 \qquad \text{for all }i\in\{1,\ldots,d\}
\]
uniformly for $z\in[-1,1]^2$ and $u\in\mathbb S^{d-2}$. Substituting this into the density yields the uniform limit
\[
p_{X_I(x_{n,z,u})}(0)
=(2\pi)^{-d/2}\prod_{i=1}^d\frac{1}{\norm{x_{n,z,u}-s_i e}}
\to[2\pi(1+\tau^2)]^{-d/2}.
\]
Lemma~\ref{lem:regularity constants} and Lemma~\ref{lem:tube} provide the remaining limits, namely
\begin{align*}
&a(x_{n,z,u})\to\vartheta_d,\quad
\sqrt{n-d}\mu(x_{n,z,u})\to-\kappa z_1e+\kappa z_2u,\quad \text{and}\quad\norm{\Sigma_{x_{n,z,u}}-\Sigma_{\tau u}}_{\operatorname{op}}\to0,
\end{align*}
uniformly for $z\in[-1,1]^2$ and $u\in\mathbb S^{d-2}$. Since $\Sigma_{\tau u}=\tau^2(I_d+2uu^\top)\succeq\tau^2I_d$, the Gaussian densities also converge uniformly. The radial factor converges uniformly to $\tau^{d-2}$. Passing to the limit on this fixed integration domain, using $(n-d)/n\to1$ and the symmetry $\varphi_\Sigma(-v)=\varphi_\Sigma(v)$, we obtain
\[
I_4=\frac{2^d\tau^{d-2}\vartheta_d}{[2\pi(1+\tau^2)]^{d/2}}n^{-d/2-1}\underbrace{\int_{\mathbb S^{d-2}}\int_{[-1,1]^2}\varphi_{\scriptscriptstyle\Sigma_{\tau u}}(-\kappa z_1e+\kappa z_2u)\dd z\,\dd\sigma(u)}_{\mathcal I_d}+o_d(n^{-d/2-1}).
\]
Combining this asymptotic formula for $I_4$ with the bound for $I_5$ proves \eqref{eq:first moment for one active set}, completing the proof.
\end{proof}

%%%%%% Second moment
\subsection{Second moment characterization}\label{sec:second moment asymptotic}\label{pf:second moment asymptotic}
For the following result and its proof, recall that $M_d$, $\lambda_d$, and $\vartheta_d$ are the constants provided by Lemma~\ref{lem:regularity constants}, and the objects $\kappa, \Sigma_{x,y}$ and $N_{d,n}$ were defined in~\eqref{eq:def-kappa}, \eqref{eq:def mean cov} and~\eqref{eq:def N}, respectively.

\begin{thm}\label{thm:second moment asymptotic appendix}
For any $d \geq 4,$ we have that
\[ n^{-d+2} \E[N_{d,n}^2] \rightarrow \gamma_d^2 \mathcal{J}_{d},\]
where \[
\begin{aligned}
\gamma_d&=\frac{2^d\tau^{d-2}\vartheta_d}{d! [2\pi(1+\tau^2)]^{d/2}}, \quad \text{and}\\
\mathcal{J}_d&= \int_{\substack{u,v\in \mathbb{S}^{d-2}\\u\neq\pm v}} \int_{[-1,1]^4} \varphi_{\scriptscriptstyle \Sigma_{\tau u,\tau v}} \left(\begin{pmatrix}-\kappa z_1\xs+\kappa z_2u\\ -\kappa w_1\xs+\kappa w_2v\end{pmatrix}\right)\, \dd z\, \dd w \,\dd \sigma(u)\, \dd \sigma(v).
\end{aligned}\]
\end{thm}
\begin{proof}
By rotational invariance of the standard Gaussian measurements, we may assume without loss of generality that $\xs=e$.
Recall that $N_{d,n}=\sum_{\abs I=d}\chi_I^n$ and that $(\chi_I^n)^2=\chi_I^n$. Expanding the square and grouping the off-diagonal terms by the overlap $r=\abs{I\cap J}$, we obtain
\[
\E[N_{d,n}^2]=\E[N_{d,n}]+\sum_{r=0}^{d-1}D_{r,n}\qquad\text{where}\qquad
D_{r,n}=\sum_{\substack{I,J\subset[n]\\\abs I=\abs J=d,\ \abs{I\cap J}=r}}\E[\chi_I^n\chi_J^n].
\]
The sum defining $D_{r,n}$ is over ordered pairs, so no factor of two is needed. By Theorem~\ref{thm:first moment asymptotic appendix}, we have $\E[N_{d,n}]=O_d(n^{d/2-1})=o_d(n^{d-2})$. It therefore suffices to establish $\mathcal J_d<\infty$ and the estimates
\begin{equation}\label{eq:second moment overlap estimates}
D_{0,n}=\gamma_d^2\mathcal J_d\,n^{d-2}+o_d(n^{d-2}), \qquad\text{and}\qquad D_{r,n}=o_d(n^{d-2})\qquad\text{for }1\le r\le d-1.
\end{equation}

The rest of the proof proceeds in three steps. We first express each $\E[\chi_I^n\chi_J^n]$ as a conditional Kac--Rice integral. We then split the integration domain into three regions and use the estimates from Appendix~\ref{app:second moment supporting proofs} to control the integral in each region. Finally, we combine these estimates.

\paragraph{Step 1: Express each pair contribution as a conditional integral.}
Fix distinct index sets $I,J\subset[n]$ with $\abs I=\abs J=d$ and $\abs{I\cap J}=r$. Let $U=(I\cup J)^c$ and $\nu=\abs U=n-2d+r$. The sums over the two sets of inactive measurements decompose as
\[
S(x;A_{I^c})=S(x;A_{J\setminus I})+S(x;A_U), \qquad\text{and}\qquad S(y;A_{J^c})=S(y;A_{I\setminus J})+S(y;A_U).
\]
Conditioning on $A_{I\cup J}$ fixes the first term in each decomposition. The rows of $A_U$ remain independent standard Gaussian vectors, since they are independent of $A_{I\cup J}$ and the constraints $X_I(x)=X_J(y)=0$ involve only rows in $I\cup J$. The sums $S(x;A_U)$ and $S(y;A_U)$ depend on the same Gaussian rows and are generally dependent, so the two inclusion events must be treated jointly.

For a deterministic matrix $B\in\R^{(2d-r)\times d}$ whose rows are indexed by $I\cup J$, let $C=(g_1,\ldots,g_\nu)^\top$, where the $g_j$ are independent $N(0,I_d)$ vectors, and define
\[
\pi_{n,I,J}(x,y;B)
=\PP_{\scriptscriptstyle g_1,\ldots,g_\nu}\left\{
\begin{aligned}
-S(x;C)&\in Q(x;B_I)(-1,1)^d+S(x;B_{J\setminus I}),\\
-S(y;C)&\in Q(y;B_J)(-1,1)^d+S(y;B_{I\setminus J})
\end{aligned}
\right\}.
\]
This is the probability appearing in Lemma~\ref{lem:limit theorem for first and second moment}(\ref{lem:limit theorem for first and second moment:II}). Thus, Lemma~\ref{lem:Kac-Rice for first and second moment}(\ref{lem:Kac-Rice for first and second moment:II}), followed by conditioning on $A_{I\cup J}$, implies
\begin{equation}\label{eq:second moment pair integral}
% Scale uniformly to keep the entire integral on one line in both formats.
\resizebox{\linewidth}{!}{$\displaystyle
\begin{aligned}
&\E[\chi_I^n\chi_J^n]\\
&\quad=\int\limits_{T_n^2\cap\mathcal R_{I,J}}p_{X_I(x),X_J(y)}(0,0)\,\E\bigl[\abs{\det A_I}\abs{\det A_J}1_{\mathcal K_x}(A_I)1_{\mathcal K_y}(A_J)\pi_{n,I,J}(x,y;A_{I\cup J})\bigm|X_I(x)=X_J(y)=0\bigr]\dd x\,\dd y,
\end{aligned}
$}
\end{equation}
where
\(
\mathcal R_{I,J}=\left\{(x,y)\colon\operatorname{Cov}\begin{pmatrix}X_I(x)\\X_J(y)\end{pmatrix}\succ0\right\},
\)
and $p_{X_I(x),X_J(y)}$ is the density of the Gaussian vector $\begin{pmatrix}X_I(x)\\X_J(y)\end{pmatrix}\in\R^{2d}$ on this region.

\paragraph{Step 2: Split the integral into three regions.}
For all sufficiently large $n$ such that $n^{-1/2}\le\tau/2$, every $x\in T_n$ has a unique representation
\[
x=\alpha_x e+(\tau+\beta_x)u_x\qquad\text{with }
\abs{\alpha_x},\abs{\beta_x}\le n^{-1/2},\qquad\text{and}\qquad u_x\in\mathbb S^{d-2}.
\]
For a pair $(x,y)$, define
\[
\delta(x,y)=\min\{\norm{u_x-u_y},\norm{u_x+u_y}\}, \;\text{and}\; \rho(x,y)=\sqrt{\bigl(\alpha_x-\sign(u_x^\top u_y)\alpha_y\bigr)^2+(\beta_x-\beta_y)^2}.
\]
The quantity $\delta$ measures angular separation up to sign. When $\delta<1$, the quantity $\rho$ measures how much the two points differ in their normal coordinates, after replacing $y$ by $-y$ if $u_y$ is closer to $-u_x$ than to $u_x$.

Fix $0<\delta_0<1$. Pairs with $\delta\ge\delta_0$ form the separated region, where the local limit estimate in Lemma~\ref{lem:limit theorem for first and second moment}(\ref{lem:limit theorem for first and second moment:II}) applies. The limiting covariance $\Sigma_{\tau u_x,\tau u_y}$ becomes singular at $u_x=\pm u_y$, so pairs with $\delta<\delta_0$ require separate estimates. We divide these nearby pairs by comparing their angular separation $\delta$ with their normal separation $\rho$. Formally, we define three regions
\[
\begin{aligned}
\Omega^{\mathrm{sep}}&=\{(x,y)\in T_n^2\cap\mathcal R_{I,J}:\delta(x,y)\ge\delta_0\},\\
\Omega^{\mathrm{tan}}&=\{(x,y)\in T_n^2\cap\mathcal R_{I,J}:\delta(x,y)<\delta_0,\ \rho(x,y)\le\delta(x,y)\},\quad\text{and}\\
\Omega^{\mathrm{nor}}&=\{(x,y)\in T_n^2\cap\mathcal R_{I,J}:\delta(x,y)<\delta_0,\ \rho(x,y)>\delta(x,y)\}.
\end{aligned}
\]
We call them the separated, tangential, and normal-dominated regions, respectively. As $u$ varies, the base point $\tau u$ moves along the saddle sphere $\mathcal S$, while $\alpha$ and $\beta$ describe displacements normal to it. In the tangential region, angular separation is at least as large as separation in the normal coordinates; in the normal-dominated region, the latter is larger.

For a fixed pair $I,J$ with overlap $r$, let $I_{r,n}^{\mathrm{sep}}$ denote the integral in~\eqref{eq:second moment pair integral} restricted to $\Omega^{\mathrm{sep}}$. Define $D_{r,n}^{\mathrm{sep}}$ by summing $I_{r,n}^{\mathrm{sep}}$ over all ordered pairs $I,J\subset[n]$ with $\abs I=\abs J=d$ and $\abs{I\cap J}=r$. Define $D_{r,n}^{\mathrm{tan}}$ and $D_{r,n}^{\mathrm{nor}}$ in the same way, using the restrictions of~\eqref{eq:second moment pair integral} to $\Omega^{\mathrm{tan}}$ and $\Omega^{\mathrm{nor}}$, respectively. We suppress the dependence of the regions and of $I_{r,n}^{\mathrm{sep}}$ on $I,J$, as well as the dependence of the regions, integrals, and sums on $\delta_0$.
Since the three regions partition the integration domain, we have
\[
D_{r,n}=D_{r,n}^{\mathrm{sep}}+D_{r,n}^{\mathrm{tan}}+D_{r,n}^{\mathrm{nor}}.
\]

Let $\mathcal J_d(\delta_0)$ denote the integral defining $\mathcal J_d$ in the theorem, with the additional restriction $\min\{\norm{u-v},\norm{u+v}\}\ge\delta_0$.
For each fixed $\delta_0$, Lemma~\ref{lem:tube} implies $\lambda_{\min}(\Sigma_{\tau u,\tau v})\ge c_{d,\delta_0}>0$ uniformly over the restricted angular domain. Since $\Sigma_{\tau u,\tau v}$ is a $2d\times2d$ covariance matrix, the Gaussian density satisfies
\[
\varphi_{\scriptscriptstyle\Sigma_{\tau u,\tau v}}(\xi)
\le\frac{1}{(2\pi)^d\sqrt{\det\Sigma_{\tau u,\tau v}}}
\le\frac{1}{(2\pi c_{d,\delta_0})^d}
\qquad\text{for every }\xi\in\R^{2d}.
\]
Thus the integrand defining $\mathcal J_d(\delta_0)$ is uniformly bounded. Since $[-1,1]^4$ has finite volume and $\mathbb S^{d-2}$ has finite surface measure, $\mathcal J_d(\delta_0)<\infty$.

Table~\ref{table:regional estimates} summarizes the estimates, proved in Appendix~\ref{app:second moment supporting proofs}, for the three regional terms. All asymptotics in the table are as $n\to\infty$ with $\delta_0$ fixed. In particular, the separated region supplies the leading term for disjoint pairs, while the tangential error can be made small by subsequently decreasing $\delta_0$.
\begin{table}[!htbp]
\centering
\renewcommand{\arraystretch}{1.5}
\begin{tabular}{@{}lll@{}}
\toprule
Term & Asymptotic estimate & Reference \\
\midrule
$D_{r,n}^{\mathrm{sep}}$ & $\begin{array}{ll}r=0:&\gamma_d^2\mathcal J_d(\delta_0)n^{d-2}+o_{d,\delta_0}(n^{d-2}),\\1\le r\le d-1:&o_{d,\delta_0}(n^{d-2})\end{array}$ & Lemma~\ref{lem:separated region} \\
\midrule
$D_{r,n}^{\mathrm{tan}}$ & $\begin{array}{ll}r=0:&\le C_dn^{d-2}\delta_0^{d/2-1}+C_dn^{d/2-1},\\1\le r\le d-1:&o_{d,\delta_0}(n^{d-2})\end{array}$ & Lemma~\ref{lem:tangential region} \\
\midrule
$D_{r,n}^{\mathrm{nor}}$ & $0\le r\le d-1:\quad o_{d,\delta_0}(n^{d-2})$ & Lemma~\ref{lem:normal dominated region} \\
\bottomrule
\end{tabular}
\caption{Estimates for the three regional terms at fixed $\delta_0$. The constant $C_d>0$ depends only on $d$, and $\gamma_d$ is as in Theorem~\ref{thm:first moment asymptotic appendix}.}\label{table:regional estimates}
\end{table}

\paragraph{Step 3: Combine the estimates and remove the angular cutoff.}
For $1\le r\le d-1$, fix any $\delta_0\in(0,1)$. All three terms in Table~\ref{table:regional estimates} are $o_{d,\delta_0}(n^{d-2})$, which proves the second estimate in~\eqref{eq:second moment overlap estimates}.

It remains to treat $r=0$. Fix $\delta_0\in(0,1)$ and multiply the estimates in Table~\ref{table:regional estimates} by $n^{2-d}$. As $n\to\infty$ with this cutoff fixed, we have
\[
n^{2-d}D_{0,n}^{\mathrm{sep}}\longrightarrow\gamma_d^2\mathcal J_d(\delta_0), \quad 0\le n^{2-d}D_{0,n}^{\mathrm{tan}}\le C_d\delta_0^{d/2-1}+C_dn^{1-d/2}, \quad\text{and}\quad n^{2-d}D_{0,n}^{\mathrm{nor}}\longrightarrow0.
\]
Since the regional terms are nonnegative and $n^{1-d/2}\to0$, we have
\begin{equation}\label{eq:delta0}
\gamma_d^2\mathcal J_d(\delta_0) \le\liminf_{n\to\infty}n^{2-d}D_{0,n} \le\limsup_{n\to\infty}n^{2-d}D_{0,n} \le\gamma_d^2\mathcal J_d(\delta_0)+C_d\delta_0^{d/2-1}.
\end{equation}
The total $D_{0,n}$ is independent of $\delta_0$. Taking $\delta_0=1/2$ in the upper bound of~\eqref{eq:delta0}, and using $\mathcal J_d(1/2)<\infty$, shows that
\[
L:=\limsup_{n\to\infty}n^{2-d}D_{0,n}<\infty.
\]
The lower bound in~\eqref{eq:delta0} then implies $\gamma_d^2\mathcal J_d(\delta_0)\le L$ for every $\delta_0\in(0,1)$. As $\delta_0\downarrow0$, the angular integration domains increase to $\{(u,v):u\ne\pm v\}$. Since the integrand is nonnegative, monotone convergence yields
\[
\lim_{\delta_0\downarrow0}\mathcal J_d(\delta_0)=\mathcal J_d\le L/\gamma_d^2<\infty.
\]
Finally, let $\delta_0\downarrow0$ in~\eqref{eq:delta0}. Since $\mathcal J_d(\delta_0)\uparrow\mathcal J_d$ and $C_d\delta_0^{d/2-1}\to0$, both the lower and upper bounds tend to $\gamma_d^2\mathcal J_d$. This proves the first estimate in~\eqref{eq:second moment overlap estimates} and completes the proof.
\end{proof}

%%%%%% compute ratio
\subsection{Constant-probability lower bound}\label{sec:probability-proof}\label{pf:probability}
For the following result, recall that $M_d$ and $\lambda_d$ are the constants provided by Lemma~\ref{lem:regularity constants}, and $N_{d,n}$ was defined in~\eqref{eq:def N}.
We show via numerical computation in Appendix~\ref{app:numerics} that we can take $c_\star=0.136$ in the following theorem.
\begin{thm}\label{thm:probability appendix}
For any $d\ge4$, we have
\[
\liminf_{n\to\infty}\PP(N_{d,n}>0)\ge c_{\star},
\]
where $c_{\star}>0$ is a universal constant independent of $d$.
\end{thm}
\begin{proof}
By rotational invariance, we may assume without loss of generality that $\xs=e$, where $e$ is the first standard basis vector.
The Paley--Zygmund inequality implies that \[
\PP(N_{d,n}>0)\ge \frac{(\E N_{d,n})^2}{\E N_{d,n}^2}.
\]
Theorems~\ref{thm:first moment asymptotic appendix}  and~\ref{thm:second moment asymptotic appendix} then imply, at fixed $d$,\[
\liminf_{n\to\infty} \PP(N_{d,n}>0)\ge \frac{\mathcal{I}_d^2}{\mathcal{J}_d},
\]
where \begin{align*}
\mathcal{I}_d &= \int_{\mathbb{S}^{d-2}}\int_{[-1,1]^2} \varphi_{\scriptscriptstyle \Sigma_{\tau u}}(-\kappa z_1 e+\kappa z_2 u) \,\dd z \dd \sigma(u),\quad\text{and}\\
\mathcal{J}_d&= \int_{\substack{u,v\in \mathbb{S}^{d-2}\\u\neq \pm v}} \int_{[-1,1]^4} \varphi_{\scriptscriptstyle \Sigma_{\tau u,\tau v}} \left(\begin{pmatrix}-\kappa z_1e+\kappa z_2u\\ -\kappa w_1e+\kappa w_2v\end{pmatrix}\right)\, \dd z\, \dd w \,\dd \sigma(u)\, \dd \sigma(v).
\end{align*}
Here $\kappa=8\tau/[\pi(1+\tau^2)^2]$, and $\varphi_\Sigma$ is the density of $N(0,\Sigma)$.
By Lemma~\ref{lem:tube}, we have
\begin{align}
\Sigma_{\tau u}&=\tau^2(I_d+2uu^\top),\quad\text{and}\label{eq:probability one point covariance}\\
\Sigma_{\tau u,\tau v}
&=\begin{pmatrix}\tau^2(I_d+2uu^\top)&C(u,v)\\ C(u,v)^\top&\tau^2(I_d+2vv^\top)\end{pmatrix},\label{eq:probability two point covariance}
\end{align}
where \[C(u,v)=\tau^2 \E\left[(a^\top u)(a^\top v)\sign \left(\tau^2(a^\top u)^2-(a^\top e)^2\right)\sign \left(\tau^2(a^\top v)^2-(a^\top e)^2\right)aa^\top\right].\]
We claim that
\begin{equation}\label{eq:probability integral bounds}
\mathcal I_d=\kappa^{-2}p_\tau\Lambda_d,
\qquad\text{and}\qquad
\mathcal J_d\le\kappa^{-4}C\Lambda_d^2,
\end{equation}
where
\begin{equation}\label{eq:probability p tau definition}
\Lambda_d=\sigma(\mathbb S^{d-2})(2\pi\tau^2)^{-(d-2)/2},
\quad\text{and}\quad
p_\tau=
\left(2\Phi\left(\frac{\kappa}{\tau}\right)-1\right)
\left(2\Phi\left(\frac{\kappa}{\sqrt3\,\tau}\right)-1\right)>0,
\end{equation}
where $\Phi$ denotes the standard Gaussian distribution function, and $C<\infty$ is a numerical constant. Since $\tau$ and $\kappa$ are numerical constants independent of $d$, these bounds will imply the theorem.

The roadmap for the rest of the proof goes as follows. We first characterize $\mathcal I_d$, then bound $\mathcal J_d$ in terms of the joint tangential density. Finally, we show that the constant in the latter bound can be chosen independently of $d$.

\paragraph{Step 1: Characterizing $\mathcal I_d$.}
For each $u,v\in\mathbb S^{d-2}$, let $(G_u,G_v)$ be an $\R^{2d}$ Gaussian vector with mean zero and covariance $\Sigma_{\tau u,\tau v}$.

For each $u$, choose an orthonormal basis $b_1,\ldots,b_{d-2}$ of $\{e,u\}^\perp$, and set
\[
H_u=(b_1^\top G_u,\ldots,b_{d-2}^\top G_u)^\top,
\quad\text{and}\quad N_u=(e^\top G_u,u^\top G_u).
\]
Equation~\eqref{eq:probability one point covariance} implies that $H_u$ and $N_u$ are independent, with
\[
H_u\sim N(0,\tau^2I_{d-2}),\qquad\text{and}\qquad
N_u\sim N\bigl(0,\operatorname{diag}(\tau^2,3\tau^2)\bigr).
\]
We note that the law of $H_u$ is independent of the choice of orthonormal basis due to rotational invariance.
Fix $u\in\mathbb S^{d-2}$. The vector $-\kappa z_1e+\kappa z_2u$ has zero component in $\{e,u\}^\perp$, and its normal coordinates are $(-\kappa z_1,\kappa z_2)$. Set $(s,t)=(-\kappa z_1,\kappa z_2)$, so that $\dd z=\kappa^{-2}\dd s\,\dd t$ and $(s,t)\in[-\kappa,\kappa]^2$. Then, we have
\[
\begin{aligned}
\int_{[-1,1]^2}
\varphi_{\scriptscriptstyle\Sigma_{\tau u}}(-\kappa z_1e+\kappa z_2u)\,\dd z
&=p_{H_u}(0)\int_{[-1,1]^2}p_{N_u}(-\kappa z_1,\kappa z_2)\,\dd z\\
&=\kappa^{-2}p_{H_u}(0)\int_{[-\kappa,\kappa]^2}p_{N_u}(s,t)\,\dd s\,\dd t\\
&=\kappa^{-2}(2\pi\tau^2)^{-(d-2)/2}
\left(2\Phi\left(\frac{\kappa}{\tau}\right)-1\right)
\left(2\Phi\left(\frac{\kappa}{\sqrt3\,\tau}\right)-1\right),
\end{aligned}
\]
where the last equality uses $p_{H_u}(0)=(2\pi\tau^2)^{-(d-2)/2}$ and the fact that the two coordinates of $N_u$ are independent with laws $N(0,\tau^2)$ and $N(0,3\tau^2)$, respectively.
The integral above is constant in $u\in\mathbb S^{d-2}$. Integrating over $u$, we have
\[
\mathcal I_d
=\frac{\sigma(\mathbb S^{d-2})}{\kappa^2(2\pi\tau^2)^{(d-2)/2}}
\left(2\Phi\left(\frac{\kappa}{\tau}\right)-1\right)
\left(2\Phi\left(\frac{\kappa}{\sqrt3\,\tau}\right)-1\right).
\]
This establishes the first equality in~\eqref{eq:probability integral bounds}.

\paragraph{Step 2: Bounding $\mathcal J_d$.}
For $u\ne\pm v$, Lemma~\ref{lem:tube}(\ref{lem:tube:III}) ensures that the joint density exists. Set $\xi(z,w)=(-\kappa z_1,\kappa z_2,-\kappa w_1,\kappa w_2)$. Since the coordinate transformation from $(G_u,G_v)$ to $(H_u,H_v,N_u,N_v)$ is orthogonal,
\[
\begin{aligned}
p_{G_u,G_v}(-\kappa z_1e+\kappa z_2u,\,-\kappa w_1e+\kappa w_2v)
&=p_{H_u,H_v,N_u,N_v}(0,0,\xi(z,w))\\
&=p_{H_u,H_v}(0,0)\,
p_{(N_u,N_v)\mid H_u=H_v=0}
(\xi(z,w)),
\end{aligned}
\]
where the second factor is the Gaussian density of $(N_u,N_v)$ conditional on $(H_u,H_v)=(0,0)$. The change of variables $(z,w)\mapsto\xi(z,w)$ has absolute Jacobian $\kappa^4$. Hence, for fixed $u,v\in\mathbb S^{d-2}$ with $u\ne\pm v$,
\[
\begin{aligned}
\int_{[-1,1]^4}
p_{(N_u,N_v)\mid H_u=H_v=0}
(\xi(z,w))
\,\dd z\,\dd w
&=\kappa^{-4}\int_{[-\kappa,\kappa]^4}
p_{(N_u,N_v)\mid H_u=H_v=0}(\xi)\,\dd\xi\\
&=\kappa^{-4}
\PP\bigl((N_u,N_v)\in[-\kappa,\kappa]^4\mid H_u=H_v=0\bigr).
\end{aligned}
\]
Therefore,
\begin{align*}
\mathcal J_d
&=\kappa^{-4}\int_{\substack{u,v\in\mathbb S^{d-2}\\u\ne\pm v}}
p_{H_u,H_v}(0,0)
\PP\bigl((N_u,N_v)\in[-\kappa,\kappa]^4\mid H_u=H_v=0\bigr)
\,\dd\sigma(u)\,\dd\sigma(v)\\
&\le\kappa^{-4}\int_{\substack{u,v\in\mathbb S^{d-2}\\u\ne\pm v}}
p_{H_u,H_v}(0,0)\,\dd\sigma(u)\,\dd\sigma(v)
=:\kappa^{-4}Q_d.
\end{align*}
It remains to prove $Q_d\le C\sigma(\mathbb S^{d-2})^2(2\pi\tau^2)^{-(d-2)}$ uniformly for $d\ge4$.

Changing the orthonormal bases of $\{e,u\}^\perp$ and $\{e,v\}^\perp$ leaves $p_{H_u,H_v}(0,0)$ unchanged. We therefore choose the following bases to simplify its computation. Fix $u,v\in\mathbb S^{d-2}$ with $u\ne\pm v$, put $c=u^\top v$ and $s=\sqrt{1-c^2}$, and choose an orthonormal basis $t_1,\ldots,t_{d-3}$ of $\{e,u,v\}^\perp$. Define
\[
h_u=\frac{v-cu}{s},\qquad\text{and}\qquad h_v=\frac{cv-u}{s}.
\]
The tuples $(h_u,t_1,\ldots,t_{d-3})$ and $(h_v,t_1,\ldots,t_{d-3})$ are orthonormal bases of $\{e,u\}^\perp$ and $\{e,v\}^\perp$, respectively. Under this choice of bases,
\[
H_u=\begin{pmatrix}h_u^\top G_u\\t_1^\top G_u\\\vdots\\t_{d-3}^\top G_u\end{pmatrix},
\qquad\text{and}\qquad
H_v=\begin{pmatrix}h_v^\top G_v\\t_1^\top G_v\\\vdots\\t_{d-3}^\top G_v\end{pmatrix}.
\]
Equation~\eqref{eq:probability one point covariance} implies $\operatorname{Cov}(H_u)=\operatorname{Cov}(H_v)=\tau^2I_{d-2}$.
For a standard Gaussian vector $a\in\R^d$, define the following real standard Gaussian variables:
\[
P=a^\top e,\enspace X=a^\top u,\enspace Y=a^\top v,\enspace
Z=a^\top h_u,\enspace T=a^\top h_v,\quad\text{and}\quad L_j=a^\top t_j
\enspace(1\le j\le d-3).
\]
The variables $P,X,Z,L_1,\ldots,L_{d-3}$ are independent, while $Y=cX+sZ$ and $T=-sX+cZ$. Also define
\[
S_X=\sign(\tau^2X^2-P^2),\qquad\text{and}\qquad
S_Y=\sign(\tau^2Y^2-P^2).
\]
By the definition of $C(u,v)$ in~\eqref{eq:probability two point covariance},
\[
\begin{aligned}
\operatorname{Cov}(h_u^\top G_u,h_v^\top G_v)
&=h_u^\top C(u,v)h_v\\
&=\tau^2\E[(a^\top u)(a^\top v)S_XS_Y(a^\top h_u)(a^\top h_v)]\\
&=\tau^2\E[XY S_XS_Y ZT].
\end{aligned}
\]
Likewise, for $1\le i,j\le d-3$, independence of the $L_j$ from $(P,X,Z)$ implies
\[
\begin{aligned}
\operatorname{Cov}(t_i^\top G_u,t_j^\top G_v)
&=t_i^\top C(u,v)t_j
=\tau^2\E[XY S_XS_Y L_iL_j]=\tau^2\E[XY S_XS_Y] \delta_{ij},\\
\operatorname{Cov}(h_u^\top G_u,t_j^\top G_v)
&=h_u^\top C(u,v)t_j
=\tau^2\E[XY S_XS_Y ZL_j]=0,\quad\text{and}\\
\operatorname{Cov}(t_i^\top G_u,h_v^\top G_v)
&=t_i^\top C(u,v)h_v
=\tau^2\E[XY S_XS_Y L_iT]=0.
\end{aligned}
\]
Since $P,X,Z$ are independent standard Gaussians and $Y=cX+sZ$, $T=-sX+cZ$, these covariances depend on $u,v$ only through their inner product $c$. We write
\[
r(c)=\E[XZ YT S_XS_Y],
\qquad\text{and}\qquad
q(c)=\E[XY S_XS_Y].
\]
Thus
\[
\operatorname{Cov}\begin{pmatrix}H_u\\H_v\end{pmatrix}
=\tau^2\begin{pmatrix}I_{d-2}&D\\D&I_{d-2}\end{pmatrix},
\qquad\text{and}\qquad D=\operatorname{diag}(r(c),q(c),\ldots,q(c)).
\]
By the Schur complement,
\[
\begin{aligned}
\det\operatorname{Cov}\begin{pmatrix}H_u\\H_v\end{pmatrix}
&=\tau^{4(d-2)}\det\begin{pmatrix}I_{d-2}&D\\D&I_{d-2}\end{pmatrix}\\
&=\tau^{4(d-2)}\det(I_{d-2}-D^2)
=\tau^{4(d-2)}(1-r(c)^2)(1-q(c)^2)^{d-3}.
\end{aligned}
\]
Since $(H_u,H_v)$ is a centered Gaussian vector in $\R^{2(d-2)}$, its density at zero is
\begin{equation}\label{eq:probability tangent density}
p_{H_u,H_v}(0,0)
=(2\pi\tau^2)^{-(d-2)}(1-r(c)^2)^{-1/2}(1-q(c)^2)^{-(d-3)/2}.
\end{equation}

\paragraph{Step 3: Uniform bound in dimension.}
We verify the scalar hypotheses of Lemma~\ref{lem:geodesic sign correlation} numerically using interval arithmetic in Section~\ref{app:numerics}. Thus Lemmas~\ref{lem:common sign correlation} and~\ref{lem:geodesic sign correlation} give, with $b=\gamma>0$,
\[
1-q(c)^2\ge(1-c^2)(1+\gamma c^2),
\qquad\text{and}\qquad
1-r(c)^2\ge b\sqrt{1-c^2}.
\]
The density formula~\eqref{eq:probability tangent density} implies
\begin{align*}
Q_d
&=(2\pi\tau^2)^{-(d-2)}
\int_{\substack{u,v\in\mathbb S^{d-2}\\u\ne\pm v}}
(1-r(u^\top v)^2)^{-1/2}
(1-q(u^\top v)^2)^{-(d-3)/2}
\,\dd\sigma(u)\,\dd\sigma(v)\\
&=\sigma(\mathbb S^{d-2})^2(2\pi\tau^2)^{-(d-2)}
\int_{\substack{u,v\in\mathbb S^{d-2}\\u\ne\pm v}}
(1-r(u^\top v)^2)^{-1/2}
(1-q(u^\top v)^2)^{-(d-3)/2}
\,\dd u\,\dd v\\
&\le b^{-1/2}\sigma(\mathbb S^{d-2})^2(2\pi\tau^2)^{-(d-2)} A_d
\int_{-1}^1(1-c^2)^{-3/4}
(1+\gamma c^2)^{-(d-3)/2}\,\dd c.
\end{align*}
Here $\dd u$ and $\dd v$ denote the uniform probability measures on $\mathbb S^{d-2}$. Under their product measure, $c=u^\top v$ has density
\[
A_d(1-c^2)^{(d-4)/2}\qquad\text{where }
A_d=\frac{\Gamma((d-1)/2)}{\sqrt\pi\,\Gamma((d-2)/2)}
\le \frac{\sqrt{d-2}}{\sqrt{2\pi}}.
\]
The last inequality follows directly from Gautschi's inequality.
In the following estimates, $C>0$ denotes a constant independent of $d$ that may change from line to line. Since $\log(1+\gamma c^2)\ge \gamma c^2/(1+\gamma/4)$ on $|c|\le1/2$, we have, for $d\ge4$,
\begin{align*}
\int_{-1}^1(1-c^2)^{-3/4}(1+\gamma c^2)^{-(d-3)/2}\,\dd c
&\le C\int_{-1/2}^{1/2}e^{-\gamma(d-3)c^2/[2(1+\gamma/4)]}\,\dd c\\
&\quad \qquad +(1+\gamma/4)^{-(d-3)/2}
\int_{1/2<|c|<1}(1-c^2)^{-3/4}\,\dd c\\
&\le C(d-3)^{-1/2}+C(1+\gamma/4)^{-(d-3)/2}
\le C(d-2)^{-1/2}.
\end{align*}
Since $A_d\le \sqrt{(d-2)/(2\pi)}$, we obtain $Q_d\le C\sigma(\mathbb S^{d-2})^2(2\pi\tau^2)^{-(d-2)}$ uniformly in $d$. This proves the second bound in~\eqref{eq:probability integral bounds}.

Finally, \eqref{eq:probability integral bounds} implies
\[
\frac{\mathcal I_d^2}{\mathcal J_d}
\ge\frac{p_\tau^2}{C}=:c_\star>0
\qquad\text{for }d\ge4.
\]
Combining this bound with the earlier Paley--Zygmund estimate proves the theorem.
\end{proof}

\subsection{Numerical certification of the probability bound}\label{app:numerics}

In this appendix, we make the probability bound in Theorem~\ref{thm:probability appendix} explicit. The preceding proof gives, for every fixed $d\ge4$,
\begin{equation}\label{eq:numerics probability reduction}
\liminf_{n\to\infty}\PP(N_{d,n}>0)
\ge\frac{\mathcal I_d^2}{\mathcal J_d}
\ge\frac{p_\tau^2}{C},
\end{equation}
where $C$ is any constant satisfying~\eqref{eq:probability integral bounds} uniformly in $d$. Combining analytic estimates with certified interval computations, we show that $C=3.054$ is admissible and that $p_\tau^2/3.054>0.1361$. This establishes the choice $c_\star=0.136$. We defer the details for the reproduction of these results to the end of this appendix.

\paragraph{Scalar parameters and correlation bounds.}
Recall that $\tau$ is the positive root of
\[
\arctan t+\frac{t}{1+t^2}-\frac\pi4=0.
\]
The derivative of the left-hand side is $2/(1+t^2)^2>0$. We enclose $\tau$ between rational endpoints by verifying, with interval arithmetic, that the left-hand side has opposite signs at the endpoints.
For the coefficient $w_0$ in the proof of Lemma~\ref{lem:common sign correlation}, we use the identity
\begin{equation}\label{eq:numerics w0 formula}
\begin{aligned}
w_0={}&1-\frac8\pi\arctan\frac1\tau
+\frac8\pi\arcsin\frac1{1+\tau^2}\\
&\quad+\frac{8\tau}{\pi}\left(
\frac1{1+\tau^2}
-\frac2{(1+\tau^2)\sqrt{\tau^2+2}}
+\frac1{(\tau^2+2)^{3/2}}\right).
\end{aligned}
\end{equation}
Evaluating this expression and $\gamma=1-w_0^2$ on the certified interval for $\tau$ gives
\begin{equation}\label{eq:numerics scalar enclosures}
\begin{aligned}
0.44161079&<\tau<0.44161080,\\
0.60010962&<w_0<0.60010963,\quad\text{and}\\
0.63986844&<\gamma<0.63986845.
\end{aligned}
\end{equation}
The displayed enclosures are rounded outward from the narrower intervals used in the computations. In particular,
\begin{equation}\label{eq:numerics formal hypotheses}
\frac25\le\tau\le\frac9{20},\qquad
w_0\ge\frac35,\qquad\text{and}\qquad
\frac4{47}<\gamma<\frac{16}{25}.
\end{equation}
Lemma~\ref{lem:common sign correlation} gives
\begin{equation}\label{eq:numerics common sign bound}
1-q(c)^2\ge(1-c^2)(1+\gamma c^2)
\qquad\text{for }-1\le c\le1.
\end{equation}
The first two bounds in~\eqref{eq:numerics formal hypotheses} allow us to apply Lemma~\ref{lem:geodesic sign correlation}, giving
\begin{equation}\label{eq:numerics sharp geodesic bound}
1-r(c)^2\ge\gamma\sqrt{1-c^2}
\qquad\text{for }-1<c<1,
\end{equation}
with equality if, and only if, $c=0$. Thus we may take $b=\gamma$ in the preceding proof.

\paragraph{Uniform bound for the angular integral.}
Combining the correlation bounds~\eqref{eq:numerics common sign bound} and~\eqref{eq:numerics sharp geodesic bound} with the density formula~\eqref{eq:probability tangent density} shows that~\eqref{eq:probability integral bounds} holds whenever
\begin{equation}\label{eq:numerics sufficient C}
C\ge\gamma^{-1/2}\sup_{d\ge4}H_d(\gamma),
\end{equation}
where
\[
H_d(\gamma)=A_d\int_{-1}^1(1-c^2)^{-3/4}(1+\gamma c^2)^{-(d-3)/2}\,\dd c,
\quad\text{and}\quad A_d=\frac{\Gamma((d-1)/2)}{\sqrt\pi\,\Gamma((d-2)/2)}.
\]
The weight $(1-c^2)^{-3/4}$ is obtained by bounding the inverse covariance factors and multiplying by the spherical correlation density $A_d(1-c^2)^{(d-4)/2}$.
Next, we consider two cases.

\emph{Case 1: $4\le d\le49$.} The substitution $u=c^2$ and the binomial expansion give
\begin{equation}\label{eq:numerics angular series}
\begin{aligned}
H_d(\gamma)&=A_d B(1/2,1/4)\sum_{n=0}^{\infty}(-1)^n a_n,\\
a_0&=1,\qquad\text{and}\qquad
a_{n+1}=a_n\gamma\frac{((d-3)/2+n)(1/2+n)}{(3/4+n)(n+1)},
\end{aligned}
\end{equation}
where $B(1/2,1/4)=\Gamma(1/2)\Gamma(1/4)/\Gamma(3/4)$. For $0\le\gamma\le16/25$, the nonnegative terms satisfy $a_{n+1}\le(4/5)a_n$ for all $n\ge100$. Consequently, for every even $N\ge100$, the partial sum $S_N=\sum_{n=0}^N(-1)^n a_n$ satisfies
\begin{equation}\label{eq:numerics angular truncation}
A_d B(1/2,1/4)(S_N-a_{N+1})
\le H_d(\gamma)\le A_d B(1/2,1/4)S_N.
\end{equation}
We evaluate the partial sums and the first omitted term using interval arithmetic. Together with~\eqref{eq:numerics angular truncation}, this accounts for both rounding and truncation errors and yields
\begin{equation}\label{eq:numerics finite maximum}
\begin{gathered}
H_d(\gamma)<H_6(\gamma)\quad\text{for }4\le d\le49\text{ with }d\ne6,\quad\text{and}\\
2.44237353<H_6(\gamma)<2.44237354.
\end{gathered}
\end{equation}
For each comparison, the upper endpoint of the interval for $H_d(\gamma)$ lies below the lower endpoint of the interval for $H_6(\gamma)$.

\emph{Case 2: $d\ge50$.} We use the analytic estimate
\begin{equation}\label{eq:numerics angular tail}
\begin{aligned}
H_d(\gamma)\le{}&
\left(\frac43\right)^{3/4}
\sqrt{\frac{1+\gamma/4}{\gamma}}\sqrt{\frac{48}{47}}\\
&+B(1/2,1/4)\sqrt{\frac{24}{\pi}}
\exp\!\left(-\frac{47}{2}\log(1+\gamma/4)\right).
\end{aligned}
\end{equation}
This bound holds when $\gamma>0$ and $48\log(1+\gamma/4)\ge1$, as ensured by~\eqref{eq:numerics formal hypotheses}. Its proof splits the integral at $|c|=1/2$, applies a Gaussian bound on the middle interval and a beta-integral bound on the remainder, and uses monotonicity in $d$. Interval evaluation of the right-hand side gives a value less than $2.131647$, below the lower bound for $H_6(\gamma)$ in~\eqref{eq:numerics finite maximum}.

Combining the two cases, we obtain
\begin{equation}\label{eq:numerics angular supremum}
\sup_{d\ge4}H_d(\gamma)=H_6(\gamma).
\end{equation}
Since $\gamma>0.639868$, the displayed enclosures imply
\[
\gamma^{-1/2}\sup_{d\ge4}H_d(\gamma)
=\frac{H_6(\gamma)}{\sqrt\gamma}
<\frac{2.44237354}{\sqrt{0.639868}}
<3.054.
\]
Thus $C=3.054$ satisfies~\eqref{eq:numerics sufficient C}.

\paragraph{Evaluation of the probability bound.}
By~\eqref{eq:probability p tau definition},
\[
p_\tau=
\operatorname{erf}\!\left(\frac{\kappa}{\sqrt2\,\tau}\right)
\operatorname{erf}\!\left(\frac{\kappa}{\sqrt6\,\tau}\right),
\qquad\text{and}\qquad \kappa=\frac{8\tau}{\pi(1+\tau^2)^2},
\]
where $\operatorname{erf}(x)=2\Phi(\sqrt2x)-1$. Both arguments are nonnegative and have squares less than $3$, as verified by interval arithmetic. In this range, the terms of the alternating series for $\operatorname{erf}$ decrease in absolute value, giving
\[
\left|\operatorname{erf}(x)
-\frac2{\sqrt\pi}\sum_{n=0}^{N}
\frac{(-1)^n x^{2n+1}}{n!(2n+1)}\right|
\le\frac2{\sqrt\pi}\frac{x^{2N+3}}{(N+1)!(2N+3)}
\qquad\text{for }N\ge0.
\]
Interval evaluation with this remainder gives
\[
0.64480748<p_\tau<0.64480749.
\]
Combining this enclosure with~\eqref{eq:numerics probability reduction}, we conclude that, for every fixed $d\ge4$,
\[
\liminf_{n\to\infty}\PP(N_{d,n}>0)
\ge\frac{p_\tau^2}{3.054}
>\frac{(0.64480748)^2}{3.054}
>0.1361.
\]
Since $N_{d,n}\le M_{d,n}$, this establishes the choice $c_\star=0.136$ in Theorems~\ref{thm:probability appendix} and~\ref{thm:probability}.

\paragraph{Verification and reproduction.}
The accompanying Lean formalization verifies the identity~\eqref{eq:numerics w0 formula}, the correlation bounds, the angular series and its remainder estimate, and the uniform bound~\eqref{eq:numerics angular tail}. The moment theorems and covariance reduction are proved in the paper but are not formalized in Lean. The numerical enclosures are computed separately, outside Lean, using interval arithmetic with outward rounding and the explicit remainder bounds above; the error-function evaluation is also outside the formalization. All displayed terminating decimals denote exact rational numbers. Supplementary quadrature checks are not used to certify the inequalities.

For the source code, full-precision certificates, and reproduction instructions, see \reproducibilitylink. The formal proofs are in \texttt{04-Lean}, and the numerical computations are in \texttt{01-numerical-lower-bound}. The repository documentation records the formal declarations and computational details. In the code and formal proofs, the dimension parameter is $d-2$ rather than $d$.

%%% Local Variables:
%%% mode: LaTeX
%%% TeX-master: "../main"
%%% End:

%%% Local Variables:
%%% mode: LaTeX
%%% TeX-master: "../main"
%%% End:

\section{Probabilistic tools for counting minima}\label{app:probability-supporting-tools}
This appendix collects the probabilistic tools used to compute the first and second moments of $N_{d,n}$ in Appendix~\ref{app:probability}, using the notation from Appendix~\ref{sec:common notation}. We prove Lemma~\ref{lem:regularity constants} and cover joint density bounds, Kac--Rice formulas, uniform local limit approximations, and mean and covariance estimates near the population saddle sphere.

\subsection{Proof of Lemma~\ref{lem:regularity constants}}\label{app:proof regularity constants}

\begin{proof}[Proof of Lemma~\ref{lem:regularity constants}]
By rotational invariance of the standard Gaussian measurements, we may assume without loss of generality that $\xs=e$, where $e$ is the first standard basis vector.
Fix $I=\{i_1<\cdots<i_d\}$. Since all the vectors $a_i$ are i.i.d.\ and every ordered active set uses the same sign pattern, the conditional expectation $a(x)$ is independent of the choice of $I$.

We first show that $a(\tau u)>0$. Fix $u\in\mathbb S^{d-2}$. Conditional on $X_I(\tau u)=0$, the vectors $a_{i_j}$ are independent standard Gaussian vectors on the hyperplanes $(\tau u-s_j e)^\perp$, for $j=1,\ldots,d$. We will construct a deterministic matrix whose rows lie in these hyperplanes and which satisfies the inequalities defining $\mathcal K_{\tau u}$ strictly.
Choose an orthonormal basis $v_1,\ldots,v_{d-2}$ of $\{e,u\}^\perp$, and let $B_0$ be the matrix with rows $b_j^\top$, where
\[
\begin{aligned}
b_1&=\tau e+u, & b_j&=\tau e+u+v_{j-1}\quad\text{for }2\le j\le d-2,\\
b_{d-1}&=-\tau e+u, & \text{and}\quad b_d&=-\tau e+u+v_{d-2}.
\end{aligned}
\]
Each vector $b_j$ satisfies $b_j^\top(\tau u-s_j e)=0$ and $b_j^\top(\tau u)=\tau$. To compute the determinant, subtract row $1$ from rows $2,\ldots,d-2$ and row $d-1$ from row $d$. Reorder the resulting rows as $b_1^\top,b_{d-1}^\top,v_1^\top,\ldots,v_{d-2}^\top$, and express them in coordinates relative to the orthonormal basis $(e,u,v_1,\ldots,v_{d-2})$. These operations preserve the absolute determinant, so
\[
\abs{\det B_0}
=\left|\det\begin{pmatrix}\tau&1&0\\-\tau&1&0\\0&0&I_{d-2}\end{pmatrix}\right|
=\left|\det\begin{pmatrix}\tau&1\\-\tau&1\end{pmatrix}\right|=2\tau.
\]
Also, $b_j^\top(\tau u)=\tau$ implies $Q(\tau u;B_0)=\tau B_0^\top$, and hence $\abs{\det Q(\tau u;B_0)}=2\tau^{d+1}$.

Since $\max_j\norm{b_j}=\sqrt{\tau^2+2}$ and $0<\tau<1$, choose $M_d>\sqrt{\tau^2+2}$ and $0<\lambda_d<2\tau^{d+1}$. All three inequalities defining $\mathcal K_{\tau u}$ are then strict at $B_0$. By continuity, they also hold for all matrices sufficiently close to $B_0$ whose rows remain in the conditioning hyperplanes. The joint conditional Gaussian density is positive on the product of these hyperplanes, so this neighborhood has positive conditional probability. Thus
\[
\PP\left(A_I\in\mathcal K_{\tau u}\;\middle|\;X_I(\tau u)=0\right)>0.
\]
On $\mathcal K_{\tau u}$, both absolute determinants in the definition of $a$ are at least $\lambda_d$, so
\[
a(\tau u)\ge\lambda_d^2\,\PP\left(A_I\in\mathcal K_{\tau u}\;\middle|\;X_I(\tau u)=0\right)>0.
\]

It remains to show that $a(\tau u)$ does not depend on $u$. For any $v\in\mathbb S^{d-2}$, choose an orthogonal map $R$ with $Re=e$ and $Ru=v$. This map sends each hyperplane $(\tau u-s_j e)^\perp$ to $(\tau v-s_j e)^\perp$ and preserves the standard Gaussian law on it. Therefore, the law of $A_I R^\top$ conditional on $X_I(\tau u)=0$ is the law of $A_I$ conditional on $X_I(\tau v)=0$. Rotation preserves the row norms and the absolute determinant of $A_I$, and
\[
Q(\tau v;A_I R^\top)=R Q(\tau u;A_I).
\]
Consequently, $A_I\in\mathcal K_{\tau u}$ if, and only if, $A_I R^\top\in\mathcal K_{\tau v}$. Since $R$ also preserves the two absolute determinant factors in the definition of $a$, applying this rotation to the conditional expectation shows that $a(\tau u)=a(\tau v)$ for all $u,v\in\mathbb S^{d-2}$. Denote this common positive value by $\vartheta_d$.

We next prove continuity near $\mathcal S$. Fix an integer $n_0$ with $n_0^{-1/2}<\tau$. The set $T_{n_0}$ is compact, its interior contains $\mathcal S$, and $x\ne\pm e$ for every $x\in T_{n_0}$. To compare the conditional laws at nearby points, let $g_1,\ldots,g_d$ be independent $N(0,I_d)$ vectors and use these same vectors to define, for every $x\in T_{n_0}$,
\[
\widehat a_j(x)=g_j-\frac{g_j^\top(x-s_j e)}{\norm{x-s_j e}^2}(x-s_j e)
\qquad\text{for }j=1,\ldots,d,
\]
and let $\widehat A(x)$ be the matrix with rows $\widehat a_j(x)^\top$. The vector $\widehat a_j(x)$ is the orthogonal projection of $g_j$ onto $(x-s_j e)^\perp$. Thus $\widehat A(x)$ has the law of $A_I$ conditional on $X_I(x)=0$, and
\[
a(x)=\E_{g_1,\ldots,g_d}\!\left[\abs{\det\widehat A(x)}\abs{\det Q(x;\widehat A(x))}1_{\mathcal K_x}\bigl(\widehat A(x)\bigr)\right].
\]
Fix $x\in T_{n_0}$ and a sequence $x_j\in T_{n_0}$ with $x_j\to x$. The projection formula implies $\widehat A(x_j)\to\widehat A(x)$ almost surely, so both determinant factors converge almost surely. To handle the indicator, we must also exclude equality in the inequalities defining $\mathcal K_x$. For this fixed $x$, each such equality requires one of the polynomials
\[
\norm{\widehat a_j(x)}^2-M_d^2,\qquad
\bigl(\det\widehat A(x)\bigr)^2-\lambda_d^2,\qquad\text{and}\qquad
\bigl(\det Q(x;\widehat A(x))\bigr)^2-\lambda_d^2
\]
in the components of $(g_1,\ldots,g_d)$ to vanish. None is identically zero: when all $g_j=0$, their values are $-M_d^2$ or $-\lambda_d^2$. Lemma~\ref{lem:zeros of polynomial} therefore shows that these equalities occur with probability zero. Since all the quantities defining $\mathcal K_x$ depend continuously on $x$ and the matrix entries, it follows that
\[
1_{\mathcal K_{x_j}}\bigl(\widehat A(x_j)\bigr)\longrightarrow 1_{\mathcal K_x}\bigl(\widehat A(x)\bigr)\qquad\text{almost surely}.
\]
Hence the full integrand defining $a(x_j)$ converges almost surely to that defining $a(x)$. To pass to expectations, note that whenever $\widehat A(x)\in\mathcal K_x$, Hadamard's inequality and the row bounds imply
\[
\abs{\det\widehat A(x)}\le M_d^d,\qquad\text{and}\qquad
\abs{\det Q(x;\widehat A(x))}
\le\prod_{j=1}^d\left(\abs{\widehat a_j(x)^\top x}\norm{\widehat a_j(x)}\right)
\le\norm{x}^dM_d^{2d}.
\]
Since $\norm{x}\le\tau+2n_0^{-1/2}$ on $T_{n_0}$, the integrand is bounded by the deterministic constant $(\tau+2n_0^{-1/2})^dM_d^{3d}$ for every $x\in T_{n_0}$. Dominated convergence now implies $a(x_j)\to a(x)$, proving continuity on $T_{n_0}$ and hence near $\mathcal S$.

Finally, continuity on the compact set $T_{n_0}$ implies uniform continuity there. For $n\ge n_0$ and $x=\alpha e+(\tau+\beta)u\in T_n$, both $x$ and $\tau u$ lie in $T_{n_0}$, and
\[
\norm{x-\tau u}=\sqrt{\alpha^2+\beta^2}\le\sqrt2n^{-1/2}.
\]
This bound holds uniformly over $x\in T_n$. Since $a(\tau u)=\vartheta_d$, uniform continuity implies
\[
\sup_{x\in T_n}\abs{a(x)-\vartheta_d}\to0\qquad\text{as }n\to\infty.
\]
\end{proof}

\subsection{Joint density bounds}\label{app:conditional-gaussian-tools}

Recall the active equations $X_I$ and their solutions $x_I$ from~\eqref{eq:def X_I} and~\eqref{eq:def x_I}. For an active set $I=\{i_1<\cdots<i_d\}$, write
\[
s_{i_k}^I=s_k\qquad\text{for }1\le k\le d
\]
for the row signs induced by the fixed pattern~\eqref{eq:sign pattern}. For distinct active sets $I,J$ and $x,y\in\R^d$, define
\[
\Gamma_{I,J}(x,y)=\operatorname{Cov}\begin{pmatrix}X_I(x)\\X_J(y)\end{pmatrix}.
\]
For $v,w\in\R^d$, set
\[
\Delta(v,w)=\sqrt{\norm{v}^2\norm{w}^2-(v^\top w)^2}.
\]

\begin{lemma}\label{lem:probability zero degenerate covariance}
Fix integers $d\ge4$ and $n\ge d$, and a deterministic vector $\xs\in\R^d$ with $\norm{\xs}=1$. Let $a_1,\ldots,a_n$ be independent $N(0,I_d)$ vectors. Then, for every pair of distinct index sets $I,J\subset[n]$ with $\abs I=\abs J=d$ and every $x,y\in\R^d$,
\begin{equation}\label{eq:joint active covariance determinant}
\det\Gamma_{I,J}(x,y)
=\prod_{i\in I\setminus J}\norm{x-s_i^I\xs}^2
\prod_{j\in J\setminus I}\norm{y-s_j^J\xs}^2
\prod_{\ell\in I\cap J}\Delta^2(x-s_\ell^I\xs,y-s_\ell^J\xs).
\end{equation}
Moreover, almost surely, this covariance map is positive definite at the random points $(x_I,x_J)$ for every such pair $I,J$:
\[
\Gamma_{I,J}(x_I,x_J)\succ 0.
\]
\end{lemma}
\begin{proof}
By rotational invariance of the standard Gaussian measurements, we may assume without loss of generality that $\xs=e$, where $e$ is the first standard basis vector. Fix $I\ne J$. For $i\in I\setminus J$, the coordinate $a_i^\top(x-s_i^I e)$ is independent of all other coordinates and has variance $\norm{x-s_i^I e}^2$; the same holds for $j\in J\setminus I$. For each $\ell\in I\cap J$, the two coordinates involving $a_\ell$ have covariance matrix
\[
\begin{pmatrix}
\norm{x-s_\ell^I e}^2 &(x-s_\ell^I e)^\top(y-s_\ell^J e)\\
(x-s_\ell^I e)^\top(y-s_\ell^J e)&\norm{y-s_\ell^J e}^2
\end{pmatrix},
\]
whose determinant is $\Delta^2(x-s_\ell^I e,y-s_\ell^J e)$. Therefore,
\[
\det\Gamma_{I,J}(x,y)
=\prod_{i\in I\setminus J}\norm{x-s_i^I e}^2
\prod_{j\in J\setminus I}\norm{y-s_j^J e}^2
\prod_{\ell\in I\cap J}\Delta^2(x-s_\ell^I e,y-s_\ell^J e).
\]

Almost surely, every $A_I$ is invertible and $a_i^\top e\ne0$ for all $i$. On this event, no $x_I$ is parallel to $e$. Indeed, since the sign pattern $s$ contains both $1$ and $-1$, if $x_I=\alpha e$, then the active equations $X_I(x)=0$ would require $\alpha$ to equal both $1$ and $-1$, which is impossible.
Therefore, almost surely, we have
\[
\prod_{i\in I\setminus J}\norm{x_I-s_i^I e}^2 \prod_{j\in J\setminus I}\norm{x_J-s_j^J e}^2>0.
\]
Since $\Gamma_{I,J}(x_I,x_J)$ is positive semidefinite, the determinant identity~\eqref{eq:joint active covariance determinant} shows that it remains to prove, almost surely,
\[
\Delta(x_I-s_\ell^I e,x_J-s_\ell^J e)>0
\qquad\text{for every }\ell\in I\cap J.
\]
There is nothing further to check if $I\cap J=\varnothing$. Otherwise, fix $\ell\in I\cap J$. We will show that the two vectors in this expression are linearly dependent with probability zero.

Condition on the rows $(a_i)_{i\in I}$, and set $x=x_I$, $\sigma=s_\ell^I$, and $\sigma'=s_\ell^J$. Then $x$ is fixed and is not parallel to $e$, while the vectors $(a_j)_{j\in J\setminus I}$ remain independent $N(0,I_d)$ vectors. Since $I\ne J$ and $\abs{I}=\abs{J}$, the set $J\setminus I$ is nonempty. Both vectors associated with $\ell$ are nonzero, so if they are linearly dependent, there is a scalar $t\ne0$ such that
\[
x_J=\sigma'e+t(x-\sigma e).
\]
For every $j\in J$, substituting this expression into $a_j^\top(x_J-s_j^J e)=0$ implies
\begin{equation}\label{eq:degenerate covariance row equation}
t\,a_j^\top(x-\sigma e)+(\sigma'-s_j^J)a_j^\top e=0.
\end{equation}
First, suppose some $j\in J\setminus I$ has $s_j^J=\sigma'$. Since $t\ne0$, equation~\eqref{eq:degenerate covariance row equation} requires $a_j^\top(x-\sigma e)=0$. Conditional on the rows in $I$, this has probability zero because $a_j$ is a standard Gaussian vector independent of those rows and $x-\sigma e\ne0$.

It remains to consider the case in which every $j\in J\setminus I$ has sign $-\sigma'$. Suppose first that $J\setminus I$ contains two distinct indices $j,k$. Conditional on the rows in $I$, we have $a_j^\top(x-\sigma e)\ne0$ almost surely, so equation~\eqref{eq:degenerate covariance row equation} for row $j$ determines
\[
t=-\frac{2\sigma'a_j^\top e}{a_j^\top(x-\sigma e)}.
\]
After also conditioning on $a_j$, this value of $t$ is fixed, and the equation for row $k$ requires
\[
a_k^\top\bigl(t(x-\sigma e)+2\sigma'e\bigr)=0.
\]
The vector in parentheses is nonzero because $x$ is not parallel to $e$. Since $a_k$ remains an independent standard Gaussian vector, this equation has conditional probability zero.

The only remaining case is $J\setminus I=\{k\}$ with $s_k^J=-\sigma'$. Each sign occurs at least twice in the fixed pattern because $d\ge4$, so there is a shared row $j\in I\cap J$ with $s_j^J=-\sigma'$. Since $a_j^\top x=s_j^I a_j^\top e$, equation~\eqref{eq:degenerate covariance row equation} for this row becomes
\[
\bigl[t(s_j^I-\sigma)+2\sigma'\bigr]a_j^\top e=0.
\]
We have $a_j^\top e\ne0$, so $s_j^I=\sigma$ would make this equality impossible. Thus $s_j^I=-\sigma$, which forces $t=\sigma'/\sigma$. Substituting this value into the equation for row $k$ implies
\[
a_k^\top\bigl((\sigma'/\sigma)x+\sigma'e\bigr)=0.
\]
The vector in parentheses is fixed under the conditioning and is nonzero because $x$ is not parallel to $e$. Since $a_k$ is a standard Gaussian vector independent of the rows in $I$, this equation also has conditional probability zero.

These cases exhaust all possibilities for $J\setminus I$. Averaging over the rows in $I$ shows that the two vectors associated with $\ell$ are linearly dependent with probability zero. Taking a finite union over $I,J$ and $\ell\in I\cap J$, we conclude that all factors in~\eqref{eq:joint active covariance determinant} are positive almost surely at $(x_I,x_J)$. Thus $\det\Gamma_{I,J}(x_I,x_J)>0$, and positive semidefiniteness implies $\Gamma_{I,J}(x_I,x_J)\succ0$ for every such pair $I,J$, almost surely.
\end{proof}

For distinct active sets $I,J$, let
\[
\mathcal R_{I,J}=\{(x,y)\in\R^{2d}:\Gamma_{I,J}(x,y)\succ0\},
\]
and denote the joint Gaussian density on this region by $p_{X_I(x),X_J(y)}$. We use the angular separation $\delta(x,y)$ from the proof of Theorem~\ref{thm:second moment asymptotic appendix}. For $x\in\R^d$, write $v_\pm=x\mp\xs$, and let $P_{v^\perp}$ denote orthogonal projection onto $v^\perp$. In part~(\ref{lem:density:II}) below, the shared-row counts are
\[
r_\pm=\#\{i\in I\cap J:s_i^I=s_i^J=\pm1\}.
\]

\Needspace{8\baselineskip}
\begin{lemma}[Joint density and angular integral bounds]\label{lem:density}
Fix integers $d\ge4$ and $0\le r\le d-1$, and index sets $I,J\subset[n]$ with $\abs I=\abs J=d$ and $\abs{I\cap J}=r$. The following estimates hold.
\begin{enumerate}[label=\textnormal{(\Roman*)},ref=\Roman*]
\item\label{lem:density:I} There exist constants $n_0$ and $C_d>0$, depending only on $d$, such that for every $n\ge n_0$ and every $(x,y)\in T_n^2\cap \mathcal{R}_{I,J}$ with $\delta(x,y)>0$, we have
\begin{equation}\label{eq:density first estimate} p_{X_I(x),X_J(y)}(0,0) \le C_d \delta(x,y)^{-r} .\end{equation}
\item\label{lem:density:II} There exist constants $n_0$ and $C_d>0$, depending only on $d$, such that for every $n\ge n_0$ and every $(x,y)\in T_n^2\cap \mathcal{R}_{I,J}$ with $\delta(x,y)>0$ and $u_x^\top u_y>0$, we have
\begin{equation}\label{eq:density second estimate} p_{X_I(x),X_J(y)}(0,0)\le\frac{C_d}{\norm{P_{v_+^\perp}(y-x)}^{r_+}\norm{P_{v_-^\perp}(y-x)}^{r_-}}. \end{equation}
\item\label{lem:density:III} For any integers $0\le r_+,r_-\le d-2$, there exist constants $n_0$ and $C_d>0$, depending only on $d$, such that for every $n\ge n_0$,
\begin{equation}\label{eq:density angular integral}
\sup_{x\in T_n}\int_{\mathbb S^{d-1}}\frac{\dd\sigma(\omega)}{\norm{P_{v_+^\perp}\omega}^{r_+}\norm{P_{v_-^\perp}\omega}^{r_-}}\le C_d.
\end{equation}
\end{enumerate}
\end{lemma}
\begin{proof}
By rotational invariance of the standard Gaussian measurements, we may assume without loss of generality that $\xs=e$, where $e$ is the first standard basis vector. Recall the notation $\Delta(a,b)$ from Lemma~\ref{lem:probability zero degenerate covariance}.

First, by Lemma~\ref{lem:probability zero degenerate covariance}, we have
\[ \det \left(\operatorname{Cov}\begin{pmatrix} X_I(x)\\X_J(y)\end{pmatrix}\right) =\prod_{i\in I\setminus J}\norm{x-s_i^Ie}^2\prod_{j\in J\setminus I}\norm{y-s_j^Je}^2 \prod_{i\in I\cap J} \Delta^2(x-s_i^I e, y-s_i^J e). \]
Indeed, for each $i \in I\setminus J$, the entry $a_i^\top (x-s_i^I e)$ is independent of the others and has variance $\norm{x-s_i^I e}^2$. The same holds for $j\in J\setminus I$. For each $\ell\in I\cap J$, the pair of Gaussian variables $a_\ell^\top (x-s_\ell^I e)$ and $a_\ell^\top (y-s_\ell^J e)$ has mean zero and covariance matrix \[
\begin{pmatrix}
\norm{x-s_\ell^I e}^2 & (x-s_\ell^I e)^\top (y-s_\ell^J e )\\
(x-s_\ell^I e)^\top (y-s_\ell^J e ) & \norm{y-s_\ell^J e}^2
\end{pmatrix},
\]
whose determinant is $\Delta^2(x-s_\ell^I e, y-s_\ell^Je)$.
Therefore \begin{align} &p_{X_I(x),X_J(y)}(0,0)=(2\pi)^{-d} \left( \det \left(\operatorname{Cov}\begin{pmatrix} X_I(x)\\X_J(y)\end{pmatrix}\right)\right)^{-1/2} \\ =&(2\pi)^{-d} \prod_{i\in I\setminus J}\norm{x-s_i^I e}^{-1}\prod_{j\in J\setminus I}\norm{y-s_j^Je}^{-1}\prod_{i\in I\cap J} \Delta^{-1}(x-s_i^Ie, y-s_i^J e).\label{eq:density formula} \end{align}
When $n$ is large enough so that $n^{-1/2}\le \tau/2$, for any $x,y\in T_n$ and any sign $\sigma\in \{-1,1\}$, we have $\norm{x-\sigma e},\norm{y-\sigma e}\ge\tau/2$, so there exists some constant $0<C_d<\infty$ such that
\begin{equation}\label{eq:nonshared factor bound}\prod_{i\in I\setminus J}\norm{x-s_i^I e}^{-1}\prod_{j\in J\setminus I}\norm{y-s_j^Je}^{-1} \le C_d.\end{equation}

For $x\in T_n$, we can write its tube parametrization \[
x=\alpha_x e + (\tau+\beta_x) u_x \qquad \text{with }u_x\in \mathbb{S}^{d-2},\quad\text{and }\abs{\alpha_x},\abs{\beta_x} \le \frac{\tau}{2}.
\]
Then for any $\ell\in I\cap J$, a direct expansion implies \begin{align}
\Delta^2(x-s_\ell^I e, y-s_\ell^Je) &=(\tau+\beta_x)^2 (\tau+\beta_y)^2 \left[1-(u_x^\top u_y)^2\right] \notag\\
&\quad+\norm{(\alpha_x-s_\ell^I)(\tau+\beta_y)u_y-(\alpha_y-s_\ell^J)(\tau+\beta_x)u_x}^2 \\
&\ge (\tau+\beta_x)^2 (\tau+\beta_y)^2 \left[1-(u_x^\top u_y)^2\right] \\
&\ge \frac{\tau^4}{16} \left[1-(u_x^\top u_y)^2\right]  \\
&=\frac{\tau^4}{16} \left(1-\left(1-\frac{\delta(x,y)^2}{2}\right)^2\right) \label{eq:identity for delta}\\
&\ge \frac{\tau^4}{32} \delta(x,y)^2\label{eq:shared factor first bound},
\end{align}
where the equality~\eqref{eq:identity for delta} follows from \[
\delta(x,y)^2 = \min\{\norm{u_x-u_y}^2,\norm{u_x+u_y}^2\}=\min\{2-2u_x^\top u_y, 2+2u_x^\top u_y\} =2(1-\abs{u_x^\top u_y}).
\]
Substituting~\eqref{eq:nonshared factor bound} and~\eqref{eq:shared factor first bound} into the density formula~\eqref{eq:density formula} proves statement (\ref{lem:density:I}).

A second estimate of $\Delta^2(x-s_\ell^I e, y-s_\ell^Je)$ yields~\eqref{eq:density second estimate}. For matching shared signs $s_i^I=s_i^J=\sigma$, we have \begin{equation}\label{eq:density matching sign bound} \Delta(x-\sigma e, y-\sigma e)=\Delta(x-\sigma e, y-x)=\norm{x-\sigma e}\norm{P_{(x-\sigma e)^\perp}(y-x)}. \end{equation}
For opposite shared signs, first consider the corresponding points on the saddle sphere $\mathcal{S}$. With $t=u_x^\top u_y>0$, a direct expansion implies \[ \Delta^2(\tau u_x-\sigma e, \tau u_y+\sigma e)=(1+\tau^2)^2-(\tau^2t-1)^2=\tau^2(1+t)\bigl[2+\tau^2(1-t)\bigr]\ge2\tau^2. \]
Since $\norm{x-\tau u_x}+\norm{y-\tau u_y}\le C_dn^{-1/2}$ for some constant $C_d$, by uniform continuity of $(x,y)\mapsto\Delta^2(x-\sigma e,y+\sigma e)$ near $\mathcal{S}$, there exists $n_0$, depending only on $d$, such that when $n\ge n_0$, \begin{equation}\label{eq:density opposite sign bound} \Delta^2(x-\sigma e, y+\sigma e) \ge \tau^2. \end{equation}
Thus opposite-sign factors are uniformly bounded away from zero. Substituting~\eqref{eq:nonshared factor bound},~\eqref{eq:density matching sign bound}, and~\eqref{eq:density opposite sign bound} into the density formula~\eqref{eq:density formula} proves statement (\ref{lem:density:II}).

To prove (\ref{lem:density:III}), fix $x\in T_n$ and recall that $v_\pm=x\mp e$. For $\omega\in\mathbb S^{d-1}$, the formula for orthogonal projection implies
\[
\norm{P_{v_\pm^\perp}\omega}^2 =\left\|\omega-\frac{v_\pm^\top\omega}{\norm{v_\pm}^2}v_\pm\right\|^2 =1-\frac{(v_\pm^\top\omega)^2}{\norm{v_\pm}^2} =\sin^2\theta_\pm(\omega),
\]
where
\[
\theta_\pm(\omega)
=\arccos\left(\frac{v_\pm^\top\omega}{\norm{v_\pm}}\right)\in[0,\pi]
\]
is the angle between $\omega$ and $v_\pm$. Since sine is nonnegative on $[0,\pi]$, we have
\[
\norm{P_{v_\pm^\perp}\omega}=\sin\theta_\pm(\omega).
\]
Thus a denominator factor becomes small only when its angle is close to $0$ or $\pi$. We divide the integral according to which of these two factors is small.

First, the directions of $v_+$ and $v_-$ remain separated uniformly in $x$. Indeed, using $x=\alpha_xe+(\tau+\beta_x)u_x$, we obtain
\[
\Delta(v_+,v_-)^2
=4\bigl(\norm{x}^2-(x^\top e)^2\bigr)
=4(\tau+\beta_x)^2.
\]
Since $\tau+\beta_x\ge\tau/2$ for sufficiently large $n$, and the norms $\norm{x\mp e}$ are uniformly bounded on $T_n$, there is a fixed $c_0\in(0,1)$, depending only on $d$, such that
\begin{equation}\label{eq:density angular separation}
\frac{\Delta(v_+,v_-)}{\norm{v_+}\norm{v_-}}
=\frac{2(\tau+\beta_x)}{\norm{x-e}\norm{x+e}}
\ge c_0.
\end{equation}
The left side is the sine of the angle between the two vectors. Hence the angle between either signed direction of $v_+$ and either signed direction of $v_-$ is at least $\arcsin c_0$.

Choose a fixed $0<\theta_0<\frac13\arcsin c_0$ and define
\[
\begin{aligned}
E_+&=\{\omega\in\mathbb S^{d-1}:\theta_+(\omega)<\theta_0
\text{ or }\theta_+(\omega)>\pi-\theta_0\},\\
E_-&=\{\omega\in\mathbb S^{d-1}:\theta_-(\omega)<\theta_0
\text{ or }\theta_-(\omega)>\pi-\theta_0\},\quad\text{and}\\
E_0&=\mathbb S^{d-1}\setminus(E_+\cup E_-).
\end{aligned}
\]
The separation bound implies $E_+\cap E_-=\varnothing$. Thus these three regions partition the sphere: $E_+$ contains directions close to the line spanned by $v_+$, $E_-$ contains directions close to the line spanned by $v_-$, and $E_0$ contains the remaining directions. Consequently,
\[
\int_{\mathbb S^{d-1}}
\frac{\dd\sigma(\omega)}{\norm{P_{v_+^\perp}\omega}^{r_+}\norm{P_{v_-^\perp}\omega}^{r_-}}
=\left(\int_{E_+}+\int_{E_-}+\int_{E_0}\right)
\frac{\dd\sigma(\omega)}{(\sin\theta_+(\omega))^{r_+}(\sin\theta_-(\omega))^{r_-}}.
\]

On $E_+$, disjointness from $E_-$ implies $\theta_-(\omega)\in[\theta_0,\pi-\theta_0]$, so $\sin\theta_-(\omega)\ge\sin\theta_0$. To integrate the other factor, use spherical coordinates around $v_+$:
\[
\omega=\cos\theta\,\frac{v_+}{\norm{v_+}}+\sin\theta\,w
\qquad\text{with }0<\theta<\pi,\quad\text{and }w\in\mathbb S^{d-1}\cap v_+^\perp.
\]
Here $\theta=\theta_+(\omega)$, the sphere containing $w$ is identified with $\mathbb S^{d-2}$, and
\[
\dd\sigma(\omega)=(\sin\theta)^{d-2}\dd\theta\,\dd\sigma(w).
\]
The region $E_+$ corresponds to $\theta\in[0,\theta_0)\cup(\pi-\theta_0,\pi]$. Since $\sin(\pi-\theta)=\sin\theta$, the two angular intervals contribute equally. Therefore,
\begin{align*}
\int_{E_+}
\frac{\dd\sigma(\omega)}{(\sin\theta_+(\omega))^{r_+}(\sin\theta_-(\omega))^{r_-}}
&\le (\sin\theta_0)^{-r_-}\int_{E_+}
(\sin\theta_+(\omega))^{-r_+}\dd\sigma(\omega)\\
&=2(\sin\theta_0)^{-r_-}\int_0^{\theta_0}\int_{\mathbb S^{d-2}}
(\sin\theta)^{d-2-r_+}\dd\sigma(w)\,\dd\theta\\
&=2(\sin\theta_0)^{-r_-}\abs{\mathbb S^{d-2}}
\int_0^{\theta_0}(\sin\theta)^{d-2-r_+}\dd\theta
\le C_d.
\end{align*}
The last inequality uses $d-2-r_+\ge0$ and $r_-\le d-2$. Interchanging $+$ and $-$ implies the same bound for the integral over $E_-$.

Finally, on $E_0$, both angles belong to $[\theta_0,\pi-\theta_0]$. Hence
\[
\int_{E_0}
\frac{\dd\sigma(\omega)}{(\sin\theta_+(\omega))^{r_+}(\sin\theta_-(\omega))^{r_-}}
\le\abs{\mathbb S^{d-1}}(\sin\theta_0)^{-(r_++r_-)}\le C_d.
\]
Adding the three bounds proves~\eqref{eq:density angular integral}. The choice of $\theta_0$ and all constants are independent of $x$ and $n$ for sufficiently large $n$, which proves (\ref{lem:density:III}).
\end{proof}

\subsection{Kac--Rice formulas}\label{app:kac-rice-tools}

We first state a weighted Kac--Rice formula for Gaussian random fields and then apply it to the active equations defining candidate minima.

\begin{lemma}[Weighted Kac--Rice formula {\cite[Theorem 6.4]{AzaisWschebor2009}}]\label{core lemma:Kac-Rice}
Let $V\subset \R^p$ be open, let $Z:V\to \R^p$ be a Gaussian random field with almost surely $C^1$ sample paths, and let $\Theta\in \R^q$ be a random vector. Suppose the following conditions hold.
\begin{enumerate}[label=\textnormal{(K\arabic*)},ref=K\arabic*]
\item\label{core lemma:Kac-Rice:H1}
The random vector $\Theta$ and the field $Z$ are jointly Gaussian, i.e., for every $n\ge 1$ and every $x_1,\ldots,x_n\in V$, the vector \( \bigl(\Theta,Z(x_1),\ldots, Z(x_n)\bigr) \in \R^{q+np} \) is Gaussian.
\item\label{core lemma:Kac-Rice:H2}
For every $x\in V$, we have
\(
\operatorname{Cov}(Z(x)) \succ 0.
\)
\item\label{core lemma:Kac-Rice:H3}
We have
\( \PP(\text{There exists an }x\in V\text{ such that }Z(x)=0\text{ and }\det DZ(x)=0 )=0.\)
\end{enumerate}
Then, for every Borel set $B\subset V$ and every nonnegative Borel function $H: V\times \R^q \to [0,\infty]$, we have \[
\E\left[\sum_{\substack {x\in B\\ Z(x)=0}} H(x,\Theta) \right] =\int_B p_{Z(x)}(0) \E\bigl[\abs{\det DZ(x)} H(x,\Theta)\bigm|Z(x)=0 \bigr] \, \dd x,
\]
where the conditional expectation uses the usual Gaussian conditional law.
\end{lemma}
\begin{proof}
The original theorem~\cite[Theorem 6.4]{AzaisWschebor2009} holds for every compact set $B$ and bounded continuous function $H$ under the same conditions. We first fix a compact set $B$; a standard monotone-class argument and the monotone convergence theorem then extend the identity to nonnegative Borel functions.

To extend the identity to a Borel set $C\subset V$, choose compact sets $K_n\subset V$ with $K_n\uparrow V$. By (\ref{core lemma:Kac-Rice:H3}), the Jacobian $DZ(x)$ is invertible at every zero of $Z$, almost surely. Since the sample paths are $C^1$, the inverse function theorem implies that every zero is isolated. The zero set is also closed in $V$ by continuity, so each compact set $K_n$ contains only finitely many zeros. Since the sets $K_n$ exhaust $V$, the full zero set is at most countable.
The sums in the formula are taken over these zeros. As $K_n$ increases, more zeros are included, and the nonnegativity of $H$ implies that, almost surely,
\[
\sum_{\substack{x\in K_n\\Z(x)=0}} 1_{x\in C}H(x,\Theta)
\rightarrow
\sum_{\substack{x\in C\\Z(x)=0}} H(x,\Theta).
\]
Apply the compact-set Kac--Rice identity on $K_n$ to the nonnegative Borel function $(x,\theta)\mapsto 1_{x\in C}H(x,\theta)$. Monotone convergence allows us to pass to the limit both in the expectation of the sum and in the integral on the right-hand side, establishing the formula for $C$.
\end{proof}

For an active set $I\subset[n]$ with $\abs I=d$, define
\[
H_I(x;A)=1_{\mathcal K_x}(A_I)1_{-S(x;A_{I^c})\in Q(x;A_I)(-1,1)^d},
\]
using the notation from Appendix~\ref{sec:common notation}. For $x\ne\pm\xs$, let $p_{X_I(x)}$ denote the Gaussian density of $X_I(x)$; the joint density and its domain $\mathcal R_{I,J}$ were defined above.

\begin{lemma}[Kac--Rice formulae]\label{lem:Kac-Rice for first and second moment}
Fix integers $d\ge4$ and $n\ge d$, and a deterministic unit vector $\xs\in\R^d$. Let $a_1,\ldots,a_n$ be independent $N(0,I_d)$ vectors, and suppose $n^{-1/2}<\tau$. Then the following identities hold.
\begin{enumerate}[label=\textnormal{(\Roman*)},ref=\Roman*]
\item\label{lem:Kac-Rice for first and second moment:I}  For any active set $I$ with $\abs{I}=d$, we have
\[
\E[\chi_I^n]=\int_{T_n} p_{X_I(x)}(0)\E\bigl[\abs{\det A_I}H_I(x;A)\bigm|X_I(x)=0 \bigr] \dd x.
\]
\item\label{lem:Kac-Rice for first and second moment:II} For any active sets $I\neq J$ with $\abs{I}=\abs{J}=d$, we have
\[
% Keep the complete integral on one line in both formats.
\ificlr
\resizebox{\linewidth}{!}{$\displaystyle
\begin{aligned}
&\E[\chi_I^n\chi_J^n]\\
&\quad=\int\limits_{T_n^2\cap\mathcal R_{I,J}}p_{X_I(x),X_J(y)}(0,0)\,\E\bigl[\abs{\det A_I}\abs{\det A_J}H_I(x;A)H_J(y;A)\bigm|X_I(x)=X_J(y)=0\bigr]\dd x\,\dd y.
\end{aligned}
$}
\else
\resizebox{\linewidth}{!}{$\displaystyle
\E[\chi_I^n\chi_J^n]=\int\limits_{T_n^2\cap\mathcal R_{I,J}}p_{X_I(x),X_J(y)}(0,0)\,\E\bigl[\abs{\det A_I}\abs{\det A_J}H_I(x;A)H_J(y;A)\bigm|X_I(x)=X_J(y)=0\bigr]\dd x\,\dd y.
$}
\fi
\]
\end{enumerate}
\end{lemma}
\begin{proof}
By rotational invariance of the standard Gaussian measurements, we may assume without loss of generality that $\xs=e$, where $e$ is the first standard basis vector.
We apply the weighted Kac--Rice formula in Lemma~\ref{core lemma:Kac-Rice}, first to $X_I$ and then to the joint field $(X_I,X_J)$. In both applications, take $\Theta$ to be the vector of all entries of $A$. For each fixed realization of $A$, the maps $x\mapsto X_I(x)$ and $(x,y)\mapsto(X_I(x),X_J(y))$ are affine, so their sample paths are smooth. At fixed inputs, their values are linear in $\Theta$, so both fields are jointly Gaussian with $\Theta$. It is also easy to verify that $H_I$ is a nonnegative Borel function.

For statement (\ref{lem:Kac-Rice for first and second moment:I}), consider $X_I$ on the open set $V=\R^d\setminus\{\pm e\}$. Write $I=\{i_1<\cdots<i_d\}$. The coordinates of $X_I(x)$ are independent Gaussian variables with mean zero, and hence
\[
\operatorname{Cov}(X_I(x))
=\operatorname{diag}\bigl(\norm{x-s_1e}^2,\ldots,\norm{x-s_de}^2\bigr)
\succ0
\qquad\text{for }x\in V.
\]
Moreover, $DX_I(x)=A_I$, which is almost surely nonsingular. Thus all three hypotheses of Lemma~\ref{core lemma:Kac-Rice} hold. The assumption $n^{-1/2}<\tau$ ensures that $T_n\subset V$. Since $X_I$ has the unique zero $x_I$, we have
\[
\sum_{\substack{x\in T_n\\X_I(x)=0}}H_I(x;A)
=1_{\{x_I\in T_n\}}H_I(x_I;A)
=\chi_I^n.
\]
Applying the formula with the Borel set $T_n$ proves statement (\ref{lem:Kac-Rice for first and second moment:I}).

For statement (\ref{lem:Kac-Rice for first and second moment:II}), consider the joint field on $\mathcal R_{I,J}$. Since covariance matrices are positive semidefinite, $\mathcal R_{I,J}$ is precisely the set where the covariance determinant is strictly positive. This determinant depends continuously on $(x,y)$, so $\mathcal R_{I,J}$ is open. The derivative of the joint field is
\[
D_{(x,y)}
\begin{pmatrix}X_I(x)\\X_J(y)\end{pmatrix}
=
\begin{pmatrix}A_I&0\\0&A_J\end{pmatrix}.
\]
Both active matrices are almost surely nonsingular, so the Jacobian is nonsingular and its absolute determinant is $\abs{\det A_I}\abs{\det A_J}$. Apply Lemma~\ref{core lemma:Kac-Rice} with the weight $H_I(x;A)H_J(y;A)$ and the Borel set $T_n^2\cap\mathcal R_{I,J}$. The right-hand side is exactly the integral in statement (\ref{lem:Kac-Rice for first and second moment:II}). Since the joint field has the unique zero $(x_I,x_J)$, the left-hand side is
\[
\E\left[\chi_I^n\chi_J^n
1_{\{(x_I,x_J)\in\mathcal R_{I,J}\}}\right].
\]
Lemma~\ref{lem:probability zero degenerate covariance} implies $(x_I,x_J)\in\mathcal R_{I,J}$ almost surely. Removing this indicator proves statement (\ref{lem:Kac-Rice for first and second moment:II}).
\end{proof}

\subsection{Uniform local limit theorems}\label{app:local-limit-tools}

We state a local limit theorem for parallelepipeds, derive a version uniform over a compact parameter set, and apply it to the probabilities appearing in the first- and second-moment formulas.
The following is a direct consequence of \cite[Theorem~1.I]{BorovkovMogulskii2006}. We state it without proof.
\begin{lemma}[Parallelepiped version of {\cite[Theorem~1.I]{BorovkovMogulskii2006}}]\label{literature local limit}
Let $Y_{n,1},\ldots,Y_{n,n}$ be independent copies of $Y_n$, with
\(
\E Y_n=0,\) and  \(\operatorname{Cov}(Y_n)=I_p.
\) Define $T_n=\sum_{j=1}^nY_{n,j}$, and
suppose the following  holds.
\begin{enumerate}
\item[\textnormal{(L)}]\namedlabel{L}{literature local limit:L} For every $\epsilon>0$, we have
\(
\lim_{n\to\infty} \E\left[ \|Y_n\|^2 1_{\{\|Y_n\|>\epsilon\sqrt n\}} \right] = 0. \)

\item[\textnormal{(R)}]\namedlabel{R}{literature local limit:R} For every $0<r<R<\infty$, we have
\(
\limsup_{n\to\infty}
\sup_{r\le\|\xi\|\le R}
\left|\mathbb Ee^{i\xi^\top Y_n}\right|<1. \)
\end{enumerate}
Then, for every compact set $\mathcal{B}\subset GL_p(\R)$, we have
\[
\lim_{n\to \infty} \sup_{\substack{a\in\mathbb R^p\\B\in\mathcal{B}}}
\left| \frac{n^{p/2}}{2^p\abs{\det B}}  \PP\left(T_n\in a+B(-1,1)^p\right) -\varphi_{\scriptscriptstyle I_p}\left(\frac{a}{\sqrt n}\right) \right|= 0,
\]
where $\varphi_{\scriptscriptstyle I_p}$ is the density of $N(0,I_p)$.
\end{lemma}

\begin{lemma}[Uniform local limit theorem over a compact parameter set]\label{core lemma:local limit theorem}
Let $K\subset \R^q$ be a nonempty compact set, and let $\{W_t: t\in K\}$ be a family of random vectors in $\R^p$. For each $t\in K$, let $W_{t,1},W_{t,2},\ldots$ be independent copies of $W_t$, and define \[
\mu_t = \E[W_t],\qquad \Sigma_t=\operatorname{Cov}(W_t)\qquad\text{and}\qquad S_{n,t}=\sum_{j=1}^n W_{t,j}.
\]
Suppose that the following hold true.
\begin{enumerate}[label=\textnormal{(H\arabic*)},ref=H\arabic*]
\item\label{cond:ullt-weak-continuity}
The family is weakly continuous, i.e., for every sequence $t_n\to t$ in $K$, we have \( W_{t_n}\xrightarrow{d} W_t.\)
\item\label{cond:ullt-third-moment}
The third moments are uniformly bounded, i.e., \(\sup_{t\in K} \E[\norm{W_t}^3]<\infty. \)
\item\label{cond:ullt-covariance}
The covariance matrices are uniformly positive definite, i.e., \( \inf_{t\in K} \lambda_{\min}(\Sigma_t)>0.\)
\item\label{cond:ullt-nonlattice}
For every $t\in K$ and every $\xi\in \R^p \setminus \{0\}$, we have \( \abs{\E[e^{i\xi^\top W_t}]}<1.\)
\end{enumerate}
Then, for every compact set $\mathcal{B}\subset GL_p(\R)$, we have
\begin{equation}\label{eq:corollary local limit}
\lim_{n\to \infty} \sup_{\substack{t \in K, a\in \R^p\\B\in \mathcal{B}}} \abs{\frac{n^{p/2}}{2^p\abs{\det B}}\PP(S_{n,t}\in a+B(-1,1)^p)-\varphi_{\scriptscriptstyle \Sigma_t}(\frac{a-n\mu_t}{\sqrt{n}})}=0,
\end{equation}
where $\varphi_{\Sigma_t}$ is the density of $N(0,\Sigma_t)$.
\end{lemma}
\begin{proof}
Suppose, for contradiction, that the limit \eqref{eq:corollary local limit} is false. Then, there exist $\epsilon>0$ and sequences \(
n_k\uparrow \infty\), \(t_k\in K\), \(a_k\in \R^p,\) and \(B_k\in\mathcal{B}\)
such that \begin{equation}\label{eq:contradiction}
\left| \frac{n_k^{p/2}} {2^p|\det B_k |} \PP\!\left( S_{n_k,t_k}\in a_k+B_k(-1,1)^p \right) - \varphi_{\scriptscriptstyle \Sigma_{t_k}} \left(\frac {a_k-n_k\mu_{t_k}}{\sqrt{n_k}}\right) \right| \ge \epsilon.
\end{equation}
Because $K$ is compact, passing to a subsequence, we may assume $t_k\to t_\star$ for some $t_{\star}\in K$. 
Condition (\ref{cond:ullt-third-moment}) implies uniform integrability of both $\{W_t\}_{t\in K}$ and $\{\norm{W_t}^2\}_{t\in K}$. Combining uniform integrability with (\ref{cond:ullt-weak-continuity}) yields $
\mu_{t_k} \rightarrow \mu_{t_\star}$ and $\Sigma_{t_k} \rightarrow \Sigma_{t_{\star}}.$ In particular, the map $t\mapsto \Sigma_t$ is continuous on $K$. Since $K$ is compact, we have \[
\sup_{t\in K} \lambda_{\max}(\Sigma_t) <\infty.
\]
Together with (\ref{cond:ullt-covariance}), this yields uniform bounds on the spectrum of $\Sigma_t$. In turn, the map $t\mapsto \Sigma_t^{-1/2}$ is continuous and uniformly bounded on $K$.

Define the random variables \[
Y_k=\Sigma_{t_k}^{-1/2} (W_{t_k}-\mu_{t_k}),\qquad \text{and} \qquad T_k = \sum_{j=1}^{n_k} Y_{k,j}.
\]
Clearly $\E Y_k=0$ and $\Cov(Y_k)=I_p$. The uniform bound on $\Sigma_t^{-1/2}$ and condition~(\ref{cond:ullt-third-moment}) imply
\(
\sup_{t\in K}\E\bigl\|\Sigma_t^{-1/2}(W_t-\mu_t)\bigr\|^3<\infty.
\)
Hence, for every $\epsilon>0$,
\[
\E\!\left[\|Y_k\|^2 1_{\{\|Y_k\|>\epsilon\sqrt{n_k}\}}\right]
\le \frac{\sup_{t\in K}\E\|\Sigma_t^{-1/2}(W_t-\mu_t)\|^3}{\epsilon\sqrt{n_k}}
\longrightarrow0,
\]
which verifies condition (\ref{literature local limit:L}) at the sample sizes $n_k$. For (\ref{literature local limit:R}), condition~(\ref{cond:ullt-weak-continuity}) and the continuity of $\mu_t$ and $\Sigma_t$ imply that
\[
(t,\xi)\longmapsto\left|\E e^{i\xi^\top\Sigma_t^{-1/2}(W_t-\mu_t)}\right|
\]
is continuous. By condition~(\ref{cond:ullt-nonlattice}) and the invertibility of $\Sigma_t$, this function is strictly less than $1$ whenever $\xi\ne0$. Therefore, for every $0<r<R<\infty$,
\[
\sup_k\sup_{r\le\|\xi\|\le R}\left|\E e^{i\xi^\top Y_k}\right|
\le\max_{\substack{t\in K\\r\le\|\xi\|\le R}}
\left|\E e^{i\xi^\top\Sigma_t^{-1/2}(W_t-\mu_t)}\right|<1.
\]
These estimates hold uniformly over $t\in K$, so Lemma~\ref{literature local limit} applies along the sample sizes $n_k$. Finally, the set \(
\mathcal{C}=\left\{ \Sigma_t^{-1/2} B: t\in K, B\in \mathcal{B} \right\} 
\)
is a compact subset of $GL_p(\R)$. Define
\[
b_k:=\Sigma_{t_k}^{-1/2}(a_k-n_k\mu_{t_k}),\qquad\text{and}\qquad
C_k:=\Sigma_{t_k}^{-1/2}B_k\in\mathcal{C}.
\]
Applying Lemma~\ref{literature local limit} to $Y_k, T_k$ and the compact set $\mathcal{C}$, and using the relation \[
\PP(S_{n_k,t_k}\in a_k+B_k (-1,1)^p)=\PP(T_k\in b_k +C_k(-1,1)^p),
\]
we have \[
\lim_{k\to\infty} \left| \frac{n_k^{p/2}} {2^p|\det B_k |} \PP\!\left( S_{n_k,t_k}\in a_k+B_k(-1,1)^p \right) - \varphi_{\scriptscriptstyle \Sigma_{t_k}} \left(\frac {a_k-n_k\mu_{t_k}}{\sqrt{n_k}}\right) \right| =0,
\]
contradicting \eqref{eq:contradiction}. Hence the result follows.
\end{proof}

We use the notation from Appendix~\ref{sec:common notation}, and write $\varphi_\Sigma$ for the density of $N(0,\Sigma)$. For the one-point probability, let $n>d$ and $C=(g_1,\ldots,g_{n-d})^\top$, where the rows are independent $N(0,I_d)$ vectors. For $x\in\R^d$ and a deterministic matrix $B\in\R^{d\times d}$, define
\[
\pi_n(x;B)=\PP_{\scriptscriptstyle g_1,\ldots,g_{n-d}}(-S(x;C)\in Q(x;B)(-1,1)^d).
\]
For the joint probability, take distinct active sets $I,J\subset[n]$ with $\abs I=\abs J=d$, and write $r=\abs{I\cap J}$ and $\nu=n-2d+r$. When $\nu\ge1$, let $C=(g_1,\ldots,g_\nu)^\top$ have independent $N(0,I_d)$ rows. For a deterministic matrix $B\in\R^{(2d-r)\times d}$ whose rows are indexed by $I\cup J$, write $B_L$ for its submatrix with row labels in $L$. For $x,y\in\R^d$, define
\[
\pi_{n,I,J}(x,y;B)=\PP_{\scriptscriptstyle g_1,\ldots,g_\nu}\left\{\begin{aligned}
-S(x;C)&\in Q(x;B_I)(-1,1)^d+S(x;B_{J\setminus I}),\\
-S(y;C)&\in Q(y;B_J)(-1,1)^d+S(y;B_{I\setminus J})
\end{aligned}\right\}.
\]
We use the angular separation $\delta(x,y)$ from the proof of Theorem~\ref{thm:second moment asymptotic appendix}.

\begin{lemma}[Local limit theorem]\label{lem:limit theorem for first and second moment}
Fix an integer $d\ge4$ and a deterministic unit vector $\xs\in\R^d$. The following limits hold.
\begin{enumerate}[label=\textnormal{(\Roman*)},ref=\Roman*]
\item\label{lem:limit theorem for first and second moment:I} We have
\[
\lim_{n\to \infty} \sup_{x\in T_n} \sup_{B\in \mathcal{K}_x} \abs{(n-d)^{d/2}\pi_n(x;B)-2^d \abs{\det Q(x;B)}\varphi_{\scriptscriptstyle \Sigma_x}(-\sqrt{n-d}\mu(x))} =0.
\]
\item\label{lem:limit theorem for first and second moment:II} Fix $0<\delta_0<1$ and an integer $0\le r\le d-1$. For index sets $I,J\subset[n]$ with $\abs I=\abs J=d$ and $\abs{I\cap J}=r$, we have
\[
\lim_{n\to\infty} \sup_{\substack{x,y \in T_n \\ \delta(x,y)\ge \delta_0}} \sup_{\substack{B\colon B_I\in \mathcal{K}_x \\B_J\in \mathcal{K}_y}} \left|\begin{aligned}
&\nu^d \pi_{n,I,J}(x,y;B)\\
&-2^{2d}\abs{\det Q(x;B_I)}\abs{\det Q(y;B_J)}\varphi_{\scriptscriptstyle \Sigma_{x,y}}\left(-\sqrt{\nu} \begin{pmatrix} \mu(x)\\ \mu(y)\end{pmatrix}\right)
\end{aligned}\right| = 0.
\]
\end{enumerate}
\end{lemma}
\begin{proof}
By rotational invariance, we may assume without loss of generality that $\xs=e$, where $e$ is the first standard basis vector.
Recall that
\[
\Sigma_x=\operatorname{Cov}(Z_x(a)),\qquad\text{and}\qquad
\Sigma_{x,y}=\operatorname{Cov}\begin{pmatrix}Z_x(a)\\Z_y(a)\end{pmatrix},
\]
where both components of the joint vector use the same $a\sim N(0,I_d)$.
We use the uniform local limit theorem in Lemma~\ref{core lemma:local limit theorem}.

To prove (\ref{lem:limit theorem for first and second moment:I}), choose a fixed integer $n_0$ sufficiently large that $n_0^{-1/2}<\tau/2$ and the bound in Lemma~\ref{lem:tube}(\ref{lem:tube:II}) holds on $T_{n_0}$, and set $K=T_{n_0}$. This compact set contains $T_n$ for all $n\ge n_0$ and excludes $\pm e$. In Lemma~\ref{core lemma:local limit theorem}, we take $W_x=Z_x(a)\in\R^d$ for $x\in K$.
For any sequence $x_j\to x$ in $K$, the sign in $Z_{x_j}(a)$ eventually agrees with that in $Z_x(a)$ unless $a^\top(x-e)=0$ or $a^\top(x+e)=0$. Since $x\ne\pm e$, the Gaussian vector $a$ lies on one of these hyperplanes with probability zero. Thus $Z_{x_j}(a)\to Z_x(a)$ almost surely and hence in distribution, verifying condition~(\ref{cond:ullt-weak-continuity}).
Moreover, since $K$ is bounded,
\[
\norm{Z_x(a)}\le\norm{x}\norm{a}^2\le C_d\norm{a}^2
\qquad\text{for }x\in K.
\]
Since $\E\norm{a}^6<\infty$, this bound implies $\sup_{x\in K}\E\norm{Z_x(a)}^3<\infty$, verifying condition~(\ref{cond:ullt-third-moment}).
Lemma~\ref{lem:tube}(\ref{lem:tube:II}) implies
\(
\inf_{x\in K}\lambda_{\min}(\Sigma_x)\ge c_d>0,
\)
verifying condition~(\ref{cond:ullt-covariance}).

We next verify condition~(\ref{cond:ullt-nonlattice}). Seeking a contradiction, suppose that
\(
|\E \exp(i\xi^\top Z_x(a))|=1
\)
for some $\xi\ne0$. Since $e^{i\xi^\top Z_x(a)}$ has modulus one, equality in the triangle inequality implies that it equals $e^{i\theta}$ almost surely for some fixed $\theta\in\R$.
Choose an open ball $U$ that does not intersect either hyperplane $a^\top(x-e)=0$ or $a^\top(x+e)=0$. On $U$, the sign in $Z_x(a)$ is fixed, so
\[
\xi^\top Z_x(a)=\sigma(a^\top x)(a^\top\xi)
\qquad\text{for }a\in U\text{ and some fixed }\sigma\in\{-1,1\}.
\]
The Gaussian density is positive everywhere, and $e^{i\xi^\top Z_x(a)}$ is continuous on $U$. Thus its almost sure equality to $e^{i\theta}$ implies equality at every point of $U$. Consequently, \[\sigma(a^\top x)(a^\top\xi)\in\theta+2\pi\mathbb Z\qquad\text{for all }a\in U.\] Since this polynomial is continuous and $U$ is connected, it cannot move between these discrete values and must be constant on $U$.
A polynomial that is constant on an open ball is constant on all of $\R^d$. Evaluating at $a=0$ shows that this constant is zero, so $(a^\top x)(a^\top\xi)=0$ for every $a\in\R^d$. This is impossible when both $x$ and $\xi$ are nonzero: one can choose $a$ outside the two hyperplanes $a^\top x=0$ and $a^\top\xi=0$. This contradiction proves~(\ref{cond:ullt-nonlattice}).
Thus all four conditions of Lemma~\ref{core lemma:local limit theorem} hold on $K$.

To apply the conclusion of Lemma~\ref{core lemma:local limit theorem} to $\pi_n(x;B)$ uniformly over $B\in\mathcal K_x$, we need the matrices $Q(x;B)$ to lie in a fixed compact subset of $GL_d(\R)$. Consider
\[
\mathcal B_1=\{Q(x;B):x\in K,\ B\in\mathcal K_x\}.
\]
The pairs $(x,B)$ in this definition form a compact set, since $K$ is compact, the rows of $B$ are bounded, and the inequalities defining $\mathcal K_x$ are nonstrict and depend continuously on $(x,B)$. The map $(x,B)\mapsto Q(x;B)$ is continuous, and $\abs{\det Q(x;B)}\ge\lambda_d$ for $B\in\mathcal K_x$. Hence $\mathcal B_1$ is a compact subset of $GL_d(\R)$; in particular, its determinants are uniformly bounded.

The set $Q(x;B)(-1,1)^d$ is symmetric, so
\[
\pi_n(x;B)=\PP\bigl(S(x;C)\in Q(x;B)(-1,1)^d\bigr).
\]
Applying Lemma~\ref{core lemma:local limit theorem} with translation zero, and multiplying its error by the uniformly bounded factor $2^d\abs{\det Q(x;B)}$, yields
\[
(n-d)^{d/2}\pi_n(x;B)
=2^d\abs{\det Q(x;B)}
\varphi_{\scriptscriptstyle\Sigma_x}(-\sqrt{n-d}\,\mu(x))+o(1),
\]
uniformly for $x\in K$ and $B\in\mathcal K_x$. Since $T_n\subset K$ for $n\ge n_0$ and $n-d\to\infty$, this proves (\ref{lem:limit theorem for first and second moment:I}).

For (\ref{lem:limit theorem for first and second moment:II}), fix $\delta_0$ and $r$ as in the statement. Choose $n_1\ge n_0$ sufficiently large, as specified below, and consider
\[
K_{\delta_0}=\{(x,y)\in T_{n_1}\times T_{n_1}:\delta(x,y)\ge\delta_0\}.
\]
For $x\in K$, the orthogonal projection of $x$ onto $e^\perp$ is nonzero, so its normalization $u_x$ depends continuously on $x$. Hence $\delta(x,y)$ is continuous on $K\times K$, and $K_{\delta_0}$ is a closed subset of the compact set $K\times K$.
For this application of Lemma~\ref{core lemma:local limit theorem}, we take
\[
W_{x,y}=\begin{pmatrix}Z_x(a)\\Z_y(a)\end{pmatrix}\in\R^{2d}\qquad\text{for }(x,y)\in K_{\delta_0}.
\]
The same convergence and domination arguments verify conditions~(\ref{cond:ullt-weak-continuity}) and~(\ref{cond:ullt-third-moment}).
The set $\{(u,v)\in(\mathbb S^{d-2})^2:\min\{\norm{u-v},\norm{u+v}\}\ge\delta_0\}$ is compact. Lemma~\ref{lem:tube}(\ref{lem:tube:III}) therefore allows $n_1$ to be chosen so that
\[
\inf_{(x,y)\in K_{\delta_0}}\lambda_{\min}(\Sigma_{x,y})
\ge c_{d,\delta_0}>0.
\]
This verifies condition~(\ref{cond:ullt-covariance}).

To verify condition~(\ref{cond:ullt-nonlattice}), fix $(x,y)\in K_{\delta_0}$ and suppose
\[
\abs{\E e^{i[\xi^\top Z_x(a)+\eta^\top Z_y(a)]}}=1
\]
for some $(\xi,\eta)\ne(0,0)$. As in (\ref{lem:limit theorem for first and second moment:I}), the random variable $e^{i[\xi^\top Z_x(a)+\eta^\top Z_y(a)]}$ must be almost surely constant.
Since $\delta(x,y)>0$, we have $u_y\ne\pm u_x$, so $e,x,y$ are linearly independent. The linear map $a\mapsto(a^\top e,a^\top x,a^\top y)$ is therefore onto $\R^3$: we can prescribe these three inner products independently. In particular, we can choose $a$ so that
\[
(a^\top e,a^\top x,a^\top y)=(1,2,2)
\qquad\text{or}\qquad (1,2,0).
\]
The residual pair $\big((a^\top x)^2-(a^\top e)^2,(a^\top y)^2-(a^\top e)^2\big)$ is $(3,3)$ for the first choice and $(3,-1)$ for the second. By continuity, the signs of its entries remain $(1,1)$ on an open ball around the first choice of $a$, and $(1,-1)$ on an open ball around the second.
On each ball, $\xi^\top Z_x(a)+\eta^\top Z_y(a)$ agrees with a quadratic polynomial. The argument from (\ref{lem:limit theorem for first and second moment:I}) shows that this polynomial is constant on the ball, hence constant on all of $\R^d$. Its value at $a=0$ is zero, so the two sign patterns imply
\[
(a^\top x)(a^\top\xi)+(a^\top y)(a^\top\eta)=0, \quad\text{and}\quad (a^\top x)(a^\top\xi)-(a^\top y)(a^\top\eta)=0
\qquad\text{for every }a\in\R^d.
\]
Adding and subtracting these identities shows that each product vanishes for every $a\in\R^d$. Since $x$ and $y$ are nonzero, this forces $\xi=\eta=0$, a contradiction. This proves~(\ref{cond:ullt-nonlattice}).
Thus all four conditions of Lemma~\ref{core lemma:local limit theorem} hold for $W_{x,y}$ on $K_{\delta_0}$.

We now apply that lemma to $\pi_{n,I,J}(x,y;B)$. Recall that $\nu=n-2d+r=n-\abs{I\cup J}$ is the number of measurements outside the two active sets. The rows $g_1,\ldots,g_\nu$ of $C$ are independent, so the event defining $\pi_{n,I,J}$ involves a sum of $\nu$ independent copies of $W_{x,y}$. By symmetry of the cube, we can write
\[
\pi_{n,I,J}(x,y;B)
=\PP\left\{\sum_{\ell=1}^{\nu}\begin{pmatrix}Z_x(g_\ell)\\Z_y(g_\ell)\end{pmatrix}
\in-b+D(-1,1)^{2d}\right\},
\]
where
\[
b=\begin{pmatrix}
S(x;B_{J\setminus I})\\S(y;B_{I\setminus J})
\end{pmatrix},
\qquad\text{and}\qquad
D=\begin{pmatrix}
Q(x;B_I)&0\\0&Q(y;B_J)
\end{pmatrix}.
\]
When $B_I\in\mathcal K_x$ and $B_J\in\mathcal K_y$, both diagonal blocks of $D$ belong to $\mathcal B_1$. Thus $D$ lies in a fixed compact subset of $GL_{2d}(\R)$, as required for uniformity over the matrices in Lemma~\ref{core lemma:local limit theorem}. We apply that lemma in dimension $2d$, with $\nu$ summands, translation $-b$, and matrix $D$. Multiplying its error by the uniformly bounded factor $2^{2d}\abs{\det D}$ yields
\[
\nu^d\pi_{n,I,J}(x,y;B) =2^{2d}\abs{\det Q(x;B_I)}\abs{\det Q(y;B_J)} \times \varphi_{\scriptscriptstyle\Sigma_{x,y}} \left(-\sqrt\nu\begin{pmatrix}\mu(x)\\\mu(y)\end{pmatrix} -\frac{b}{\sqrt\nu}\right)+o(1),
\]
uniformly for $(x,y)\in K_{\delta_0}$ and matrices $B$ satisfying $B_I\in\mathcal K_x$ and $B_J\in\mathcal K_y$.

To obtain the expression in (\ref{lem:limit theorem for first and second moment:II}), we remove the shift $b/\sqrt\nu$ from the Gaussian argument. Every row of $B$ has norm at most $M_d$, and each sum defining $b$ contains at most $d$ terms. The bound $\norm{Z_x(a)}\le C_d\norm{a}^2$ therefore implies $\norm{b}\le C_d$, uniformly over the index sets and matrices $B$ satisfying $B_I\in\mathcal K_x$ and $B_J\in\mathcal K_y$. Moreover, the uniform covariance lower bound and compactness imply
\[
\sup_{\substack{(x,y)\in K_{\delta_0}\\z\in\R^{2d}}}
\norm{\nabla\varphi_{\scriptscriptstyle\Sigma_{x,y}}(z)}<\infty.
\]
Since $\norm{b}/\sqrt\nu\le C_d/\sqrt\nu$ and the determinant factors are uniformly bounded, removing $b/\sqrt\nu$ from the Gaussian argument changes the displayed expression by $o(1)$, uniformly. Restricting to $x,y\in T_n$ and using $\nu=n-2d+r\to\infty$ proves (\ref{lem:limit theorem for first and second moment:II}).
\end{proof}

\subsection{Estimates in the shrinking tube}\label{app:shrinking-tube-tools}

Recall $Z_x(a)$ from~\eqref{eq:def Z} and its mean $\mu(x)$ and covariance matrices $\Sigma_x$ and $\Sigma_{x,y}$ from~\eqref{eq:def mean cov}. The following lemma gives the uniform estimates near the saddle sphere needed for the local limit approximations. Recall that we set $\kappa=8\tau/[\pi(1+\tau^2)^2]$.

\begin{lemma}[Mean and covariance near the saddle sphere]\label{lem:tube}
For any $n>\max\{d,\tau^{-2}\}$, write
\[
x_{n,z,u}=\frac{z_1}{\sqrt n}\xs+\left(\tau+\frac{z_2}{\sqrt n}\right)u,
\qquad (z,u)\in[-1,1]^2\times\mathbb S^{d-2}.
\]
Then the following hold.
\begin{enumerate}[label=\textnormal{(\Roman*)},ref=\Roman*]
\item\label{lem:tube:I} The mean satisfies
\[
\sup_{\substack{z\in[-1,1]^2\\u\in\mathbb S^{d-2}}}
\norm{\sqrt n\,\mu(x_{n,z,u})+\kappa z_1\xs-\kappa z_2u}\to0
\qquad\text{as }n\to\infty.
\]
\item\label{lem:tube:II} The map $x\mapsto\Sigma_x$ is continuous near $\mathcal S$, and for every $x\in T_n$,
\begin{equation}\label{eq:tube covariance formula}
\Sigma_x=\norm{x}^2I_d+2xx^\top-\mu(x)\mu(x)^\top.
\end{equation}
In particular, $\Sigma_{\tau u}=\tau^2(I_d+2uu^\top)\succ0$ for $u\in\mathbb S^{d-2}$. Consequently,
\[
\inf_{x\in T_n}\lambda_{\min}(\Sigma_x)\ge c_d>0
\]
for all sufficiently large $n$, with the constant and threshold depending only on $d$.
\item\label{lem:tube:III} The map $(x,y)\mapsto\Sigma_{x,y}$ is continuous near $\mathcal S\times\mathcal S$, and $\Sigma_{\tau u,\tau v}\succ0$ whenever $u,v\in\mathbb S^{d-2}$ and $u\ne\pm v$. Consequently, for every compact set $K\subset\{(u,v)\in(\mathbb S^{d-2})^2:u\ne\pm v\}$,
\[
\inf_{\substack{z,w\in[-1,1]^2\\(u,v)\in K}}
\lambda_{\min}\bigl(\Sigma_{x_{n,z,u},x_{n,w,v}}\bigr)\ge c_{d,K}>0
\]
for all sufficiently large $n$, with the constant and threshold depending only on $d$ and $K$.
\end{enumerate}
\end{lemma}
\begin{proof}
By rotational invariance, we may assume without loss of generality that $\xs=e$, where $e$ is the first standard basis vector.
For (\ref{lem:tube:I}), recall that the population function is
\[
\F(x)=\E\abs{(a^\top x)^2-(a^\top e)^2}.
\]
For each fixed $x$ near $\mathcal S$, the random function $x\mapsto\abs{(a^\top x)^2-(a^\top e)^2}$ is differentiable at $x$ almost surely, with gradient $2Z_x(a)$. For $\norm{h}=1$ and $0<\abs{t}\le1$, its difference quotient satisfies
\[
\begin{aligned}
\left|
\frac{\abs{(a^\top(x+th))^2-(a^\top e)^2}
-\abs{(a^\top x)^2-(a^\top e)^2}}{t}
\right|&\le 2\abs{a^\top x}\abs{a^\top h}
+\abs{t}(a^\top h)^2\\
&\le (2\norm{x}+1)\norm{a}^2
\le C_d\norm{a}^2,
\end{aligned}
\]
uniformly for $x$ in a fixed bounded neighborhood of $\mathcal S$. Since $\E[\norm{a}^2]=d<\infty$, dominated convergence implies
\[
\nabla\F(x)=2\E[Z_x(a)]=2\mu(x).
\]
Proposition~\ref{prop:population landscape} therefore implies
\[
\mu(\tau u)=0,
\qquad\text{and}\qquad
D\mu(\tau u)=\kappa(uu^\top-ee^\top).
\]
Since $D\mu$ is continuous near the compact set $\mathcal S$, Taylor expansion implies
\[
\mu(x_{n,z,u})
=\frac{\kappa}{\sqrt n}(-z_1e+z_2u)+o(n^{-1/2}),
\]
uniformly for $z\in[-1,1]^2$ and $u\in\mathbb S^{d-2}$. Multiplying by $\sqrt n$ proves (\ref{lem:tube:I}).

For (\ref{lem:tube:II}), recall that
\[
\Sigma_x=\E[Z_x(a)Z_x(a)^\top]-\mu(x)\mu(x)^\top.
\]
For $x$ near $\mathcal S$, the sign factor in $Z_x(a)$ squares to one almost surely, and a direct computation gives
\[
\E[Z_x(a)Z_x(a)^\top]
=\E[(a^\top x)^2aa^\top]
=\norm{x}^2I_d+2xx^\top.
\]
Combining the two displayed identities proves the formula for $\Sigma_x$ in \eqref{eq:tube covariance formula}. Since $\mu=\nabla\F/2$ is continuous near $\mathcal S$, this formula shows that $\Sigma_x$, and hence $\lambda_{\min}(\Sigma_x)$, is continuous near $\mathcal{S}$. In particular, setting $x=\tau u$ and using $\mu(\tau u)=0$, we obtain \[ \Sigma_{\tau u}=\tau^2(I_d+2uu^\top)\succeq\tau^2I_d.\] Fix $n_0$ with $n_0^{-1/2}<\tau$. Since $T_{n_0}$ is compact and $T_n\subset T_{n_0}$ for $n\ge n_0$, uniform continuity and $\sup_{z,u}\norm{x_{n,z,u}-\tau u}\le\sqrt2n^{-1/2}\to0$ imply that, for all sufficiently large $n$,
\[
\inf_{\substack{z\in[-1,1]^2\\u\in\mathbb S^{d-2}}}\lambda_{\min}\bigl(\Sigma_{x_{n,z,u}}\bigr)\ge c_d:=\tau^2/2.
\]

For (\ref{lem:tube:III}), choose a compact set $U\subset\R^d$ whose interior contains $\mathcal S$ and which contains neither $e$ nor $-e$. Set $W_{x,y}(a)=\begin{pmatrix}Z_x(a)\\Z_y(a)\end{pmatrix}$. For any sequence $(x_j,y_j)\to(x,y)$ in $U\times U$, we have $W_{x_j,y_j}(a)\to W_{x,y}(a)$ almost surely. Moreover, for $x,y\in U$,
\[
\norm{W_{x,y}(a)}^2
=\norm{Z_x(a)}^2+\norm{Z_y(a)}^2
\le(\norm{x}^2+\norm{y}^2)\norm{a}^4
\le2R^2\norm{a}^4,
\]
where $R=\max_{w\in U}\norm{w}$. Since $\E[\norm{a}^4]<\infty$, dominated convergence applied to each entry of the mean vector and second-moment matrix implies $\Sigma_{x_j,y_j}\to\Sigma_{x,y}$.

We next prove that $\Sigma_{\tau u,\tau v}\succ0$ for $u,v\in\mathbb S^{d-2}$ with $u\ne\pm v$. Fix such $u,v$, and set $x=\tau u$ and $y=\tau v$. For $\xi,\eta\in\R^d$, suppose
\[
\operatorname{Var}\bigl(\xi^\top Z_x(a)+\eta^\top Z_y(a)\bigr)=0.
\]
Then $\xi^\top Z_x(a)+\eta^\top Z_y(a)$ is almost surely constant. Since $e,x,y$ are linearly independent, the same polynomial argument used to verify condition~(\ref{cond:ullt-nonlattice}) in the proof of Lemma~\ref{lem:limit theorem for first and second moment}(\ref{lem:limit theorem for first and second moment:II}) implies $\xi=\eta=0$. Thus $\Sigma_{\tau u,\tau v}\succ0$.

For any compact set
\(
K\subset\{(u,v)\in\mathbb S^{d-2}\times\mathbb S^{d-2}:u\ne\pm v\},
\)
the function $(u,v)\mapsto\lambda_{\min}(\Sigma_{\tau u,\tau v})$ is continuous and positive on $K$. Hence
\(
\inf_{(u,v)\in K}\lambda_{\min}(\Sigma_{\tau u,\tau v})>0.
\)
The same uniform-continuity argument as in (\ref{lem:tube:II}), using $x_{n,z,u}\to\tau u$ and $x_{n,w,v}\to\tau v$ uniformly, implies that for all sufficiently large $n$,
\[
\inf_{\substack{z,w\in[-1,1]^2\\(u,v)\in K}}\lambda_{\min}\bigl(\Sigma_{x_{n,z,u},x_{n,w,v}}\bigr)\ge c_{d,K}>0.
\]
\end{proof}

%%% Local Variables:
%%% mode: LaTeX
%%% TeX-master: "../main"
%%% End:

\section{Three ingredients for second moment}\label{app:second moment supporting proofs}
This section contains the supporting lemmas used in the proof of Theorem~\ref{thm:second moment asymptotic appendix}.

In the proofs below, $C_d,c_d>0$ denote constants depending only on $d$, used for upper and lower bounds, respectively. Their values may change from line to line; additional dependencies are indicated by subscripts.

\subsection{Separated region}
We use the separated-region contribution $D_{r,n}^{\mathrm{sep}}$ and the truncated angular integral $\mathcal J_d(\delta_0)$ defined in the proof of Theorem~\ref{thm:second moment asymptotic appendix}, and the constant $\gamma_d$ from Theorem~\ref{thm:first moment asymptotic appendix}.

\begin{lemma}[Separated region]\label{lem:separated region}
Fix $0<\delta_0<1$. As $n\to\infty$, we have $D_{0,n}^{\mathrm{sep}}=\gamma_d^2\mathcal J_d(\delta_0)n^{d-2}(1+o_{d,\delta_0}(1))$ and
\[
D_{r,n}^{\mathrm{sep}}=O_{d,\delta_0}(n^{d-2-r})=o_{d,\delta_0}(n^{d-2}) \qquad\text{for }1\le r\le d-1.
\]
\end{lemma}

\begin{proof}
The quantity $D_{r,n}^{\mathrm{sep}}$ is a sum over ordered pairs of active sets with overlap $r$. We first estimate the contribution of one such pair, and then sum over the pairs.

Fix $I\ne J$ with $\abs I=\abs J=d$ and $\abs{I\cap J}=r$. Recall from the proof of Theorem~\ref{thm:second moment asymptotic appendix} that
\[
\begin{aligned}
I_{r,n}^{\mathrm{sep}}
&=\int_{\Omega^{\mathrm{sep}}}p_{X_I(x),X_J(y)}(0,0)\\
&\quad\times\E\left[\begin{aligned}
&\abs{\det A_I}\abs{\det A_J}1_{\mathcal K_x}(A_I)1_{\mathcal K_y}(A_J)\\
&\quad\times\pi_{n,I,J}(x,y;A_{I\cup J})
\end{aligned}\;\middle|\;X_I(x)=X_J(y)=0\right]\dd x\,\dd y.
\end{aligned}
\]
We first state the bounds used in both cases. For all sufficiently large $n$, Lemmas~\ref{lem:density}(\ref{lem:density:I}) and~\ref{lem:tube}(\ref{lem:tube:III}) imply, uniformly over $\Omega^{\mathrm{sep}}$,
\begin{equation}\label{eq:sep covariance}
p_{X_I(x),X_J(y)}(0,0)\le C_d\delta(x,y)^{-r}, \qquad\text{and}\qquad \lambda_{\min}(\Sigma_{x,y})\ge c_{d,\delta_0}>0.
\end{equation}
The covariance bound follows from the continuity and positive definiteness of $(x,y)\mapsto\Sigma_{x,y}$ on the compact set of limiting pairs $(\tau u_x,\tau u_y)$ with $\delta(x,y)\ge\delta_0$. In particular,
\begin{equation}\label{eq:sep gaussian density bound}
\sup_{v\in\R^{2d}}\varphi_{\scriptscriptstyle\Sigma_{x,y}}(v)
=\frac{1}{(2\pi)^d\sqrt{\det\Sigma_{x,y}}}\le C_{d,\delta_0}.
\end{equation}
Also, the row norm bounds defining $\mathcal K_x$ and $\mathcal K_y$, together with boundedness of $T_n$, imply
\begin{equation}\label{eq:second moment determinant bound}
\abs{\det A_I}\abs{\det A_J}\le M_d^{2d},\qquad\text{and}\qquad
\abs{\det Q(x;A_I)}\abs{\det Q(y;A_J)}\le C_d
\end{equation}
whenever $x,y\in T_n$, $A_I\in\mathcal K_x$, and $A_J\in\mathcal K_y$.

Let $\nu=n-2d+r=n-\abs{I\cup J}$, the number of rows outside the two active sets. For $A_I\in\mathcal K_x$ and $A_J\in\mathcal K_y$, Lemma~\ref{lem:limit theorem for first and second moment}(\ref{lem:limit theorem for first and second moment:II}) implies
\begin{equation}\label{eq:separated local limit expansion}
\begin{aligned}
\pi_{n,I,J}(x,y;A_{I\cup J})
&=\frac{2^{2d}}{\nu^d}\abs{\det Q(x;A_I)}\abs{\det Q(y;A_J)}
\varphi_{\scriptscriptstyle\Sigma_{x,y}}\left(-\sqrt\nu\begin{pmatrix}\mu(x)\\\mu(y)\end{pmatrix}\right)\\
&\quad+\nu^{-d}r_n(x,y;A_{I\cup J}),
\end{aligned}
\end{equation}
where, for some deterministic sequence $\eta_n\to0$,
\begin{equation}\label{eq:separated local limit remainder}
\sup_{(x,y)\in\Omega^{\mathrm{sep}}}\sup_{\substack{B\in\R^{(2d-r)\times d}\\B_I\in\mathcal K_x,\ B_J\in\mathcal K_y}}
\abs{r_n(x,y;B)}\le\eta_n.
\end{equation}
We now use \eqref{eq:separated local limit expansion} to bound the terms with $r\ge1$ and evaluate the leading term with $r=0$.

\paragraph{Case 1: $1\le r\le d-1$.}
Combining~\eqref{eq:sep gaussian density bound}, \eqref{eq:second moment determinant bound}, \eqref{eq:separated local limit expansion}, and~\eqref{eq:separated local limit remainder}, we obtain
\begin{equation}\label{eq:sep pair probability bound}
\sup_{\substack{B:B_I\in\mathcal K_x\\B_J\in\mathcal K_y}}
\pi_{n,I,J}(x,y;B)\le\nu^{-d}(C_{d,\delta_0}+\eta_n)\le C_{d,\delta_0}\nu^{-d}.
\end{equation}
Substituting \eqref{eq:sep pair probability bound} into the definition of $I_{r,n}^{\mathrm{sep}}$ and using \eqref{eq:sep covariance}, we obtain
\begin{equation}\label{eq:sep overlapping pair integral bound}
I_{r,n}^{\mathrm{sep}} \le C_{d,\delta_0}\nu^{-d}\int_{\Omega^{\mathrm{sep}}}p_{X_I(x),X_J(y)}(0,0)\dd x\,\dd y \le C_{d,\delta_0}\nu^{-d}\delta_0^{-r}\operatorname{Vol}(T_n)^2 \le C_{d,\delta_0}n^{-d-2},
\end{equation}
since $\nu/n\to1$ and $\operatorname{Vol}(T_n)=O_d(n^{-1})$. The bound \eqref{eq:sep overlapping pair integral bound} is uniform over all pairs with overlap $r$. There are $\binom nd\binom dr\binom{n-d}{d-r}=O_d(n^{2d-r})$ such ordered pairs, so
\[
D_{r,n}^{\mathrm{sep}}\le C_{d,\delta_0}n^{2d-r}n^{-d-2}
=C_{d,\delta_0}n^{d-2-r}=o_{d,\delta_0}(n^{d-2}).
\]

\paragraph{Case 2: $r=0$.}
Now $I\cap J=\varnothing$ and $\nu=n-2d$. Write $I=\{i_1<\cdots<i_d\}$. For fixed $x\in T_n$ and each $k=1,\ldots,d$, the $k$-th coordinate of $X_I(x)$ is $a_{i_k}^\top(x-s_k\xs)$, which is Gaussian with mean zero and variance $\norm{x-s_k\xs}^2$. These coordinates are independent, so
\begin{equation}\label{eq:separated marginal density}
p_{X_I(x)}(0)=(2\pi)^{-d/2}\prod_{k=1}^d\frac{1}{\norm{x-s_k\xs}}.
\end{equation}
For sufficiently large $n$, these variances are bounded away from zero on $T_n$, since $\norm{x-s_k\xs}\ge\tau-n^{-1/2}\ge\tau/2$. The formula \eqref{eq:separated marginal density} also holds for $X_J(y)$. Since the two active sets are disjoint, $X_I(x)$ and $X_J(y)$ are independent; hence $T_n\times T_n\subset\mathcal R_{I,J}$ and
\begin{equation}\label{eq:sep joint density factorization}
p_{X_I(x),X_J(y)}(0,0)=p_{X_I(x)}(0)p_{X_J(y)}(0).
\end{equation}
Define
\[
a(x)=\E\bigl[\abs{\det A_I}\abs{\det Q(x;A_I)}1_{\mathcal K_x}(A_I)\bigm|X_I(x)=0\bigr],
\]
as in Lemma~\ref{lem:regularity constants}. This function is independent of $I$, since every ordered active set uses the same sign pattern. Moreover, $A_I$ and $A_J$ remain independent under the conditioning $X_I(x)=X_J(y)=0$, because each constraint concerns only its own block. Consequently,
\begin{equation}\label{eq:sep conditional factorization}
\E\left[\begin{aligned}
&\abs{\det A_I}\abs{\det Q(x;A_I)}1_{\mathcal K_x}(A_I)\\
&\quad\times\abs{\det A_J}\abs{\det Q(y;A_J)}1_{\mathcal K_y}(A_J)
\end{aligned}\;\middle|\;X_I(x)=X_J(y)=0\right]=a(x)a(y).
\end{equation}

We now substitute the expansion~\eqref{eq:separated local limit expansion} for $\pi_{n,I,J}(x,y;A_{I\cup J})$ in the integral defining $I_{0,n}^{\mathrm{sep}}$. For fixed $x,y$, the Gaussian density in that expansion is deterministic, so it can be taken outside the conditional expectation. The remaining expectation in the leading term factorizes as in \eqref{eq:sep conditional factorization}. Together with \eqref{eq:sep joint density factorization}, this yields $I_{0,n}^{\mathrm{sep}}=I_1+I_2$, where
\[
I_1=\frac{2^{2d}}{\nu^d}\int_{\Omega^{\mathrm{sep}}}p_{X_I(x)}(0)p_{X_J(y)}(0)a(x)a(y)
\varphi_{\scriptscriptstyle\Sigma_{x,y}}\left(-\sqrt\nu\begin{pmatrix}\mu(x)\\\mu(y)\end{pmatrix}\right)\dd x\,\dd y,
\]
and $I_2$ contains the remainder $r_n$. By~\eqref{eq:separated local limit remainder}, \eqref{eq:separated marginal density}, and~\eqref{eq:second moment determinant bound},
\begin{equation}\label{eq:sep remainder integral bound}
\abs{I_2}\le C_d\eta_n\nu^{-d}\operatorname{Vol}(T_n)^2=o_{d,\delta_0}(n^{-d-2}).
\end{equation}
It remains to evaluate $I_1$.

Define the parametrization
\[
x_n(z,u)=\frac{z_1}{\sqrt n}\xs+\left(\tau+\frac{z_2}{\sqrt n}\right)u
\qquad\text{for }z\in[-1,1]^2\text{ and }u\in\mathbb S^{d-2}.
\]
We use this map for both integration variables, writing $x=x_n(z,u)$ and $y=x_n(w,v)$, where $z,w\in[-1,1]^2$ and $u,v\in\mathbb S^{d-2}$.
For each of $x$ and $y$, the two normal coordinates each contribute a factor $n^{-1/2}$, and polar coordinates in $(\xs)^\perp$ contribute the radial factor $(\tau+z_2/\sqrt n)^{d-2}$. Consequently,
\[
\dd x\,\dd y=n^{-2}\left(\tau+\frac{z_2}{\sqrt n}\right)^{d-2}
\left(\tau+\frac{w_2}{\sqrt n}\right)^{d-2}
\dd z\,\dd w\,\dd\sigma(u)\,\dd\sigma(v).
\]
The restriction $\delta(x,y)\ge\delta_0$ becomes $\min\{\norm{u-v},\norm{u+v}\}\ge\delta_0$, so the transformed integration domain is compact and independent of $n$.

We now let $n\to\infty$. Since $\norm{\tau u-s_k\xs}^2=1+\tau^2$, equation~\eqref{eq:separated marginal density} and Lemma~\ref{lem:regularity constants} imply
\begin{equation}\label{eq:sep marginal limits}
\begin{aligned}
p_{X_I(x_n(z,u))}(0)&\to[2\pi(1+\tau^2)]^{-d/2},\qquad
a(x_n(z,u))\to\vartheta_d,\\
p_{X_J(x_n(w,v))}(0)&\to[2\pi(1+\tau^2)]^{-d/2},\qquad\text{and}\qquad
a(x_n(w,v))\to\vartheta_d.
\end{aligned}
\end{equation}
Also, Lemma~\ref{lem:tube}(\ref{lem:tube:I}) and~(\ref{lem:tube:III}), together with $\nu/(n-d)\to1$, implies
\begin{equation}\label{eq:sep mean covariance limits}
\begin{aligned}
\sqrt\nu\,\mu(x_n(z,u))&\to-\kappa z_1\xs+\kappa z_2u,\\
\sqrt\nu\,\mu(x_n(w,v))&\to-\kappa w_1\xs+\kappa w_2v,\quad\text{and}\\
\norm{\Sigma_{x_n(z,u),x_n(w,v)}-\Sigma_{\tau u,\tau v}}_{\mathrm{op}}&\to0.
\end{aligned}
\end{equation}
The limits in \eqref{eq:sep marginal limits} and \eqref{eq:sep mean covariance limits} are uniform over the transformed integration domain. By~\eqref{eq:sep covariance}, the covariance matrices stay uniformly positive definite, so the Gaussian densities also converge uniformly. We may therefore pass to the limit in the integral. Using $\varphi_\Sigma(-v)=\varphi_\Sigma(v)$, we obtain
\begin{equation}\label{eq:sep leading integral limit}
n^2\nu^d I_1\longrightarrow
\frac{2^{2d}\tau^{2d-4}\vartheta_d^2}{[2\pi(1+\tau^2)]^d}\mathcal J_d(\delta_0)
=(d!)^2\gamma_d^2\mathcal J_d(\delta_0).
\end{equation}
Combining \eqref{eq:sep leading integral limit}, \eqref{eq:sep remainder integral bound}, and $\nu/n\to1$, we obtain
\begin{equation}\label{eq:sep disjoint pair asymptotic}
I_{0,n}^{\mathrm{sep}}=(d!)^2\gamma_d^2\mathcal J_d(\delta_0)n^{-d-2}(1+o_{d,\delta_0}(1)).
\end{equation}
Relabeling the Gaussian rows shows that all disjoint pairs have the same integral. Multiplying \eqref{eq:sep disjoint pair asymptotic} by their number,
$\binom nd\binom{n-d}{d}=n^{2d}(1+O_d(n^{-1}))/(d!)^2$, proves the asserted formula for $D_{0,n}^{\mathrm{sep}}$.
\end{proof}

%%%%%%%%%%%%%%%%%%%%%%%%%%%%%%%%%Tangential region
\subsection{Tangential region}
We continue to use the notation $\Omega^{\mathrm{tan}}$ and $D_{r,n}^{\mathrm{tan}}$ from the proof of Theorem~\ref{thm:second moment asymptotic appendix}.

\begin{lemma}[Tangential region]\label{lem:tangential region}
Fix $d\ge4$ and $0<\delta_0<1$. For all sufficiently large $n$, we have \[ D_{0,n}^{\mathrm{tan}}\le C_dn^{d-2}\delta_0^{d/2-1}+C_dn^{d/2-1}. \]
For $1\le r\le d-1$, \[ D_{r,n}^{\mathrm{tan}}=o_d(n^{d-2}). \]
The constants $C_d$ are independent of $\delta_0$.
\end{lemma}
\begin{proof}
Since there are a total of $\binom nd\binom dr\binom{n-d}{d-r}=O_d(n^{2d-r})$ ordered pairs $(I,J)$ with overlap $r$, it suffices to show that, for any fixed $I\ne J$ with $\abs{I}=\abs{J}=d$ and $\abs{I\cap J}=r$, the integral
\begin{equation}\label{eq:tangential pair integral}
\begin{aligned}
I_{r,n}^{\mathrm{tan}}
&:=\int_{\Omega^{\mathrm{tan}}} p_{X_I(x),X_J(y)}(0,0)\\
&\quad\times\E\left[\begin{aligned}
&\abs{\det A_I}\abs{\det A_J}1_{\mathcal K_x}(A_I)1_{\mathcal K_y}(A_J)\\
&\quad\times\pi_{n,I,J}(x,y;A_{I\cup J})
\end{aligned}\;\middle|\;X_I(x)=X_J(y)=0\right]\dd x\,\dd y,
\end{aligned}
\end{equation}
satisfies
\begin{equation}\label{eq:tan disjoint pair bound}
I_{0,n}^{\mathrm{tan}}\le C_dn^{-3d/2-1}+C_dn^{-d-2}\delta_0^{d/2-1},
\end{equation}
and
\begin{equation}\label{eq:tan overlapping pair bound}
I_{r,n}^{\mathrm{tan}}=o_d(n^{-d+r-2})\qquad\text{for }1\le r\le d-1.
\end{equation}

Since the set $\{\delta(x,y)=0\}$ has Lebesgue measure zero, we first estimate the integral over
\[
\Omega_+^{\mathrm{tan}}=\Omega^{\mathrm{tan}}\cap\{(x,y):u_x^\top u_y>0,\ \delta(x,y)>0\}.
\]
For sufficiently large $n$, we have $n^{-1}<\delta_0$. We split at $\delta(x,y)=n^{-1}$, where the prefactor $n^{-d}\delta(x,y)^{-d/2}$ in the joint probability estimate equals the marginal scale $n^{-d/2}$. For larger separations, we use the joint estimate and its Gaussian decay in the normal displacement; for smaller separations, we use the marginal estimate and the small volume imposed by $\rho(x,y)\le\delta(x,y)$. Define
\[
\Omega_+^{(1)}= \Omega_+^{\mathrm{tan}} \cap \{(x,y):n^{-1}\le \delta(x,y)\le \delta_0\} \quad\text{and}\quad \Omega_+^{(2)}= \Omega_+^{\mathrm{tan}} \cap \{(x,y):0< \delta(x,y)< n^{-1}\},
 \]
and denote by $I_1$ and $I_2$  the integrals over $\Omega_+^{(1)}$ and $\Omega_+^{(2)}$, respectively.

\paragraph{Separations \texorpdfstring{$\delta(x,y)\ge n^{-1}$}{delta(x,y) >= 1/n}.}
We claim that, uniformly on $\Omega_+^{\mathrm{tan}}$,
\begin{equation}\label{eq:tan pair probability claim}
\sup_{\substack{B:B_I\in\mathcal K_x\\B_J\in\mathcal K_y}}\pi_{n,I,J}(x,y;B)\le C_dn^{-d}\delta(x,y)^{-d/2}\exp\left[-c_dn\left(\alpha_x^2+\beta_x^2+\frac{\rho(x,y)^2}{\delta(x,y)}\right)\right].
\end{equation}
We defer the proof of \eqref{eq:tan pair probability claim} to Lemma~\ref{lem:tangential region estimate}. By \eqref{eq:tan pair probability claim}, Lemma~\ref{lem:density}(\ref{lem:density:I}), and \eqref{eq:second moment determinant bound}, we obtain \begin{equation}\label{eq:tan density} I_1\le C_dn^{-d}\int_{\substack{(x,y)\in\Omega_+^{\mathrm{tan}}\\n^{-1}\le\delta(x,y)<\delta_0}}\delta(x,y)^{-r-d/2}\exp\left[-c_dn\left(\alpha_x^2+\beta_x^2+\frac{\rho(x,y)^2}{\delta(x,y)}\right)\right]\dd x\,\dd y. \end{equation}
We use tube coordinates, in which \[ \dd x\,\dd y=(\tau+\beta_x)^{d-2}(\tau+\beta_y)^{d-2}\dd\alpha_x\,\dd\beta_x\,\dd\alpha_y\,\dd\beta_y\,\dd\sigma(u_x)\,\dd\sigma(u_y). \]
Since $\delta(x,y)$ depends only on $u_x$ and $u_y$, \eqref{eq:tan density} implies
\ificlr
\[
\begin{aligned}
I_1&\le C_dn^{-d}\\
&\quad\times\int_{\Omega_u}\delta^{-r-d/2}
\underbrace{\int_{\Omega_{\alpha,\beta}(\delta)}
e^{-c_dn(\alpha_x^2+\beta_x^2+\rho^2/\delta)}
\prod_{q\in\{x,y\}}\bigl[(\tau+\beta_q)^{d-2}\,\dd\alpha_q\,\dd\beta_q\bigr]}_{I_3}
\dd\sigma(u_x)\,\dd\sigma(u_y),
\end{aligned}
\]
\else
\[
I_1\le C_dn^{-d}\int_{\Omega_u}\delta^{-r-d/2}
\underbrace{\int_{\Omega_{\alpha,\beta}(\delta)}
e^{-c_dn(\alpha_x^2+\beta_x^2+\rho^2/\delta)}
\prod_{q\in\{x,y\}}\bigl[(\tau+\beta_q)^{d-2}\,\dd\alpha_q\,\dd\beta_q\bigr]}_{I_3}
\dd\sigma(u_x)\,\dd\sigma(u_y),
\]
\fi
where $\delta=\delta(x,y)$, $\rho=\rho(x,y)$, and, with $u_x,u_y\in\mathbb S^{d-2}$,
\[
\Omega_u=\{(u_x,u_y):u_x^\top u_y>0,\ n^{-1}\le\delta<\delta_0\},
\quad\text{and}\quad
\Omega_{\alpha,\beta}(\delta)=[-n^{-1/2},n^{-1/2}]^4\cap\{\rho\le\delta\}.
\]
The change of variables $\zeta=\alpha_x-\alpha_y$ and $\eta=\beta_x-\beta_y$, under which $\rho(x,y)^2=\zeta^2+\eta^2$, has Jacobian of absolute value one. The tube Jacobian factors are bounded, so after enlarging $\Omega_{\alpha,\beta}(\delta)$, we obtain \begin{equation}\label{eq:tan normal coordinate integral}
I_3\le C_d\int_{\R^2}e^{-c_dn(\alpha_x^2+\beta_x^2)}\dd\alpha_x\,\dd\beta_x\int_{\R^2}e^{-c_dn(\zeta^2+\eta^2)/\delta}\dd\zeta\,\dd\eta \le C_dn^{-1}\frac{\delta}{n}=C_d\delta n^{-2}.
\end{equation}
Substituting \eqref{eq:tan normal coordinate integral} and then using Lemma~\ref{lem:delta integral}, we obtain \begin{equation}\label{eq:tan final integral}
I_1\le C_dn^{-d-2}\int_{\substack{u_x,u_y\in\mathbb S^{d-2}\\n^{-1}\le\delta(u_x,u_y)<\delta_0}}\delta(u_x,u_y)^{1-r-d/2}\dd\sigma(u_x)\,\dd\sigma(u_y) \le C_dn^{-d-2}\int_{1/n}^{\delta_0}t^{d/2-2-r}\dd t.
\end{equation}
Evaluating the integral in \eqref{eq:tan final integral} gives \begin{equation}\label{eq:tan first integral bound} I_1\le\begin{cases} C_dn^{-d-2}\delta_0^{d/2-1-r},&r<d/2-1,\\ C_dn^{-d-2}\log n,&r=d/2-1,\\ C_dn^{-3d/2-1+r},&r>d/2-1.\end{cases} \end{equation}

\paragraph{Separations \texorpdfstring{$0<\delta(x,y)<n^{-1}$}{0 < delta(x,y) < 1/n}.}
Condition on $A_{I\cup J}=B$, and let $C_\nu$ contain the $\nu=n-2d+r$ remaining independent Gaussian rows. The constraint at $x$ becomes \[ S(x;C_\nu)\in-Q(x;B_I)(-1,1)^d-S(x;B_{J\setminus I}). \]
When $B_I\in\mathcal K_x$ and $B_J\in\mathcal K_y$, the target set on the right-hand side lies in $B(0,R_d)$ for a fixed $R_d$. Applying the local limit argument for $S(x;C_\nu)$ from the proof of Lemma~\ref{lem:limit theorem for first and second moment}(\ref{lem:limit theorem for first and second moment:I}) to a fixed cube containing this ball and using the uniform covariance lower bound, we obtain \begin{equation}\label{eq:tan marginal probability bound} \sup_{\substack{B:B_I\in\mathcal K_x\\B_J\in\mathcal K_y}}\pi_{n,I,J}(x,y;B)\le C_d\nu^{-d/2}\le C_dn^{-d/2}. \end{equation}
Using Lemma~\ref{lem:density}(\ref{lem:density:I}) and the bounded tube Jacobian factors, we obtain
\[
I_2\le C_dn^{-d/2}\int_{\substack{u_x,u_y\in\mathbb S^{d-2}\\0<\delta(u_x,u_y)<1/n}}\delta^{-r}
\underbrace{\int_{\Omega_{\alpha,\beta}(\delta)}
\dd\alpha_x\,\dd\beta_x\,\dd\alpha_y\,\dd\beta_y}_{I_4}
\dd\sigma(u_x)\,\dd\sigma(u_y).
\]
The change of variables $\zeta=\alpha_x-\alpha_y$ and $\eta=\beta_x-\beta_y$ has Jacobian of absolute value one and transforms $\rho\le\delta$ into $\zeta^2+\eta^2\le\delta^2$. After enlarging the image of $\Omega_{\alpha,\beta}(\delta)$ under this change of variables, we obtain \[ I_4\le\int_{[-n^{-1/2},n^{-1/2}]^2}\dd\alpha_x\,\dd\beta_x\int_{\zeta^2+\eta^2\le\delta^2}\dd\zeta\,\dd\eta=\frac{4\pi\delta^2}{n}. \]
Consequently, Lemma~\ref{lem:delta integral} yields \begin{equation}\label{eq:tan second integral bound}\begin{aligned} I_2&\le C_dn^{-d/2-1}\int_{\substack{u_x,u_y\in\mathbb S^{d-2}\\0<\delta(u_x,u_y)<1/n}}\delta(u_x,u_y)^{2-r}\dd\sigma(u_x)\,\dd\sigma(u_y)\\ &\le C_dn^{-d/2-1}\int_0^{1/n}t^{d-1-r}\dd t\le C_dn^{-3d/2-1+r}, \end{aligned}\end{equation}
where the last integral is finite because $r\le d-1$.

\paragraph{The antipodal contribution.}
For pairs with $u_x^\top u_y<0$, replace $y$ by $y'=-y$ and reverse all active signs assigned to $J$. Each active equation for $J$ is then multiplied by $-1$. Moreover,
\[
Z_{y'}(a)=-Z_y(a),\qquad Q(y';A_J)=-Q(y;A_J),\qquad\text{and}\qquad \mathcal K_{y'}=\mathcal K_y.
\]
Thus the joint density at zero, the conditional law of the active rows, and the determinant factors are unchanged, while the second score constraint is multiplied by $-1$. The antipodal contribution is therefore the corresponding integral in the positive chart with reversed signs for $J$. Indeed, $y'$ has coordinates $(-\alpha_y,\beta_y,-u_y)$ and $\dd y'=\dd y$. The two terms defining $\delta(x,y)$ are interchanged, and, since $u_x^\top u_y<0$, both $\rho(x,y)^2$ and $\rho(x,y')^2$ equal $(\alpha_x+\alpha_y)^2+(\beta_x-\beta_y)^2$.
The estimates above are uniform under this sign reversal. The probability bounds for $\pi_{n,I,J}$ use only the boundedness of the target sets when $B_I\in\mathcal K_x$ and $B_J\in\mathcal K_y$. The density bounds in Lemma~\ref{lem:density}(\ref{lem:density:I})--(\ref{lem:density:II}) follow from the same rowwise calculation for the reversed signs. For the angular bound in Lemma~\ref{lem:density}(\ref{lem:density:III}), the signs on $I$ remain fixed, so the numbers of shared rows with matching signs satisfy $r_+\le d-2$ and $r_-\le2\le d-2$, even after reversing the signs on $J$.

Combining \eqref{eq:tan first integral bound} and \eqref{eq:tan second integral bound} with the reflected antipodal contribution proves \eqref{eq:tan disjoint pair bound} for $r=0$ and \eqref{eq:tan overlapping pair bound} for $1\le r\le d-1$, uniformly over $I,J$. This completes the proof.
\end{proof}

\begin{lemma}[Distance identity]\label{lem:delta rho distance identity}
For all sufficiently large $n$ and every $x,y\in T_n$ with $u_x^\top u_y>0$,
\[
\norm{y-x}^2
=\rho(x,y)^2+(\tau+\beta_x)(\tau+\beta_y)\delta(x,y)^2.
\]
\end{lemma}
\begin{proof}
Since $\xs\perp u_x,u_y$ and $u_x^\top u_y>0$, we have $\delta(x,y)=\norm{u_x-u_y}$ and $\rho(x,y)^2=(\alpha_y-\alpha_x)^2+(\beta_y-\beta_x)^2$. Therefore
\begin{align*}
\norm{y-x}^2
&=(\alpha_y-\alpha_x)^2+
\norm{(\tau+\beta_y)u_y-(\tau+\beta_x)u_x}^2\\
&=(\alpha_y-\alpha_x)^2+(\beta_y-\beta_x)^2
+(\tau+\beta_x)(\tau+\beta_y)\norm{u_y-u_x}^2\\
&=\rho(x,y)^2+(\tau+\beta_x)(\tau+\beta_y)\delta(x,y)^2.
\end{align*}
\end{proof}

\begin{lemma}[Tangential-region probability estimate]\label{lem:tangential region estimate}
Fix $d\ge4$ and active sets $I\ne J$ with $\abs{I}=\abs{J}=d$. There are constants $C_d,c_d>0$ such that, for all sufficiently large $n$ and every $x,y\in T_n$ satisfying $u_x^\top u_y>0$, $0<\delta(x,y)<1$, and $\rho(x,y)\le\delta(x,y)$, \[ \sup_{\substack{B:B_I\in\mathcal K_x\\B_J\in\mathcal K_y}}\pi_{n,I,J}(x,y;B)\le C_dn^{-d}\delta(x,y)^{-d/2}\exp\left[-c_dn\left(\alpha_x^2+\beta_x^2+\frac{\rho(x,y)^2}{\delta(x,y)}\right)\right]. \]
\end{lemma}
\begin{proof}
By rotational invariance of the standard Gaussian measurements, we may assume without loss of generality that $\xs=e$, where $e$ is the first standard basis vector.
Let $r=\abs{I\cap J}$ and $\nu=n-2d+r$.
Put $\delta=\delta(x,y)$, $h=y-x$, and $t=\norm{h}$. Lemma~\ref{lem:delta rho distance identity} and $\rho(x,y)\le\delta$ imply
\[
t^2=\rho(x,y)^2+(\tau+\beta_x)(\tau+\beta_y)\delta^2,
\qquad\text{and}\qquad c_d\delta\le t\le C_d\delta.
\]
Choose $0<\delta_*<1$, depending only on $d$, sufficiently small that $C_d\delta_*\le t_1$, where $C_d$ is the constant in $t\le C_d\delta$ and $t_1$ is as in Lemma~\ref{lem:tangential concentration}.

First consider $0<\delta\le\delta_*$. Then $t\le t_1$, and the tangential-region assumptions of Lemma~\ref{lem:tangential concentration} hold, so it applies with constants depending only on $d$.
Conditional on the active rows $A_{I\cup J}=B$, relabel the remaining independent standard Gaussian rows as $g_1,\ldots,g_\nu$. The target sets for the sums at $x$ and $y$ are
\[
\mathcal B_x=-Q(x;B_I)(-1,1)^d-S(x;B_{J\setminus I}),
\qquad\text{and}\qquad
\mathcal B_y=-Q(y;B_J)(-1,1)^d-S(y;B_{I\setminus J}).
\]
When $B_I\in\mathcal K_x$ and $B_J\in\mathcal K_y$, both targets lie in $B(0,R_d)$ for a fixed $R_d$. To combine the exponential estimate in Lemma~\ref{lem:tangential one row}(\ref{lem:tangential one row:II}) with the concentration estimate in Lemma~\ref{lem:tangential concentration}(\ref{lem:tangential concentration:II}), we use two independent groups of rows. Take the first $\nu_1=\lfloor\nu/2\rfloor$ rows for $C_1$ and the rest for $C_2$:
\[
C_1=(g_1,\ldots,g_{\nu_1})^\top,\qquad\text{and}\qquad
C_2=(g_{\nu_1+1},\ldots,g_\nu)^\top.
\]
Write $\nu_2=\nu-\nu_1$; both groups have $n/2+O_d(1)$ rows. Define
\[
U_j=S(x;C_j),\qquad\text{and}\qquad V_j=S(y;C_j)-S(x;C_j)\qquad\text{for }j=1,2.
\]
Thus $\pi_{n,I,J}(x,y;B)=\PP(F)$, where
\[
F=\{U_1+U_2\in\mathcal B_x,\quad U_1+U_2+V_1+V_2\in\mathcal B_y\}.
\]

To apply Lemma~\ref{lem:tangential one row}(\ref{lem:tangential one row:II}) to each group, let
\[
H_1=\sum_{\ell=1}^{\nu_1}L(g_\ell),
\qquad\text{and}\qquad
H_2=\sum_{\ell=\nu_1+1}^{\nu}L(g_\ell),
\]
where $L$ is defined in that lemma. Thus
\[
H_j=\mu(x)^\top U_j+
\frac{[\mu(y)-\mu(x)]^\top V_j}{t}
\qquad\text{for }j=1,2.
\]
Lemma~\ref{lem:tangential one row}(\ref{lem:tangential one row:I}) implies
\[
\norm{\mu(y)-\mu(x)}
=\norm{\E[Z_y(a)-Z_x(a)]}\le C_dt.
\]
Also, $\norm{\mu(x)}\le C_d$ on $T_n$. On the event $F$,
the bounds $\mathcal B_x,\mathcal B_y\subset B(0,R_d)$ imply $\norm{U_1+U_2}\le R_d$ and $\norm{V_1+V_2}\le2R_d$. Hence
\[
H_1+H_2\le R_d\norm{\mu(x)}+
2R_d\frac{\norm{\mu(y)-\mu(x)}}{t}\le R'_d
\]
for a fixed $R'_d$. Thus at least one of the events $E_j=\{H_j\le R'_d/2\}$ occurs. Lemma~\ref{lem:tangential one row}(\ref{lem:tangential one row:II}) implies
\[
\PP(E_j)
\le e^{\lambda R'_d/2-c_d\nu_jq}
\le C_de^{-c_dnq},
\]
where $q=\norm{\mu(x)}^2+\norm{\mu(y)-\mu(x)}^2/t$.

To apply Lemma~\ref{lem:tangential concentration}(\ref{lem:tangential concentration:II}) to the other group, define
\[
\mathcal T=\{(u,v):u\in\mathcal B_x,\ u+v\in\mathcal B_y\}
\subset B(0,\sqrt5R_d).
\]
For $j\in\{1,2\}$, the event $E_j$ is determined by $C_j$. Since $C_{3-j}$ is independent of $C_j$, Lemma~\ref{lem:tangential concentration}(\ref{lem:tangential concentration:II}) applies uniformly to the translated set and implies
\begin{align*}
\PP(F\cap E_j)
&=\E\bigl[1_{E_j}\PP(F\mid C_j)\bigr]\\
&=\E\bigl[1_{E_j}\PP((U_{3-j},V_{3-j})\in\mathcal T-(U_j,V_j)\mid C_j)\bigr]\\
&\le C_d\nu_{3-j}^{-d}t^{-d/2}\PP(E_j)
\le C_dn^{-d}t^{-d/2}\PP(E_j).
\end{align*}
Summing over $j=1,2$ proves
\begin{equation}\label{eq:tan mean penalty}
\pi_{n,I,J}(x,y;B)
\le C_dn^{-d}t^{-d/2}
\exp\left[-c_dn\left(\norm{\mu(x)}^2+
\frac{\norm{\mu(y)-\mu(x)}^2}{t}\right)\right].
\end{equation}

We express the exponent in tube coordinates. Put $\widetilde y=\alpha_y e+(\tau+\beta_y)u_x$, so $\norm{\widetilde y-x}=\rho(x,y)$. By Proposition~\ref{prop:population landscape}, $\mu(\tau u_x)=0$ and
\[
D\mu(\tau u_x)=\kappa(u_xu_x^\top-ee^\top),
\qquad\text{and}\qquad\kappa=\frac{8\tau}{\pi(1+\tau^2)^2}>0.
\]
The map $u_xu_x^\top-ee^\top$ is an isometry on $\operatorname{span}\{e,u_x\}$. Continuity of $D\mu$ near the compact saddle sphere therefore implies, for sufficiently large $n$,
\begin{equation}\label{eq:tan mean coordinate lower bounds}
\norm{\mu(x)}\ge c_d\norm{z_x},
\qquad\text{and}\qquad
\norm{\mu(\widetilde y)-\mu(x)}\ge c_d\rho(x,y),
\end{equation}
where $z_x=(\alpha_x,\beta_x)$ and $z_y=(\alpha_y,\beta_y)$. For every orthogonal map $R$ fixing $e$, rotational invariance implies $\mu(Rz)=R\mu(z)$. For any unit vector $v\perp e,u_x$, the reflection across $v^\perp$ fixes $\widetilde y$, so $v^\top\mu(\widetilde y)=0$. Thus $\mu(\widetilde y)\in\operatorname{span}\{e,u_x\}$. Choosing $R$ with $Re=e$ and $Ru_x=u_y$ then yields
\[
\mu(y)-\mu(\widetilde y)
=(u_x^\top\mu(\widetilde y))(u_y-u_x).
\]
Since $\mu(\tau u_x)=0$ and $D\mu$ is uniformly bounded near $\mathcal S$, we have $\norm{\mu(\widetilde y)}\le C_d\norm{z_y}$. Hence
\begin{equation}\label{eq:tan angular mean difference}
\norm{\mu(y)-\mu(\widetilde y)}
\le\norm{\mu(\widetilde y)}\delta
\le C_d\norm{z_y}\delta.
\end{equation}
Combining \eqref{eq:tan mean coordinate lower bounds} and \eqref{eq:tan angular mean difference} implies
\[
\norm{\mu(y)-\mu(x)}^2
\ge c_d\rho(x,y)^2-C_d\norm{z_y}^2\delta^2.
\]
Since $\norm{z_y}^2\le2\norm{z_x}^2+2\rho(x,y)^2$ and $t\asymp\delta$, after reducing $\delta_*$ if necessary, we obtain
\[
q\ge c_d\left(\alpha_x^2+\beta_x^2+
\frac{\rho(x,y)^2}{\delta}\right)
\qquad\text{for }0<\delta\le\delta_*.
\]
Together with \eqref{eq:tan mean penalty} and $t\asymp\delta$, this lower bound proves \eqref{eq:tan pair probability claim} for $0<\delta\le\delta_*$.

For $\delta_*\le\delta<1$, Lemma~\ref{lem:limit theorem for first and second moment}(\ref{lem:limit theorem for first and second moment:II}) implies $\pi_{n,I,J}(x,y;B)\le C_dn^{-d}$ uniformly for $B$ satisfying $B_I\in\mathcal K_x$ and $B_J\in\mathcal K_y$. In this range, $n(\alpha_x^2+\beta_x^2+\rho(x,y)^2/\delta)$ is bounded by a constant depending only on $d$. Enlarging $C_d$ proves \eqref{eq:tan pair probability claim} in this range.
\end{proof}

\begin{lemma}\label{lem:tangential one row}
Fix $d\ge4$. There exist an integer $n_0$ and constants $t_0,c_d,C_d>0$, all depending only on $d$, such that the following holds for every $x,y\in T_{n_0}$ with $0<t=\norm{y-x}\le t_0$. For $a\sim N(0,I_d)$, let
\[
V=Z_y(a)-Z_x(a),\qquad\text{and}\qquad W=\begin{pmatrix}Z_x(a)\\V\end{pmatrix}.
\]
\begin{enumerate}[label=\textnormal{(\Roman*)},ref=\Roman*]
\item\label{lem:tangential one row:I} We have
\[
\E\norm{V}^j\le C_dt\qquad\text{for }j=1,2,3.
\]
In particular, $\norm{\E V}\le C_dt$.

\item\label{lem:tangential one row:II} There is a constant $\lambda>0$, depending only on $d$, such that, for every $k\ge1$, independent copies $a_1,\ldots,a_k$ of $a$, and every $R\in\R$,
\[
\PP\left(\sum_{\ell=1}^kL(a_\ell)\le R\right)
\le e^{\lambda R-c_dkq},
\]
where
\[
L(a)=\mu(x)^\top Z_x(a)+\frac{[\mu(y)-\mu(x)]^\top V}{t},
\quad\text{and}\quad
q=\E L(a)=\norm{\mu(x)}^2+\frac{\norm{\mu(y)-\mu(x)}^2}{t}.
\]

\end{enumerate}
\end{lemma}
\begin{proof}
By rotational invariance of the standard Gaussian measurements, we may assume without loss of generality that $\xs=e$, where $e$ is the first standard basis vector.

For (\ref{lem:tangential one row:I}), put $w=(y-x)/t$. Let $s_x(a)=\sign((a^\top x)^2-(a^\top e)^2)$, and similarly define $s_y$. Then
\begin{equation}\label{eq:tan increment decomposition}
V=t s_y(a)(a^\top w)a+[s_y(a)-s_x(a)](a^\top x)a.
\end{equation}
Thus
\begin{equation}\label{eq:tan increment pointwise bound}
\norm{Z_x(a)}\le C_d\norm{a}^2,
\qquad\text{and}\qquad
\norm{V}\le C_d\norm{a}^2\bigl(t+1_{\{s_x\ne s_y\}}\bigr).
\end{equation}
To bound the sign-change term in \eqref{eq:tan increment pointwise bound}, for $\sigma\in\{-1,1\}$ define
\[
A_{\sigma,t}=\{a:(a^\top b_\sigma)(a^\top b_\sigma+t a^\top w)<0\},
\qquad\text{and}\qquad b_\sigma=x-\sigma e.
\]
The factorization
\[
(a^\top x)^2-(a^\top e)^2
=\bigl(a^\top(x-e)\bigr)\bigl(a^\top(x+e)\bigr)
\]
and $y=x+tw$ imply that, if $s_x(a)\ne s_y(a)$, at least one of these two factors changes sign. Thus $\{s_x\ne s_y\}\subseteq A_{+,t}\cup A_{-,t}$ up to Gaussian null sets.

If $a\in A_{\sigma,t}$, continuity along the segment from $b_\sigma$ to $b_\sigma+tw$ implies $a^\top(b_\sigma+st w)=0$ for some $s\in[0,1]$. Writing $a=\norm{a}\omega$, we obtain $\abs{\omega^\top b_\sigma}\le t$. Since $\norm{b_\sigma}$ is bounded below uniformly on $T_{n_0}$, a direct computation of the surface measure implies $\PP(A_{\sigma,t})\le C_dt$. Therefore,
\[
\PP(s_x(a)\ne s_y(a))
\le\PP(A_{+,t})+\PP(A_{-,t})
\le C_dt.
\]
The events $A_{\sigma,t}$ depend only on $\omega$, which is independent of $\norm{a}$. Taking $t_0\le1$ and using finite Gaussian moments, we obtain from \eqref{eq:tan increment pointwise bound}, for $j=1,2,3$,
\begin{equation}\label{eq:tan increment moment bound}
\E\norm{V}^j
\le C_d\left(t^j\E\norm{a}^{2j}
+\E\bigl[\norm{a}^{2j}1_{A_{+,t}\cup A_{-,t}}\bigr]\right)
\le C_dt.
\end{equation}
In particular, $\norm{\E V}\le\E\norm{V}\le C_dt$, which proves (\ref{lem:tangential one row:I}).

To prove (\ref{lem:tangential one row:II}), we use $e^{-s}\le1-s+s^2e^{\abs{s}}/2$ to estimate $\E e^{-\lambda L(a)}$, so we first bound $\E[L(a)^2e^{\lambda_0\abs{L(a)}}]$ for some fixed $\lambda_0>0$.
Compactness of the tube and \eqref{eq:tan increment moment bound} imply $\norm{\mu(x)}+\norm{\mu(y)-\mu(x)}/t\le C_d$.
The coefficient bound implies
\begin{align*}
\abs{L(a)}
&\le \norm{\mu(x)}\norm{Z_x(a)}
+\frac{\norm{\mu(y)-\mu(x)}}{t}\norm{V}\\
&\le\left(\norm{\mu(x)}+
\frac{\norm{\mu(y)-\mu(x)}}{t}\right)\norm{W}
\le C_d\norm{W}.
\end{align*}
For the exponential weight, the independence of $\omega$ and $\norm{a}$, together with the $\chi_d^2$ distribution of $\norm{a}^2$, implies
\begin{equation}\label{eq:tan weighted crossing}
\E\left[(1+\norm{a})^6e^{\norm{a}^2/4}1_{A_{\sigma,t}}\right]\le C_dt.
\end{equation}
Choose $\lambda_0$ small enough that $\lambda_0\norm{W}\le\norm{a}^2/4$. Then \eqref{eq:tan increment pointwise bound} and \eqref{eq:tan weighted crossing} imply
\begin{equation}\label{eq:tan exponential moments}
\E\left[\norm{Z_x(a)}^2e^{\lambda_0\norm{W}}\right]\le C_d,
\qquad\text{and}\qquad
\E\left[\norm{V}^2e^{\lambda_0\norm{W}}\right]\le C_dt.
\end{equation}
After reducing $\lambda_0$ if necessary, \eqref{eq:tan exponential moments} remains valid with the weight $e^{\lambda_0\abs{L(a)}}$. Using $\abs{a+b}^2\le2(\abs{a}^2+\abs{b}^2)$, we obtain
\begin{align*}
\E\left[L(a)^2e^{\lambda_0\abs{L(a)}}\right]
&\le 2\norm{\mu(x)}^2
\E\left[\norm{Z_x(a)}^2e^{\lambda_0\abs{L(a)}}\right]\\
&\quad+\frac{2\norm{\mu(y)-\mu(x)}^2}{t^2}
\E\left[\norm{V}^2e^{\lambda_0\abs{L(a)}}\right]\\
&\le C_d\norm{\mu(x)}^2+
C_d\frac{\norm{\mu(y)-\mu(x)}^2}{t^2}\,t\\
&=C_d\left(\norm{\mu(x)}^2+
\frac{\norm{\mu(y)-\mu(x)}^2}{t}\right)=C_dq.
\end{align*}
Since $\E L(a)=q$, for a sufficiently small fixed $0<\lambda\le\lambda_0$, we obtain
\[
\E e^{-\lambda L(a)}
\le1-\lambda q+C_d\lambda^2q
\le e^{-c_dq}.
\]
By Markov's inequality, we have
\begin{align*}
\PP\left(\sum_{\ell=1}^kL(a_\ell)\le R\right)
&\le e^{\lambda R}\E\exp\left(-\lambda\sum_{\ell=1}^kL(a_\ell)\right)\\
&=e^{\lambda R}\prod_{\ell=1}^k\E e^{-\lambda L(a_\ell)}
\le e^{\lambda R-c_dkq},
\end{align*}
which proves (\ref{lem:tangential one row:II}).
\end{proof}

\begin{lemma}[Covariance and concentration in the tangential region]\label{lem:tangential concentration}
Fix $d\ge4$, and let $n_0,t_0$ be as in Lemma~\ref{lem:tangential one row}. There exists $0<t_1\le t_0$, depending only on $d$, such that, for every $x,y\in T_{n_0}$ satisfying $\norm{y-x}\le t_1$ and
\[
u_x^\top u_y>0,\qquad
0<\delta(x,y)<1,\qquad\text{and}\qquad
\rho(x,y)\le\delta(x,y),
\]
the following conclusions hold, where $t=\norm{y-x}$.
\begin{enumerate}[label=\textnormal{(\Roman*)},ref=\Roman*]
\item\label{lem:tangential concentration:I} For $a\sim N(0,I_d)$, let $V=Z_y(a)-Z_x(a)$ and $W=(Z_x(a)^\top,V^\top)^\top$. After increasing $n_0$ if necessary, there is a constant $c_d>0$, depending only on $d$, such that
\begin{equation}\label{eq:tan one row covariance}
\operatorname{Cov}(W)\succeq c_d\begin{pmatrix}I_d&0\\0&tI_d\end{pmatrix}.
\end{equation}

\item\label{lem:tangential concentration:II} Let $C_k$ contain $k$ independent standard Gaussian rows. After increasing $n_0$ if necessary, for every fixed $R>0$ we have
\[
\sup_{z\in\R^{2d}}\PP\left(\begin{pmatrix}S(x;C_k)\\S(y;C_k)-S(x;C_k)\end{pmatrix}\in z+B(0,R)\right)
\le C_{d,R}k^{-d}t^{-d/2}
\qquad\text{for all }k\ge1.
\]
The constant $C_{d,R}$ is uniform over $x,y,k$ satisfying these assumptions.
\end{enumerate}
\end{lemma}
\begin{proof}
By rotational invariance of the standard Gaussian measurements, we may assume without loss of generality that $\xs=e$, where $e$ is the first standard basis vector.

To prove the covariance bound \eqref{eq:tan one row covariance} in (\ref{lem:tangential concentration:I}), it suffices to show that, for every $\xi,\eta\in\R^d$,
\[
\operatorname{Var}\bigl(\xi^\top Z_x(a)+\eta^\top V\bigr)
\ge c_d\bigl(\norm{\xi}^2+t\norm{\eta}^2\bigr).
\]
Expanding the variance, we have
\begin{equation}\label{eq:tan covariance decomposition}
\operatorname{Var}\bigl(\xi^\top Z_x(a)+\eta^\top V\bigr)
=\xi^\top\Sigma_x\xi+\eta^\top\operatorname{Cov}(V)\eta+2\xi^\top\operatorname{Cov}(Z_x(a),V)\eta.
\end{equation}
By continuity of $\Sigma_x$ near the compact saddle sphere and \eqref{eq:tube covariance formula} with $x=\tau u$, choosing $n_0$ sufficiently large ensures $c_dI_d\preceq\Sigma_x\preceq C_dI_d$ uniformly on $T_{n_0}$. By Young's inequality, it therefore suffices to prove
\begin{equation}\label{eq:tan covariance targets}
\operatorname{Cov}(V)\succeq c_dtI_d,
\qquad\text{and}\qquad
\norm{\operatorname{Cov}(Z_x(a),V)}_{\mathrm{op}}\le C_dt,
\end{equation}
and then choose $t_1$ sufficiently small.

The first estimate gives a lower bound of order $t$ for the variance of $V$ in every unit direction. The second controls its interaction with $Z_x(a)$, allowing us to establish (\ref{lem:tangential concentration:I}). We then combine this bound with the moment estimates of Lemma~\ref{lem:tangential one row} to prove the concentration bound in (\ref{lem:tangential concentration:II}).

\paragraph{Increment covariance.}
Lemma~\ref{lem:tangential one row}(\ref{lem:tangential one row:I}) implies
\[
\operatorname{Var}(\zeta^\top V)
=\E(\zeta^\top V)^2-(\zeta^\top\E V)^2
\ge\E(\zeta^\top V)^2-C_dt^2
\qquad\text{for }\norm{\zeta}=1.
\]
Thus it suffices to find a uniform positive lower bound for $t^{-1}\E(\zeta^\top V)^2$ over unit vectors $\zeta$. The contribution of order $t$ comes from individual sign crossings; the noncrossing and simultaneous-crossing terms are of order $t^2$. To decompose this second moment, put $w=(y-x)/t$, which is a unit vector, let $s_x(a)=\sign((a^\top x)^2-(a^\top e)^2)$, and similarly define $s_y$. By \eqref{eq:tan increment decomposition}, on $\{s_x=s_y\}$ we have $V=t s_x(a)(a^\top w)a$, while on $\{s_x\ne s_y\}$ we have $V=-s_x(a)(2a^\top x+t a^\top w)a$.

For $\sigma\in\{-1,1\}$, recall the events from the proof of Lemma~\ref{lem:tangential one row}:
\[
A_{\sigma,t}=\{a:(a^\top b_\sigma)(a^\top b_\sigma+t a^\top w)<0\},
\qquad\text{and}\qquad b_\sigma=x-\sigma e.
\]
The factorization
\[
(a^\top x)^2-(a^\top e)^2
=\bigl(a^\top(x-e)\bigr)\bigl(a^\top(x+e)\bigr)
\]
and $y=x+tw$ imply that, if $s_x(a)\ne s_y(a)$, at least one of these two factors changes sign. Thus $\{s_x\ne s_y\}\subseteq A_{+,t}\cup A_{-,t}$ up to Gaussian null sets. If both factors change sign, their product keeps its sign. Hence $\{s_x\ne s_y\}$ occurs if and only if exactly one of $A_{+,t}$ and $A_{-,t}$ occurs. Therefore,
\begin{equation}\label{eq:tan crossing second moment decomposition}
\begin{aligned}
\E(\zeta^\top V)^2
&=t^2\E\left[(a^\top w)^2(a^\top\zeta)^2 1_{\{s_x=s_y\}}\right]\\
&\quad+\sum_{\sigma=\pm1}\E\left[(2a^\top x+t a^\top w)^2(a^\top\zeta)^2 1_{A_{\sigma,t}}\right]\\
&\quad-2\E\left[(2a^\top x+t a^\top w)^2(a^\top\zeta)^2 1_{A_{+,t}\cap A_{-,t}}\right].
\end{aligned}
\end{equation}
The first term in \eqref{eq:tan crossing second moment decomposition} is $O_d(t^2)$ because Gaussian fourth moments are finite.

To bound the intersection term, we first check the geometry of $b_+$ and $b_-$.
Write $x=\alpha e+(\tau+\beta)u$ and $H_\sigma=b_\sigma^\perp$. A direct computation shows
\begin{equation}\label{eq:tan hyperplane geometry}
\norm{P_{H_\sigma}x}=\frac{\tau+\beta}{\norm{b_\sigma}},
\qquad\text{and}\qquad
\norm{P_{H_\sigma}b_{-\sigma}}=2\norm{P_{H_\sigma}x}.
\end{equation}
Thus the normal vectors $b_+=x-e$ and $b_-=x+e$ of $H_+$ and $H_-$, respectively, are uniformly nonparallel, their lengths are bounded above and below, and $\norm{P_{H_\sigma}x}$ is bounded below.

To bound simultaneous crossings, write $a=\norm{a}\omega$ and choose orthonormal vectors $v_1=b_+/\norm{b_+}$ and $v_2=P_{H_+}b_-/\norm{P_{H_+}b_-}$. On $A_{+,t}\cap A_{-,t}$, the defining inequalities of $A_{+,t}$ and $A_{-,t}$ imply
\[
\abs{\omega^\top v_1}\le C_dt,
\qquad\text{and}\qquad
\abs{\omega^\top v_2}
=\frac{\abs{\omega^\top b_--(b_-^\top v_1)(\omega^\top v_1)}}{\norm{P_{H_+}b_-}}
\le C_dt,
\]
where the bounds on $\norm{b_+}$ and $\norm{P_{H_+}b_-}$ follow from \eqref{eq:tan hyperplane geometry}. The pair $(\omega^\top v_1,\omega^\top v_2)$ has density proportional to $(1-z_1^2-z_2^2)^{(d-4)/2}$ on the unit disk. Since $d\ge4$, this density is bounded, and hence $\PP(A_{+,t}\cap A_{-,t})\le C_dt^2$.
Independence of $\omega$ and $\norm{a}$ then implies
\begin{equation}\label{eq:tan weighted simultaneous crossing}
\E\left[(1+\norm{a})^6e^{\norm{a}^2/4}1_{A_{+,t}\cap A_{-,t}}\right]\le C_dt^2.
\end{equation}
Since $(2a^\top x+t a^\top w)^2(a^\top\zeta)^2\le C_d\norm{a}^4$, \eqref{eq:tan weighted simultaneous crossing} bounds the intersection term in \eqref{eq:tan crossing second moment decomposition} by $O_d(t^2)$. Consequently,
\begin{equation}\label{eq:tan crossing second moment expansion}
\E(\zeta^\top V)^2
=\sum_{\sigma=\pm1}\E\left[(2a^\top x+t a^\top w)^2(a^\top\zeta)^2 1_{A_{\sigma,t}}\right]+O_d(t^2).
\end{equation}

By \eqref{eq:tan crossing second moment expansion}, it now suffices to compute the limit of the sum divided by $t$ and show that it is uniformly positive. For each $\sigma\in\{-1,1\}$ in this sum, let $b=b_\sigma$, put $\widehat b=b/\norm{b}$, and write $a=z+v\widehat b$ with $z\perp b$. For sufficiently small $t$ and fixed $z$, the crossing interval in $v$ has endpoints $0$ and
\[
v_t(z)=-\frac{t z^\top w}{\norm{b}+t \widehat b^\top w}.
\]
Thus, for a continuous function $F$ bounded by a polynomial times a Gaussian,
\begin{equation}\label{eq:tan moving interval}
\int_{\{(a^\top b)(a^\top(b+tw))<0\}}F(a)\dd a
=\int_{b^\perp}\int_0^1F(z+\ell v_t(z)\widehat b)\abs{v_t(z)}\dd\ell\,\dd z.
\end{equation}
Here and below, $\dd z$ denotes Lebesgue measure in orthonormal coordinates on the hyperplane.
Moreover, $\abs{v_t(z)}\le C_dt\norm{z}$ and $\norm{z+\ell v_t(z)\widehat b}\ge\norm{z}$. Dividing by $t$ and applying dominated convergence gives the leading contribution of this crossing:
\begin{equation}\label{eq:tan crossing limit}
\frac1t\int_{\{(a^\top b)(a^\top(b+tw))<0\}}F(a)\dd a
\longrightarrow\frac1{\norm{b}}\int_{b^\perp}\abs{z^\top w}F(z)\dd z.
\end{equation}
We also need \eqref{eq:tan crossing limit} to hold uniformly as $x,w,\zeta$ vary. The same dominated-convergence argument applies along convergent sequences of parameters $b,w,F$ with a common Gaussian bound, using orthonormal coordinates on $b^\perp$ that vary continuously with $b$.
Applying \eqref{eq:tan crossing limit} to $(2a^\top x+t a^\top w)^2(a^\top\zeta)^2\varphi_{\scriptscriptstyle I_d}(a)$ and using \eqref{eq:tan crossing second moment expansion}, we therefore obtain
\begin{equation}\label{eq:tan crossing quadratic form}
\frac1t\E(\zeta^\top V)^2
\longrightarrow4\sum_{\sigma=\pm1}\frac1{\norm{b_\sigma}}\int_{H_\sigma}\abs{z^\top w}(z^\top x)^2(z^\top\zeta)^2\varphi_{\scriptscriptstyle I_d}(z)\dd z.
\end{equation}

To obtain the covariance lower bound, we show that the limit in \eqref{eq:tan crossing quadratic form} is uniformly positive. Put $\delta=\delta(x,y)$. Lemma~\ref{lem:delta rho distance identity} and $\rho(x,y)\le\delta$ imply
\[
t^2=\rho(x,y)^2+(\tau+\beta_x)(\tau+\beta_y)\delta^2
\le C_d\delta^2.
\]
For $\sigma=\pm1$, \eqref{eq:shared factor first bound} and \eqref{eq:density matching sign bound}, together with $\norm{x-\sigma e}\le C_d$, therefore imply
\begin{equation}\label{eq:tan transversality bound}
\norm{P_{H_\sigma}w}
=\frac{\norm{P_{(x-\sigma e)^\perp}(y-x)}}{t}
\ge\frac{c_d\delta}{t}\ge c_d.
\end{equation}
On each $H_\sigma$, \eqref{eq:tan transversality bound} and \eqref{eq:tan hyperplane geometry} imply that the restrictions of $z\mapsto z^\top w$ and $z\mapsto z^\top x$ are nonzero linear forms. Thus the corresponding integral vanishes exactly when $\zeta$ is parallel to $b_\sigma$. Since $b_+$ and $b_-$ are nonparallel, the sum of the two integrals is positive for every unit $\zeta$.

The triples $(x,\zeta,w)$, with $x\in T_{n_0}$ and unit vectors $\zeta,w$ satisfying the fixed positive lower bound in \eqref{eq:tan transversality bound}, form a compact set. The right side of \eqref{eq:tan crossing quadratic form} is continuous there, by Gaussian domination in local orthonormal coordinates, so it has a positive minimum depending only on $d$.

The representation \eqref{eq:tan moving interval} implies convergence along every convergent sequence of parameters, hence uniform convergence in \eqref{eq:tan crossing quadratic form}. The positive lower bound therefore transfers to $t^{-1}\E(\zeta^\top V)^2$ for sufficiently small $t$. Since $\norm{\E V}\le C_dt$, after choosing $t_1\le t_0$ sufficiently small, we obtain, uniformly over unit $\zeta$,
\[
\operatorname{Var}(\zeta^\top V)
=\E(\zeta^\top V)^2-(\zeta^\top\E V)^2
\ge c_dt-C_dt^2
\ge \frac{c_d}{2}t.
\]
Lemma~\ref{lem:tangential one row}(\ref{lem:tangential one row:I}) with $j=2$ implies $\E\norm{V}^2\le C_dt$, which proves the upper covariance bound. Thus
\begin{equation}\label{eq:tan increment covariance}
c_dtI_d\preceq\operatorname{Cov}(V)\preceq C_dtI_d.
\end{equation}
This proves the first estimate in \eqref{eq:tan covariance targets}.

\paragraph{Mixed covariance.}
Since $\norm{Z_x(a)}\le C_d\norm{a}^2$ and $\norm{\mu(x)}\le C_d$, \eqref{eq:tan increment pointwise bound} and \eqref{eq:tan weighted crossing} imply
\begin{equation}\label{eq:tan cross covariance bound}
\begin{aligned}
\norm{\operatorname{Cov}(Z_x(a),V)}_{\mathrm{op}}
&\le \E[\norm{Z_x(a)}\norm{V}]+\norm{\mu(x)}\norm{\E V}\\
&\le C_d\left(t\E\norm{a}^4
+\E\bigl[\norm{a}^4 1_{A_{+,t}\cup A_{-,t}}\bigr]\right)+C_dt
\le C_dt.
\end{aligned}
\end{equation}
This proves the second estimate in \eqref{eq:tan covariance targets}.
Substituting \eqref{eq:tan increment covariance} and \eqref{eq:tan cross covariance bound} into \eqref{eq:tan covariance decomposition} and applying Young's inequality, we obtain
\[
\begin{aligned}
\operatorname{Var}(\xi^\top Z_x(a)+\eta^\top V)
&\ge c_d\norm{\xi}^2+c_dt\norm{\eta}^2-C_dt\norm{\xi}\norm{\eta}\\
&\ge\frac{c_d}{2}\norm{\xi}^2+(c_dt-C_dt^2)\norm{\eta}^2\\
&\ge c_d(\norm{\xi}^2+t\norm{\eta}^2),
\end{aligned}
\]
where the last inequality holds after reducing $t_1$ once more. This proves (\ref{lem:tangential concentration:I}).

\paragraph{Concentration estimate.}
For any fixed $T>0$, Lemma~\ref{lem:Esseen concentration inequality} with $m=2d$ and $\varepsilon=RT$, together with independence, implies
\begin{equation}\label{eq:tan concentration reduction}
\sup_{z\in\R^{2d}}\PP(X\in z+B(0,R))
\le C_{2d,RT}R^{2d}\int_{\norm{(\xi,\eta)}\le T}
\abs{\chi(\xi,\eta)}^k\dd\xi\,\dd\eta,
\end{equation}
where
\[
X=\begin{pmatrix}S(x;C_k)\\S(y;C_k)-S(x;C_k)\end{pmatrix},
\qquad\text{and}\qquad
\chi(\xi,\eta)=\E e^{i[\xi^\top Z_x(a)+\eta^\top V]}.
\]
It therefore suffices to prove that, for some fixed $T>0$ depending only on $d$,
\begin{equation}\label{eq:tan small frequency}
\abs{\chi(\xi,\eta)}\le e^{-c_d(\norm{\xi}^2+t\norm{\eta}^2)}
\qquad\text{when }\norm{\xi},\norm{\eta}\le T,
\end{equation}
since the resulting Gaussian integral in \eqref{eq:tan concentration reduction} has order $k^{-d}t^{-d/2}$.

To prove \eqref{eq:tan small frequency}, we separate the variance from the third-order remainder.
Let $Y=\xi^\top[Z_x(a)-\mu(x)]+\eta^\top[V-\E V]$, and let $Y'$ be an independent copy. We use the same symmetrization idea as in \cite[proof of Proposition~8.2]{bobkov2023refinements}. Independence and the inequality $\cos s\le1-s^2/2+\abs{s}^3/6$ imply \[ \abs{\chi(\xi,\eta)}^2=\E\cos(Y-Y')\le1-\operatorname{Var}(Y)+\frac16\E\abs{Y-Y'}^3. \]
For the remainder, Jensen's inequality, $\abs{a+b}^3\le4(\abs{a}^3+\abs{b}^3)$, and \eqref{eq:tan increment pointwise bound} imply
\begin{align*}
\E\abs{Y-Y'}^3
&\le8\E\abs{Y}^3\\
&\le32\norm{\xi}^3\E\norm{Z_x(a)-\mu(x)}^3
+32\norm{\eta}^3\E\norm{V-\E V}^3\\
&\le C_d\left(\norm{\xi}^3\E\norm{a}^6
+\norm{\eta}^3\E\norm{V}^3\right)\\
&\le C_d(\norm{\xi}^3+t\norm{\eta}^3),
\end{align*}
where the last inequality uses Lemma~\ref{lem:tangential one row}(\ref{lem:tangential one row:I}) with $j=3$ and the finiteness of Gaussian sixth moments.
Since part (\ref{lem:tangential concentration:I}) implies $\operatorname{Var}(Y)\ge c_d(\norm{\xi}^2+t\norm{\eta}^2)$, we obtain \[ \abs{\chi(\xi,\eta)}^2\le1-c_d(\norm{\xi}^2+t\norm{\eta}^2)+C_d(\norm{\xi}^3+t\norm{\eta}^3). \]
For $\norm{\xi},\norm{\eta}\le T$, we have ${\norm{\xi}^3+t\norm{\eta}^3\le T(\norm{\xi}^2+t\norm{\eta}^2)}$. Choosing $T>0$ sufficiently small depending only on $d$ therefore allows us to absorb the cubic terms into the quadratic terms. Taking square roots and using $1-s\le e^{-s}$ proves \eqref{eq:tan small frequency}.
Finally, substituting \eqref{eq:tan small frequency} into \eqref{eq:tan concentration reduction} and enlarging the integration domain, we obtain
\begin{align*}
\sup_{z\in\R^{2d}}\PP(X\in z+B(0,R))
&\le C_{d,R}\int_{\R^d}\int_{\R^d}
e^{-c_dk(\norm{\xi}^2+t\norm{\eta}^2)}\dd\xi\,\dd\eta\\
&\le C_{d,R}k^{-d/2}(kt)^{-d/2}
=C_{d,R}k^{-d}t^{-d/2}.
\end{align*}
This proves (\ref{lem:tangential concentration:II}).
\end{proof}

%%%%%%%%%%%%%%%%%%%%%%%%%%%%%%%%%% Normal dominated region
\subsection{Normal-dominated region}
We continue to use the notation $\Omega^{\mathrm{nor}}$ and $D_{r,n}^{\mathrm{nor}}$ from the proof of Theorem~\ref{thm:second moment asymptotic appendix}.

\begin{lemma}[Normal-dominated region]\label{lem:normal dominated region}
Fix $d\ge4$ and $0<\delta_0<1$. As $n\to\infty$, we have \[ D_{r,n}^{\mathrm{nor}}=O_d(n^{d/2-1})=o_d(n^{d-2}) \qquad\text{for }0\le r\le d-1. \]
\end{lemma}
\begin{proof}
By rotational invariance of the standard Gaussian measurements, we may assume without loss of generality that $\xs=e$, where $e$ is the first standard basis vector.

There are $\binom nd\binom dr\binom{n-d}{d-r}=O_d(n^{2d-r})$ ordered pairs $(I,J)$ with overlap $r$. It therefore suffices to show that, for any fixed $I\ne J$ with $\abs{I}=\abs{J}=d$ and $\abs{I\cap J}=r$, the integral
\begin{equation}\label{eq:normal pair integral}
\begin{aligned}
I_{r,n}^{\mathrm{nor}}
&:=\int_{\Omega^{\mathrm{nor}}}p_{X_I(x),X_J(y)}(0,0)\\
&\quad\times\E\left[\begin{aligned}
&\abs{\det A_I}\abs{\det A_J}1_{\mathcal K_x}(A_I)1_{\mathcal K_y}(A_J)\\
&\quad\times\pi_{n,I,J}(x,y;A_{I\cup J})
\end{aligned}\;\middle|\;X_I(x)=X_J(y)=0\right]\dd x\,\dd y
\end{aligned}
\end{equation}
satisfies
\begin{equation}\label{eq:normal pair bound}
I_{r,n}^{\mathrm{nor}}\le C_dn^{-3d/2-1+r}\qquad\text{for }0\le r\le d-1.
\end{equation}
We first estimate the integral over
\[
\Omega_+^{\mathrm{nor}}=\Omega^{\mathrm{nor}}\cap\{(x,y):u_x^\top u_y>0,\ \delta(x,y)>0\}.
\]
Denote this contribution by $I_+$.
By Lemma~\ref{lem:normal pi}, for pairs in $\Omega_+^{\mathrm{nor}}$,
\begin{equation}\label{eq:normal pair probability claim}
\sup_{\substack{B:B_I\in\mathcal K_x\\B_J\in\mathcal K_y}}\pi_{n,I,J}(x,y;B)\le C_dn^{-d/2}e^{-c_dn\norm{y-x}}.
\end{equation}
When $A_I\in\mathcal K_x$ and $A_J\in\mathcal K_y$, $\abs{\det A_I}\abs{\det A_J}\le M_d^{2d}$. Therefore, \eqref{eq:normal pair probability claim} implies \[ I_+\le C_dn^{-d/2}\int_{\Omega_+^{\mathrm{nor}}}p_{X_I(x),X_J(y)}(0,0)e^{-c_dn\norm{y-x}}\dd x\,\dd y. \]

Write $h=y-x$. Lemma~\ref{lem:delta rho distance identity} implies \[ \norm{h}^2=\rho(x,y)^2+(\tau+\beta_x)(\tau+\beta_y)\delta(x,y)^2>0, \]
where positivity follows from $\rho(x,y)>\delta(x,y)$ and $\tau+\beta_x,\tau+\beta_y>0$ for sufficiently large $n$.
By \eqref{eq:density second estimate}, \[ p_{X_I(x),X_J(y)}(0,0)\le \frac{C_d}{\norm{P_{v_+^\perp}h}^{r_+}\norm{P_{v_-^\perp}h}^{r_-}}, \]
where $v_\pm=x\mp e$ and $r_\pm$ count the shared rows with matching signs $\pm1$. Put $q=r_++r_-\le r$. The chosen sign pattern implies $r_+\le d-2$ and $r_-\le2\le d-2$. The disjoint case is included by taking $r_+=r_-=0$.

After enlarging the displacement domain to $\R^d\setminus\{0\}$ and changing variables $h=t\omega$, we obtain
\[
I_+\le C_dn^{-d/2}\int_{T_n}
\underbrace{\int_{\mathbb S^{d-1}}\frac{\dd\sigma(\omega)}{\norm{P_{v_+^\perp}\omega}^{r_+}\norm{P_{v_-^\perp}\omega}^{r_-}}}_{I_1}
\underbrace{\int_0^\infty e^{-c_dnt}t^{d-1-q}\dd t}_{I_2}
\dd x.
\]
By Lemma~\ref{lem:density}(\ref{lem:density:III}), the first integral satisfies $I_1\le C_d$ uniformly in $x$. For the second integral, the change of variables $s=nt$ implies
\begin{equation}\label{eq:normal radial integral bound}
I_2=n^{-(d-q)}\int_0^\infty e^{-c_ds}s^{d-1-q}\dd s\le C_dn^{-(d-q)},
\end{equation}
since $q\le r\le d-1$. Substituting the bounds on $I_1$ and $I_2$, we obtain
\begin{align}
I_+&\le C_dn^{-d/2}n^{-(d-q)}\int_{T_n}\dd x\notag\\
&=C_dn^{-3d/2+q}\int_{\mathbb S^{d-2}}\int_{[-n^{-1/2},n^{-1/2}]^2}
(\tau+\beta_x)^{d-2}\dd\alpha_x\,\dd\beta_x\,\dd\sigma(u_x)\notag\\
&\le C_dn^{-3d/2-1+q}\le C_dn^{-3d/2-1+r},\label{eq:integral for omega}
\end{align}
where we used the change of variables $x=\alpha_xe+(\tau+\beta_x)u_x$ and $q\le r$.

\paragraph{The antipodal contribution.}
For pairs with $u_x^\top u_y<0$, replace $y$ by $y'=-y$ and reverse all active signs assigned to $J$, as at the end of the proof of Lemma~\ref{lem:tangential region}. This change preserves $\delta(x,y)$, $\rho(x,y)$, and Lebesgue measure, and transforms the antipodal contribution into the corresponding integral over $\Omega_+^{\mathrm{nor}}$ with reversed signs for $J$.
The proof of Lemma~\ref{lem:normal pi} applies under this sign reversal, so \eqref{eq:normal pair probability claim} holds with $y'$ in place of $y$. The density bound \eqref{eq:density second estimate} and the angular bound in Lemma~\ref{lem:density}(\ref{lem:density:III}) also remain valid: the signs on $I$ are unchanged, so $r_+\le d-2$ and $r_-\le2\le d-2$. Hence the antipodal contribution satisfies \eqref{eq:integral for omega} as well. Since the sets $\{\delta(x,y)=0\}$ and $\{u_x^\top u_y=0\}$ have Lebesgue measure zero, combining the two contributions proves \eqref{eq:normal pair bound}, uniformly over $I,J$, and completes the proof.
\end{proof}

\begin{lemma}[Probability estimate in the normal-dominated region]\label{lem:normal pi}
Fix $d\ge4$ and active sets $I\ne J$ with $\abs{I}=\abs{J}=d$. There are constants $C_d,c_d>0$ such that, for all sufficiently large $n$ and every $x,y\in T_n$ satisfying $u_x^\top u_y>0$ and $0<\delta(x,y)<\rho(x,y)$, \[ \sup_{\substack{B:B_I\in\mathcal K_x\\B_J\in\mathcal K_y}}\pi_{n,I,J}(x,y;B)\le C_dn^{-d/2}e^{-c_dn\norm{y-x}}. \]
\end{lemma}
\begin{proof}
Fix active rows $B$ satisfying $B_I\in\mathcal K_x$ and $B_J\in\mathcal K_y$, let $r=\abs{I\cap J}$ and $\nu=n-2d+r$, and relabel the remaining independent standard Gaussian rows as $g_1,\ldots,g_\nu$.
We split the remaining rows into two independent groups to obtain the two factors in the desired bound. The joint event forces at least one group's increment to have a projection bounded above by a fixed constant in the direction $\mu(y)-\mu(x)$, yielding the exponential factor $e^{-c_dn\norm{y-x}}$. Conditional on that group, a uniform concentration estimate for the other yields the factor $n^{-d/2}$.
Take the first $\nu_1=\lfloor\nu/2\rfloor$ rows for $C_1$ and the rest for $C_2$:
\[
C_1=(g_1,\ldots,g_{\nu_1})^\top,\qquad\text{and}\qquad
C_2=(g_{\nu_1+1},\ldots,g_\nu)^\top.
\]
With $\nu_2=\nu-\nu_1$, both groups have $n/2+O_d(1)$ rows. As in the proof of Lemma~\ref{lem:tangential region estimate}, the target sets for the sums at $x$ and $y$ are
\[
\mathcal B_x=-Q(x;B_I)(-1,1)^d-S(x;B_{J\setminus I}),
\qquad\text{and}\qquad
\mathcal B_y=-Q(y;B_J)(-1,1)^d-S(y;B_{I\setminus J}),
\]
and both lie in $B(0,R_d)$ by the row norm bounds on $B$. Write
\[
U_j=S(x;C_j),\qquad\text{and}\qquad V_j=S(y;C_j)-S(x;C_j)\qquad\text{for }j=1,2.
\]
Then $\pi_{n,I,J}(x,y;B)=\PP(F)$ for the pair event
\[
F=\{U_1+U_2\in\mathcal B_x,\quad U_1+U_2+V_1+V_2\in\mathcal B_y\}.
\]

Put $h=y-x$ and $w=[\mu(y)-\mu(x)]/\norm{\mu(y)-\mu(x)}$. For a single row $a\sim N(0,I_d)$, write $V=Z_y(a)-Z_x(a)$. Let $\lambda_0$ be as in Lemma~\ref{lem:normal increment moments}. For a sufficiently small fixed $0<\lambda\le\lambda_0$,
\begin{equation}\label{eq:normal increment exponential bound}
\begin{aligned}
\E e^{-\lambda w^\top V}
&\le 1-\lambda\E[w^\top V]
 +\frac{\lambda^2}{2}\E\left[(w^\top V)^2e^{\lambda\abs{w^\top V}}\right]\\
&\le 1-\lambda\norm{\mu(y)-\mu(x)}
 +\frac{\lambda^2}{2}\E\left[\norm{V}^2e^{\lambda_0\norm{V}}\right]\\
&\le 1-\lambda c_d\norm{h}+C_d\lambda^2\norm{h}\\
&\le e^{-c_d\norm{h}}.
\end{aligned}
\end{equation}
The first inequality in \eqref{eq:normal increment exponential bound} uses $e^{-t}\le1-t+t^2e^{\abs t}/2$; the third uses \eqref{eq:normal increment moment bounds}. The last inequality in \eqref{eq:normal increment exponential bound} follows by choosing $\lambda$ sufficiently small.
On $F$, the difference sum $V_1+V_2$ belongs to $\mathcal B_y-\mathcal B_x\subset B(0,2R_d)$. Thus $F$ implies at least one of the events $E_j=\{w^\top V_j\le R_d\}$. For $j\in\{1,2\}$, the exponential Markov inequality implies
\[
\PP(E_j)\le e^{\lambda R_d}
\left(\E e^{-\lambda w^\top[Z_y(a)-Z_x(a)]}\right)^{\nu_j}
\le C_de^{-c_dn\norm{h}}.
\]

The local limit argument for $S(x;C_j)$ in the proof of Lemma~\ref{lem:limit theorem for first and second moment}(\ref{lem:limit theorem for first and second moment:I}), applied to a fixed cube containing $B(0,R_d)$, implies
\begin{equation}\label{eq:normal sum concentration}
\sup_{z\in\R^d}\PP\bigl(S(x;C_j)\in z+B(0,R_d)\bigr)
\le C_d\nu_j^{-d/2}
\end{equation}
uniformly for $x\in T_n$ and all sufficiently large $n$.
For $j\in\{1,2\}$, conditional on $C_j$, the event $F\cap E_j$ requires $U_{3-j}\in\mathcal B_x-U_j$. Since $C_{3-j}$ is independent of $C_j$ and $\mathcal B_x\subset B(0,R_d)$, \eqref{eq:normal sum concentration} implies
\begin{align*}
\PP(F\cap E_j)
&=\E\bigl[1_{E_j}\PP(F\mid C_j)\bigr]\\
&\le\E\bigl[1_{E_j}\PP(U_{3-j}\in\mathcal B_x-U_j\mid C_j)\bigr]\\
&\le C_d\nu_{3-j}^{-d/2}\PP(E_j)
\le C_dn^{-d/2}\PP(E_j).
\end{align*}
Summing over $j=1,2$ proves \eqref{eq:normal pair probability claim}. The argument uses only the row bounds on $B$, so it also applies after reversing the active signs at an antipodal endpoint.
\end{proof}

\begin{lemma}[Moments of increments in the normal-dominated region]\label{lem:normal increment moments}
Fix $d\ge4$ and let $a\sim N(0,I_d)$. There exist $n_0$ and constants $c_d,C_d,\lambda_0>0$ depending only on $d$ such that, for every $n\ge n_0$ and $x,y\in T_n$ with $u_x^\top u_y>0$ and $\delta(x,y)<\rho(x,y)$,
\begin{equation}\label{eq:normal increment moment bounds}
\ificlr
\begin{aligned}
\norm{\mu(y)-\mu(x)}&\ge c_d\norm{y-x},\quad\text{and}\\
\E\left[\norm{Z_y(a)-Z_x(a)}^2e^{\lambda_0\norm{Z_y(a)-Z_x(a)}}\right]&\le C_d\norm{y-x}.
\end{aligned}
\else
\norm{\mu(y)-\mu(x)}\ge c_d\norm{y-x},\quad\text{and}\quad
\E\left[\norm{Z_y(a)-Z_x(a)}^2e^{\lambda_0\norm{Z_y(a)-Z_x(a)}}\right]\le C_d\norm{y-x}.
\fi
\end{equation}
\end{lemma}
\begin{proof}
By rotational invariance of the standard Gaussian measurements, we may assume without loss of generality that $\xs=e$, where $e$ is the first standard basis vector.
Put $h=y-x$.
Since $u_x^\top u_y>0$, we have
\[
\delta=\norm{u_x-u_y},\qquad\text{and}\qquad
\rho^2=(\alpha_y-\alpha_x)^2+(\beta_y-\beta_x)^2.
\]
Lemma~\ref{lem:delta rho distance identity} and the condition $\delta<\rho$ of the normal-dominated region imply
\[
\norm{h}^2=\rho^2+(\tau+\beta_x)(\tau+\beta_y)\delta^2
\le C_d\rho^2\le C_dn^{-1}.
\]
The segment from $x$ to $y$ therefore lies within $O_d(n^{-1/2})$ of $\tau u_x$. Since $\mu=\nabla\F/2$, Proposition~\ref{prop:population landscape} and continuity of $D\mu$ near $\mathcal S$ give, uniformly over these pairs,
\begin{equation}\label{eq:normal mean expansion}
\mu(y)-\mu(x)
=\kappa(u_xu_x^\top-ee^\top)h+o(\norm{h})
\qquad\text{with }
\kappa=\frac{8\tau}{\pi(1+\tau^2)^2}>0.
\end{equation}
Moreover,
\[
e^\top h=\alpha_y-\alpha_x,\qquad\text{and}\qquad
u_x^\top h=\beta_y-\beta_x-\frac{\tau+\beta_y}{2}\delta^2.
\]
Hence,
\begin{equation}\label{eq:normal projected displacement}
\begin{aligned}
\norm{(u_xu_x^\top-ee^\top)h}
&=\sqrt{(\alpha_y-\alpha_x)^2+
\left(\beta_y-\beta_x-\frac{\tau+\beta_y}{2}\delta^2\right)^2}\\
&\ge\rho-\frac{\tau+\beta_y}{2}\delta^2\\
&\ge\rho-C_d\rho^2
\ge\frac{\rho}{2}
\ge c_d\norm{h}.
\end{aligned}
\end{equation}
The first equality in \eqref{eq:normal projected displacement} uses that $e$ and $u_x$ are orthogonal unit vectors. The bound $\rho-C_d\rho^2\ge\rho/2$ in \eqref{eq:normal projected displacement} holds for sufficiently large $n$, since $\rho\le C_dn^{-1/2}$.
The expansion \eqref{eq:normal mean expansion} and \eqref{eq:normal projected displacement} therefore imply, for sufficiently large $n$,
\[
\norm{\mu(y)-\mu(x)}
\ge\kappa\norm{(u_xu_x^\top-ee^\top)h}-o(\norm{h})
\ge c_d\norm{h}.
\]

To prove the second bound in \eqref{eq:normal increment moment bounds}, write $s_x(a)=\sign((a^\top x)^2-(a^\top e)^2)$ and similarly for $s_y$. Then
\[
Z_y(a)-Z_x(a)
=s_y(a)(a^\top h)a+[s_y(a)-s_x(a)](a^\top x)a,
\]
so
\begin{equation}\label{eq:normal row increment bound}
\norm{Z_y(a)-Z_x(a)}
\le C_d\norm{a}^2\bigl(\norm{h}+1_{\{s_x(a)\ne s_y(a)\}}\bigr).
\end{equation}
If $s_x(a)\ne s_y(a)$, continuity of $(a^\top z)^2-(a^\top e)^2=a^\top(z-e)a^\top(z+e)$ along the segment from $x$ to $y$ implies that $a^\top(x+th-\sigma e)=0$ for some $t\in[0,1]$ and $\sigma\in\{-1,1\}$. Writing $a=\norm{a}\omega$, we obtain
\[
\abs{\omega^\top(x-\sigma e)}
=t\abs{\omega^\top h}\le\norm{h}.
\]
For each fixed $\sigma$, a direct computation of the surface measure of the spherical band, together with $\norm{x-\sigma e}\ge\tau/2$, shows that
\[
\PP\bigl(\omega\in\mathbb S^{d-1}:
\abs{\omega^\top(x-\sigma e)}\le\norm{h}\bigr)
\le C_d\norm{h}.
\]
Taking the union over $\sigma\in\{-1,1\}$ yields $\PP(s_x(a)\ne s_y(a))\le C_d\norm{h}$. This event depends only on $\omega$, which is independent of $\norm{a}$. Since $\norm{a}^2$ has a $\chi_d^2$ distribution, for $0<\lambda_0<1/(2C_d)$ we have
\[
\E\bigl[\norm{a}^4e^{C_d\lambda_0\norm{a}^2}\bigr]
=d(d+2)(1-2C_d\lambda_0)^{-d/2-2}<\infty.
\]
Thus, for sufficiently large $n$, \eqref{eq:normal row increment bound} implies
\[
\E\left[\norm{Z_y(a)-Z_x(a)}^2e^{\lambda_0\norm{Z_y(a)-Z_x(a)}}\right]
\le C_d(\norm{h}^2+\norm{h})\le C_d\norm{h}.
\]
\end{proof}

%%% Local Variables:
%%% mode: LaTeX
%%% TeX-master: "../main"
%%% End:

\section{Auxiliary lemmas}\label{app:auxiliary-lemmas}
We collect auxiliary lemmas used throughout the paper. The proof of the first lemma is simple, so we omit it.

\begin{lemma}\label{lem:zeros of polynomial}
Let $G\in\R^d$ be a nondegenerate Gaussian vector, and let $p:\R^d\to \R$ be a polynomial that is not identically zero. Then, \(\PP(p(G)=0)=0\).
\end{lemma}

\begin{lemma}[Almost sure nonsingularity of projected matrices]
	\label{lem:nonsingularity of projected hessian}
	Suppose $n> d$, and let $(a_i)_{i=1}^n$ be i.i.d.\ standard Gaussian vectors in $\R^d$. Then, almost surely, for every index set $I\subset [n]$ with $\abs{I}<d$, every sign pattern $\sigma\in \{-1,1\}^{I^c}$, and every matrix $V$ whose columns form an orthonormal basis of the orthogonal complement of $\operatorname{span}\{a_i\}_{i\in I}$,
	\[\text{the matrix}\quad
	V^\top \left(\sum_{i\notin I}\sigma_i a_i a_i^\top\right) V \quad \text{is nonsingular}.
	\]
\end{lemma}

\begin{proof}
	Fix an index set $I\subset[n]$ with $\abs{I}<d$ and a sign pattern $\sigma$. Almost surely, the orthogonal complement of \(\operatorname{span}\{a_i\}_{i\in I}\) equals that of \(\operatorname{span}\{u_i\}_{i\in I}\), where \(u_i=a_i/\norm{a_i}\).
	We can therefore choose an orthonormal basis matrix $V=(v_1,\ldots, v_{d-\abs{I}})$ as a measurable function of the directions $(u_i)_{i\in I}$ alone. Define \(H=\sum_{i\notin I}\sigma_i a_i a_i^\top\).
	Since any orthonormal basis can be written as \(
	\tilde{V}=VQ\)  with \(Q^\top Q=I
	\),
	the singularity of $V^\top H V$ does not depend on the choice of $V$.

	The vectors $(a_i)_{i\notin I}$ are independent of the directions $(u_i)_{i\in I}$. Conditional on these directions, $\det(V^\top H V)$ is a polynomial in the components of $(a_i)_{i\notin I}$. This polynomial is not identically zero. Indeed, choose distinct inactive indices $j_1,\ldots, j_{d-\abs{I}}$ from $I^c$ and set \[
	a_{j_l}=v_l,\qquad\text{and}\qquad a_j=0 \quad \text{for all remaining indices in }I^c.
	\]
	For this choice, $V^\top H V=\diag(\sigma_{j_1},\ldots, \sigma_{j_{d-\abs{I}}})$, so its determinant is nonzero. Lemma~\ref{lem:zeros of polynomial} gives
	\[ \PP(\det(V^\top H V)=0\mid(u_i)_{i\in I})=0 \quad\text{almost surely}.\]
	Taking expectations and then a finite union over all choices of $I$ and $\sigma$ completes the proof.
\end{proof}

\begin{lemma}[{\cite[Proposition~7.2]{drusvyatskiy2014optimality}}]\label{lem:restricted local minimality}
Let $f$ be closed and weakly convex, and let $\cM$ be an active manifold at a critical point $\bar x$. Then $\bar x$ is a (strict) local minimizer of $f|_{\cM}$ if, and only if, it is a (strict) local minimizer of $f$.
\end{lemma}

\begin{lemma}\label{lem:common sign correlation}
Let $\tau>0$ be the unique solution of $\arctan\tau+\tau/(1+\tau^2)=\pi/4$. Let $P\sim N(0,1)$ be independent of jointly standard Gaussian variables $X,Y$ with $\E[XY]=c\in[-1,1]$, and define
\[
q(c)=\E\!\left[XY
\sign(\tau^2X^2-P^2)\sign(\tau^2Y^2-P^2)\right].
\]
Then there exists a constant $\gamma>0$ such that
\[
1-q(c)^2\ge(1-c^2)(1+\gamma c^2)
\qquad\text{for }-1\le c\le1.
\]
\end{lemma}
\begin{proof}
For $p\in\R$, put $f_p(t)=t\sign(\tau^2t^2-p^2)$. This function is odd and has unit norm in $L^2(N(0,1))$. Expanding $f_p$ in the orthonormal Hermite basis, applying the Gaussian correlation identity, and averaging over $P$ yields
\begin{equation}\label{eq:common sign Hermite expansion}
q(c)=\sum_{j\ge0}w_jc^{2j+1}
\qquad\text{with }w_j\ge0,\qquad\text{and}\qquad\sum_{j\ge0}w_j=1.
\end{equation}
Here $w_j$ is the expectation of the squared Hermite coefficient of degree $2j+1$. In particular,
\[
w_0=\E_P\!\left[\left(\E_X[X^2\sign(\tau^2X^2-P^2)]\right)^2\right]<1,
\]
because the inner expectation lies strictly between $-1$ and $1$ whenever $0<|P|<\infty$. Hence
\[
|q(c)|\le w_0|c|+(1-w_0)|c|^3.
\]
Since $[w_0+(1-w_0)c^2]^2\le w_0^2+(1-w_0^2)c^2$, the claim follows with $\gamma=1-w_0^2>0$.
\end{proof}

For the next lemma, let $\tau>0$ and let $X,Z,P$ be independent standard Gaussian variables. For $-1<c<1$, set
\[
s=\sqrt{1-c^2},\qquad Y_c=cX+sZ,\qquad\text{and}\qquad T_c=-sX+cZ,
\]
and define
\[
r(c)=\E\!\left[XZY_cT_c
\sign(\tau^2X^2-P^2)\sign(\tau^2Y_c^2-P^2)\right].
\]
As in the proof of Lemma~\ref{lem:common sign correlation}, write
\[
w_0=\E_P\!\left[\left(\E_X[X^2\sign(\tau^2X^2-P^2)]\right)^2\right],
\qquad\text{and}\qquad \gamma=1-w_0^2.
\]
\begin{lemma}\label{lem:geodesic sign correlation}
Suppose that $2/5\le\tau\le9/20$ and $w_0\ge3/5$. Then $\gamma>0$ and
\[
1-r(c)^2\ge\gamma\sqrt{1-c^2}
\qquad\text{for }-1<c<1,
\]
with equality if, and only if, $c=0$.
\end{lemma}
\begin{proof}
We first obtain global bounds for $r(c)$ from the Hermite expansion, then prove a quantitative endpoint estimate, and finally combine the two bounds.

\emph{Global bounds.}
Use the functions $f_p(x)=x\sign(\tau^2x^2-p^2)$ and Hermite weights $w_j$ from the proof of Lemma~\ref{lem:common sign correlation}, whose construction applies for every $\tau>0$. In particular, $w_j\ge0$, $\sum_{j\ge0}w_j=1$, and $w_0<1$. Write $h_i=\mathrm{He}_i/\sqrt{i!}$ for the orthonormal probabilists' Hermite polynomials. For $i,j\ge1$, the Gaussian moment generating function gives
\begin{equation}\label{eq:geodesic Hermite rotation}
\E[h_i(X)Z\,h_j(Y_c)T_c]
=\mathbf{1}_{\{i=j\}}\bigl(c^{i+1}-i(1-c^2)c^{i-1}\bigr).
\end{equation}
Indeed, for real $a,b$,
\[
\E\!\left[ZT_c e^{aX-a^2/2}e^{bY_c-b^2/2}\right]
=\bigl(c-(1-c^2)ab\bigr)e^{cab},
\]
and comparing coefficients proves~\eqref{eq:geodesic Hermite rotation}. Applying this identity to the Hermite expansion of $f_p$ and averaging over $P$, we obtain, with $u=c^2$,
\begin{equation}\label{eq:geodesic Hermite expansion}
r(c)=\sum_{j\ge0}w_j\bigl((2j+2)u-(2j+1)\bigr)u^j.
\end{equation}
To justify the passage to the infinite expansion, multiplication by the independent $Z$ preserves the $L^2$ norm, as does the orthogonal rotation giving $(Y_c,T_c)$. Thus finite Hermite sums converge in the required expectations. Cauchy--Schwarz bounds each diagonal expectation in~\eqref{eq:geodesic Hermite rotation} by one in absolute value, so $\sum_jw_j=1$ also justifies averaging over $P$.

For $0\le u<1$, each summand polynomial satisfies
\[
\bigl((2j+2)u-(2j+1)\bigr)u^j
=u^j\bigl(2u-1-2j(1-u)\bigr)\le\max\{0,2u-1\}.
\]
For $j\ge1$, its negative is at most $(2j+1)(1-u)u^j\le3u$, since
\[
j(1-u)u^{j-1}\le(1-u)\sum_{i=0}^{j-1}u^i=1-u^j\le1,
\qquad\text{and}\qquad \frac{2j+1}{j}\le3.
\]
Separating the $j=0$ term in~\eqref{eq:geodesic Hermite expansion} gives
\begin{equation}\label{eq:geodesic global bounds}
-w_0+(5w_0-3)c^2\le r(c)\le\max\{0,2c^2-1\}.
\end{equation}
In particular, $r(c)\ge-w_0$ because $w_0\ge3/5$.

\emph{Endpoint estimate.}
We next show that
\begin{equation}\label{eq:geodesic quantitative endpoint}
1-r(c)\ge\frac25\sqrt{1-c^2}
\qquad\text{for }31/32\le c<1.
\end{equation}
Write $c=\cos(2\delta)$ with $0<\delta<\pi/2$, and put $h=\sin\delta$ and $\ell=\cos\delta$. Then $0<h\le1/8$, $\ell^2=1-h^2\ge63/64$, and $s=2h\ell$. Conditioning on $X,Z$ gives
\[
r(c)=\E\!\left[XZY_cT_c\left(1-4\left|
\Phi(\tau|X|)-\Phi(\tau|Y_c|)\right|\right)\right],
\]
where $\Phi$ is the standard Gaussian distribution function. Define the radial integral
\[
\Psi(a)=\int_0^\infty \rho^5e^{-\rho^2/2}
\left(\Phi(\tau\rho a)-\frac12\right)\,\dd\rho
\qquad\text{for }0\le a\le1.
\]
The integrable bound $\tau\rho^6e^{-\rho^2/2}/\sqrt{2\pi}$ justifies differentiation under the integral. The Gaussian sixth moment, together with $1+\tau^2\le(11/10)^2$ and $(11/10)^7<2$, yields
\begin{equation}\label{eq:geodesic radial derivative}
\frac{15\tau}{4}\le\Psi'(a)
=\frac{15\tau}{2(1+\tau^2a^2)^{7/2}}
\le\frac{15\tau}{2}
\qquad\text{for }0\le a\le1.
\end{equation}
Using polar coordinates $X=\rho\cos\phi$ and $Z=\rho\sin\phi$, followed by $\theta=\phi-\delta$, we obtain
\begin{equation}\label{eq:geodesic polar identity}
r(c)=2c^2-1-\frac4\pi I,
\qquad\text{with }I=2\int_0^{\pi/2}K(\theta)D(\theta)\,\dd\theta,
\end{equation}
where
\[
K(\theta)=\sin^2\theta\cos^2\theta-h^2\ell^2,
\quad\text{and}\quad
D(\theta)=\Psi(\cos(\theta-\delta))-\Psi(|\cos(\theta+\delta)|).
\]
Here $\E[XZY_cT_c]=2c^2-1$, the Gaussian density and polar Jacobian leave the radial weight $\rho^5e^{-\rho^2/2}/(2\pi)$, and the original angular integrand is $\pi$-periodic and invariant under $\theta\mapsto\pi-\theta$. These identities give the factor $4/\pi$ in~\eqref{eq:geodesic polar identity}.

For $0\le\theta\le\pi/2$, the two arguments of $\Psi$ differ by $2\min\{h\sin\theta,\ell\cos\theta\}$. Consequently,~\eqref{eq:geodesic radial derivative} gives
\[
\frac{15\tau}{2}\min\{h\sin\theta,\ell\cos\theta\}
\le D(\theta)\le15\tau\min\{h\sin\theta,\ell\cos\theta\}.
\]
Moreover, since $\ell\ge1/2$,
\[
\sin^2\theta\cos^2\theta\min\{h\sin\theta,\ell\cos\theta\}
\ge h\sin^3\theta\cos^2\theta-\frac{h^2}{2}\sin^4\theta\cos\theta.
\]
When the minimum is $h\sin\theta$, this is immediate. Otherwise it follows from
\[
\ell\cos^2\theta-h\sin\theta\cos\theta+\frac{h^2}{2}\sin^2\theta
\ge\frac12(\cos\theta-h\sin\theta)^2\ge0.
\]
Apply the lower bound for $D(\theta)$ to the nonnegative term $\sin^2\theta\cos^2\theta$ in $K(\theta)$ and its upper bound to $-h^2\ell^2$. Using
\[
\int_0^{\pi/2}\sin^3\theta\cos^2\theta\,\dd\theta=\frac2{15},\quad
\int_0^{\pi/2}\sin^4\theta\cos\theta\,\dd\theta=\frac15,\quad\text{and}\quad
\int_0^{\pi/2}\sin\theta\,\dd\theta=1,
\]
we obtain
\[
I\ge15\tau h\left(\frac2{15}-\frac h{10}-2h^2\ell^2\right)
\ge h\left(\frac45-\frac35h-12h^2\ell^2\right)\ge0.
\]
The first bracket is positive because $h\le1/8$ and $\ell^2\le1$, and the second inequality uses $\tau\ge2/5$. Since $\pi<4$,~\eqref{eq:geodesic polar identity} implies
\[
1-r(c)=8h^2\ell^2+\frac4\pi I
\ge\frac45h+h^2\left(8\ell^2-\frac35-12h\ell^2\right)
\ge\frac45h\ge\frac25s.
\]
The last two inequalities use $\ell^2\ge63/64$, $h\ell^2\le1/8$, and $s=2h\ell\le2h$. This proves~\eqref{eq:geodesic quantitative endpoint}.

\emph{The minimum.}
Equation~\eqref{eq:geodesic Hermite expansion} shows that $r(-c)=r(c)$ and $r(0)=-w_0$, so it suffices to consider $0\le c<1$. We have $0<\gamma=1-w_0^2\le16/25$. If $r(c)\le0$,~\eqref{eq:geodesic global bounds} gives $|r(c)|\le w_0$, and hence
\[
\frac{1-r(c)^2}{s}\ge\frac{\gamma}{s}\ge\gamma,
\]
with strict inequality when $c>0$.

Suppose that $r(c)>0$ and $c\le31/32$. By~\eqref{eq:geodesic global bounds}, $u=c^2>1/2$ and $r(c)\le2u-1$, so
\[
\frac{1-r(c)^2}{s}\ge4u\sqrt{1-u}.
\]
For $u\le3/4$, the right-hand side is at least $1$. For $u\ge3/4$, we have $s\ge\sqrt{63}/32\ge7/32$, so it is at least $21/32>16/25\ge\gamma$.

Finally, suppose that $r(c)>0$ and $31/32\le c<1$. The endpoint estimate~\eqref{eq:geodesic quantitative endpoint} gives $r(c)\le1-(2/5)s$, whence
\[
\frac{1-r(c)^2}{s}
\ge\frac{1-(1-(2/5)s)^2}{s}
=\frac45-\frac4{25}s>\frac{16}{25}\ge\gamma.
\]
At $c=0$, the identity $r(0)=-w_0$ gives equality, completing the proof.
\end{proof}

\begin{lemma}[Angular integration]\label{lem:delta integral}
Fix $d\ge4$. Define $g(u,v)=\min\{\norm{u-v},\norm{u+v}\}$ for $u,v\in\mathbb S^{d-2}$.
There exists a constant $0<C_d<\infty$, depending only on $d$, such that for every nonnegative measurable function $f$ on $[0,\infty)$ and every $0\le a<b$,\[ \int_{\substack{u,v\in\mathbb S^{d-2}\\a\le g(u,v)\le b}}f(g(u,v))\dd\sigma(u)\,\dd\sigma(v)\le C_d\int_a^b f(t)t^{d-3}\dd t. \]
Here $\dd\sigma$ denotes the surface measure on $\mathbb S^{d-2}$.
\end{lemma}
\begin{proof}
Fix $u$. Reflecting $v$ through the origin produces a factor of two and reduces the angular integral to $0\le\theta\le\pi/2$, where $u^\top v=\cos\theta$. Write $v=\cos\theta\,u+\sin\theta\,w$, where $w\in\mathbb S^{d-2}\cap u^\perp$, a sphere identified with $\mathbb S^{d-3}$. On the hemisphere $\{v:u^\top v\ge0\}$, \[ g(u,v)=\norm{u-v}=\sqrt{2-2\cos\theta} = 2\sin(\theta/2),\qquad\text{and}\qquad \dd\sigma(v)=(\sin\theta)^{d-3}\dd\theta\,\dd\sigma(w). \]
Integrating over $v$ implies
\begin{align*}
\int_{\substack{v\in\mathbb S^{d-2}\\a\le g(u,v)\le b}}f(g(u,v))\dd\sigma(v)
&=2\int_{\substack{0\le\theta\le\pi/2\\a\le2\sin(\theta/2)\le b}}\int_{\mathbb S^{d-3}}f(2\sin(\theta/2))(\sin\theta)^{d-3}\dd\sigma(w)\,\dd\theta\\
&=2\abs{\mathbb S^{d-3}}\int_{\substack{0\le\theta\le\pi/2\\a\le2\sin(\theta/2)\le b}}f(2\sin(\theta/2))(\sin\theta)^{d-3}\dd\theta\\
&=2\abs{\mathbb S^{d-3}}\int_{[a,b]\cap[0,\sqrt2]}f(t)t^{d-3}(1-t^2/4)^{(d-4)/2}\dd t,
\end{align*}
where the last equality uses the change of variables $t=2\sin(\theta/2)$.
Integrating over $u$ as well yields \[ \int_{\substack{u,v\in\mathbb S^{d-2}\\a\le g(u,v)\le b}}f(g(u,v))\dd\sigma(u)\,\dd\sigma(v)=2\abs{\mathbb S^{d-2}}\abs{\mathbb S^{d-3}}\int_{[a,b]\cap[0,\sqrt2]}f(t)t^{d-3}(1-t^2/4)^{(d-4)/2}\dd t. \]
Since $d\ge4$, the last factor is at most one. The bound therefore holds with $C_d=2\abs{\mathbb S^{d-2}}\abs{\mathbb S^{d-3}}$.
\end{proof}

\begin{lemma}[Ess\'een concentration inequality {\cite[Exercise~2.2.11, p.~85]{tao2012topics}}]\label{lem:Esseen concentration inequality}
For every integer $m\ge1$ and every $\varepsilon>0$, there exists a constant $C_{m,\varepsilon}>0$, depending only on $m$ and $\varepsilon$, such that for every random vector $X\in\R^m$ and every $R>0$,
\[
\sup_{z\in\R^m}\PP\bigl(X\in z+B(0,R)\bigr)
\le C_{m,\varepsilon}R^m
\int_{\norm{v}\le\varepsilon/R}\abs{\E e^{iv^\top X}}\dd v.
\]
\end{lemma}

%%% Local Variables:
%%% mode: LaTeX
%%% TeX-master: "../main"
%%% End:

\section{Numerical experiment details}\label{app:experiments}
This appendix details how we test the criticality condition~\eqref{eq:criticality} numerically for the candidates in the first experiment of Section~\ref{sec:experiments}. The signal is fixed at $\xs=e_1$, and the dimension, measurement ratios, and trial counts are specified in the main text.

\paragraph{Candidate construction.}
For each enumerated index set $I=\{i_1<\cdots<i_d\}$ and nonconstant sign vector $s\in\{\pm1\}^d$, we impose the active equations $a_{i_j}^\top\widehat x=s_j a_{i_j}^\top e_1$ for $j=1,\ldots,d$ and form the candidate $\widehat x=x_{I,s}$ from~\eqref{eq:candidate}. Let $S_I(\widehat x)$ and $Q_I(\widehat x)$ be as defined in Section~\ref{sec:proof-strategy}, with the columns of $Q_I$ ordered by $i_1,\ldots,i_d$.

\paragraph{Numerical verification.}
We recompute each candidate selected for verification and its multiplier $\lambda\in\RR^d$ using 60-digit decimal arithmetic and then apply the following checks.
\begin{enumerate}[leftmargin=*]
\item \emph{Consistency of the active set.} We solve the selected active equations and require every unselected measurement to have a nonzero residual:
\[
\min_{i\notin I}|(a_i^\top\widehat x)^2-(a_i^\top e_1)^2|>10^{-30}.
\]
This check numerically separates the prescribed active measurements from the remaining ones.

\item \emph{Criticality with an interior multiplier.} We solve $Q_I(\widehat x)\lambda=-S_I(\widehat x)$ and require
\[
1-\|\lambda\|_\infty>10^{-30}.
\]
This condition keeps the multiplier strictly inside $(-1,1)^d$, as required by~\eqref{eq:strategy-criticality}.
\end{enumerate}

We count a problem instance as positive if at least one candidate passes both checks.

%%% Local Variables:
%%% mode: LaTeX
%%% TeX-master: "../main"
%%% End:

\end{document}